\documentclass[12pt]{article}
\usepackage{amssymb}
\usepackage{amsmath}
\usepackage{amsthm}
\usepackage{latexsym}
\usepackage{amscd}

\newtheorem{definition}{Definition}[section]
\newtheorem{theorem}[definition]{Theorem}
\newtheorem{lemma}[definition]{Lemma}
\newtheorem{corollary}[definition]{Corollary}

\newtheorem{example}[definition]{Example}

\newtheorem{note}[definition]{Note}

\newtheorem{proposition}[definition]{Proposition}

\def\F{\mathbb F}

\def\R{\mathbb R}
\def\C{\mathbb C}
\def\K{\mathbb K}

\def\K{\mathbb K}
\def\fld{\mathbb K}

\begin{document}
\title{Dual polar graphs, the quantum algebra $U_q(\mathfrak{sl}_2)$, and Leonard systems of dual $q$-Krawtchouk type}
\author{Chalermpong Worawannotai}
\date{}
\maketitle

\begin{abstract}
In this paper we consider how the following three objects are related:
(i) the dual polar graphs; (ii) the quantum algebra $U_q(\mathfrak{sl}_2)$; (iii) the Leonard systems of dual $q$-Krawtchouk type.
For convenience we first describe how (ii) and (iii) are related.
For a given Leonard system of dual $q$-Krawtchouk type,
we obtain two $U_q(\mathfrak{sl}_2)$-module structures on its underlying vector space.
We now describe how (i) and (iii) are related.
Let $\Gamma$ denote a dual polar graph.
Fix a vertex $x$ of $\Gamma$
and let $T = T(x)$ denote the corresponding subconstituent algebra.
By definition $T$ is generated by the adjacency matrix $A$ of $\Gamma$ and a certain diagonal matrix $A^* = A^*(x)$ called the dual adjacency matrix that corresponds to $x$.
By construction the algebra $T$ is semisimple.
We show that for each irreducible $T$-module $W$ the restrictions of $A$ and $A^*$ to $W$ induce a Leonard system of dual $q$-Krawtchouk type.
We now describe how (i) and (ii) are related.
We obtain two $U_q(\mathfrak{sl}_2)$-module structures on the standard module of $\Gamma$.
We describe how these two $U_q(\mathfrak{sl}_2)$-module structures are related.
Each of these $U_q(\mathfrak{sl}_2)$-module structures induces a $\mathbb{C}$-algebra homomorphism $U_q(\mathfrak{sl}_2) \rightarrow T$.
We show that in each case $T$ is generated by the image together with the center of $T$.
Using the combinatorics of $\Gamma$ we obtain a generating set $L, F, R, K$ of $T$ along with some attractive relations satisfied by these generators.

\end{abstract}

\section{Introduction}

In this paper we investigate a topic that involves algebraic graph theory, quantum groups, and linear algebra.
We focus on three objects that turn out to be closely related.
The first object is called a {\it dual polar graph} \cite{bannai,bcn,cameron}.
Dual polar graphs are distance-regular \cite{bcn}.
The second object is the quantized universal enveloping algebra $U_q(\mathfrak{sl}_2)$.
We refer the reader to \cite{jantzen,kassel} for background information on $U_q(\mathfrak{sl}_2)$.
The third object is called a {\it Leonard system of dual $q$-Krawtchouk type}.
In order to explain this concept we start with a more basic notion called a {\it Leonard pair} \cite{LS99}.
Roughly speaking a Leonard pair consists of two diagonalizable linear transformations on a finite-dimensional vector space, each of which acts in an irreducible tridiagonal fashion on an eigenbasis for the other one.
A Leonard system is an ``oriented'' version of a Leonard pair.
The Leonard systems are classified up to isomorphism \cite{LS99,madrid09}.
We will focus on a family of Leonard systems said to have {\it dual $q$-Krawtchouk type}.

Our central results are about how the above three objects are related.
Shortly we will summarize these results.
First we describe the three objects more precisely.
We begin with the definition of a Leonard system.
For convenience we take the underlying field to be the complex number field $\C$.

\begin{definition} \label{def:intro-ls} \rm \cite{LS99}
Let $d$ denote a nonnegative integer
and let $V$ denote a vector space over $\C$ with dimension $d+1$.
By a {\em Leonard system} on $V$ we mean a sequence
\begin{eqnarray*}
\Phi = (A;\{E_i\}_{i=0}^d;A^*;\{E^*_i\}_{i=0}^d)
\end{eqnarray*}
that satisfies (i)--(v) below.
\begin{enumerate}
\item[\rm (i)] Each of $A,A^*$ is a multiplicity-free element in ${\rm End}(V)$.
\item[\rm (ii)] $\{E_i\}_{i=0}^d$ is an ordering of the primitive idempotents of $A$.
\item[\rm (iii)] $\{E^*_i\}_{i=0}^d$ is an ordering of the primitive idempotents of $A^*$.
\item[\rm (iv)] $E_iA^*E_j =
\begin{cases}
0 &\mbox{if } |i-j| > 1\\
\neq 0 &\mbox{if } |i-j| = 1
\end{cases}
\qquad (0 \leq i, j \leq d).$
\item[\rm (v)] $E^*_iAE^*_j =
\begin{cases}
0 &\mbox{if } |i-j| > 1\\
\neq 0 &\mbox{if } |i-j| = 1
\end{cases}
\qquad (0 \leq i, j \leq d).$
\end{enumerate}
\end{definition}

\begin{definition}\rm
Referring to Definition \ref{def:intro-ls}, for $0 \leq i \leq d$ let $\theta_i$ (resp. $\theta^*_i$) denote the eigenvalue of $A$ (resp. $A^*$) associated with $E_i$ (resp. $E^*_i$).
\end{definition}

Referring to the Leonard system $\Phi$ from Definition \ref{def:intro-ls}, assume $\Phi$ has dual $q$-Krawtchouk type.
As we will see, the eigenvalues of $A$ and $A^*$ have the form
\begin{align}
\theta_i= h + \kappa q^{d-2i} + \upsilon q^{2i-d},\qquad
\theta_i^*= h^* + \kappa^*q^{d-2i}\qquad (0 \leq i \leq d)\label{eq:intro-ls-eval}
\end{align}
where $q,h,h^*,\kappa,\kappa^*,\upsilon$ are scalars in $\C$ with $q^2 \neq 1$ and $\kappa,\kappa^*,\upsilon$ nonzero.

We now recall the algebra $U_q(\mathfrak{sl}_2)$.
We will use the equitable presentation \cite{equit1,uaw-uqsl2}.

\begin{definition} {\rm \cite{equit1}}
Let $U_q(\mathfrak{sl}_2)$ denote the $\C$-algebra with generators $x,y,z^{\pm 1}$ and relations $zz^{-1} = z^{-1}z = 1,$
\begin{eqnarray*}
\frac{qxy-q^{-1}yx}{q-q^{-1}} = 1,\qquad
\frac{qyz-q^{-1}zy}{q-q^{-1}} = 1,\qquad
\frac{qzx-q^{-1}xz}{q-q^{-1}} = 1.
\end{eqnarray*}
\end{definition}

We now show how a Leonard system of dual $q$-Krawtchouk type is related to $U_q(\mathfrak{sl}_2)$.
The following two theorems are our main results along this line.

\begin{theorem}\label{thm:intro-ls-u1}
Fix $\epsilon \in \{1,-1\}$.
Let $\Phi$ denote a Leonard system on $V$ as in Definition \ref{def:intro-ls}.
Assume $\Phi$ has dual $q$-Krawtchouk type and let $h,h^*,\kappa,\kappa^*,\upsilon$ denote the corresponding parameters from {\rm (\ref{eq:intro-ls-eval})}.
Then there exists a unique $U_q(\mathfrak{sl}_2)$-module structure on $V$ such that on $V$,
\begin{align*}
A & = h1 + \epsilon \kappa x + \epsilon \upsilon y,\\
A^* & = h^*1 + \epsilon \kappa^*z.
\end{align*}
\end{theorem}

\begin{theorem}\label{thm:intro-ls-u2}
Fix $\epsilon \in \{1,-1\}$.
Let $\Phi$ denote a Leonard system on $V$ as in Definition \ref{def:intro-ls}.
Assume $\Phi$ has dual $q$-Krawtchouk type and let $h,h^*,\kappa,\kappa^*,\upsilon$ denote the corresponding parameters from {\rm (\ref{eq:intro-ls-eval})}.
Then there exists a unique $U_q(\mathfrak{sl}_2)$-module structure on $V$ such that on $V$,
\begin{align*}
A & = h1 + \epsilon \kappa y + \epsilon \upsilon x,\\
A^* & = h^*1 + \epsilon \kappa^*z.
\end{align*}
\end{theorem}

We now recall the definition of a dual polar graph.
Let $b$ denote a prime power.
Let $\F_b$ denote the finite field of order $b$.
Let $U$ denote  a finite-dimensional vector space over $\F_b$ endowed with one of the following forms: $C_D(b),$ $B_D(b),$ $D_D(b),$ $^2D_{D+1}(b),$ $^2A_{2D}(q)$, $^2A_{2D-1}(q),$ where $b = q^2$ \cite[p.~274]{bcn}.
A subspace $W$ of $U$ is called {\it isotropic} whenever the form vanishes completely on $W$.  By \cite[Theorem~6.3.1]{cameron} each maximal isotropic subspace of $U$ has dimension $D$.
Define a graph $\Gamma$ as follows.
The vertex set $X$ of $\Gamma$ consists of the maximal isotropic subspaces of $U$.
Vertices $y, z$ in $X$ are adjacent in $\Gamma$ whenever ${\rm dim}(y \cap z) = D-1$.
Let $\partial$ denote the path-length distance function for $\Gamma$.
By \cite[p.~276]{bcn} for $y,z \in X$, $\partial(y,z) = D-{\rm dim}(y \cap z)$.
By \cite[p.~274]{bcn} the graph $\Gamma$ is distance-regular with diameter $D$.
We call $\Gamma$ the {\it dual polar graph} associated with $U$.

We have a few comments about $\Gamma$.
By the {\it standard module} of $\Gamma$ we mean the vector space $V = \C^X$ of column vectors with rows indexed by $X$.
Let $A \in {\rm Mat}_X(\C)$ denote the adjacency matrix of $\Gamma$.
Let $\theta_0 > \theta_1 > \cdots > \theta_D$ denote the eigenvalues of $A$.
For $0 \leq i \leq D$ let $E_i$ denote the projection onto the eigenspace of $A$ associated with the eigenvalue $\theta_i$.
For the rest of this section fix a vertex $x \in X$.
Let $A^* = A^*(x)$ denote the diagonal matrix in ${\rm Mat}_X(\C)$ whose diagonal is obtained by rotating row $x$ of $|X|E_1$ by 45 degrees.
For $0 \leq i \leq D$, by the {\it $i^{th}$ subconstituent} of $\Gamma$ we mean the subspace of $V$  spanned by the vertices at distance $i$ from $x$.
The subconstituents are the eigenspaces of $A^*$;
for $0 \leq i \leq D$ let $E^*_i = E^*_i(x)$ (resp. $\theta^*_i$) denote the corresponding projection (resp. eigenvalue).
By \cite[Theorem~8.4.2, Theorem~9.4.3]{bcn} the eigenvalues of $A$ and $A^*$ have the form
 \begin{align*}
\theta_i= h + \kappa q^{D-2i} + \upsilon q^{2i-D},\qquad
\theta_i^*= h^* + \kappa^*q^{D-2i}\qquad (0 \leq i \leq D)
\end{align*}
where $h,h^*,\kappa,\kappa^*,\upsilon$ are in $\C$ with $\kappa,\kappa^*,\upsilon$ nonzero.
Let $T = T(x)$ denote the subalgebra of ${\rm Mat}_X(\C)$ generated by $A, A^*$.
The algebra $T$ is called the {\it subconstituent algebra} or {\it Terwilliger algebra} with respect to $x$ \cite{terwSub1}.
By \cite[p.~157]{CR} the algebra $T$ is semisimple.
By a $T$-module we mean a subspace $W \subseteq V$ such that $BW \subseteq W$ for all $B \in T$.
Let $W$ denote an irreducible $T$-module.
By the {\it endpoint} of $W$ we mean
$\mbox{min}\lbrace i |0\leq i \leq D, \; E^*_iW\not=0\rbrace $.
By the {\it dual endpoint} of $W$ we mean
$\mbox{min}\lbrace i |0\leq i \leq D, \; E_iW\not=0\rbrace $.
By the {\it diameter} of $W$ we mean
$ |\lbrace i | 0 \leq i \leq D,\; E^*_iW\not=0 \rbrace |-1 $.
We now show how $\Gamma$ is related to the Leonard systems of dual $q$-Krawtchouk type.
The following is our third main result.

\begin{theorem}\label{thm:intro-ls-on-w}
Let $W$ denote an irreducible $T$-module.
Let $r, t, d$ denote the endpoint, dual endpoint, and diameter of $W$, respectively.
Then $(A|_W;\{E_{t+i}|_W\}_{i=0}^d;A^*|_W;\{E^*_{r+i}|_W\}_{i=0}^d)$
is a Leonard system of dual $q$-Krawtchouk type.
\end{theorem}

We now summarize how $\Gamma$ is related to $U_q(\mathfrak{sl}_2)$.
Since $T$ is semisimple, $V$ is a direct sum of irreducible $T$-modules.
For each irreducible $T$-module $W$ in the sum, combining Theorem \ref{thm:intro-ls-on-w} with Theorem \ref{thm:intro-ls-u1} and Theorem \ref{thm:intro-ls-u2} we obtain two $U_q(\mathfrak{sl}_2)$-module structures on $W$.
This gives two $U_q(\mathfrak{sl}_2)$-module structures on $V$.
In order to describe them in a coherent fashion we introduce some elements $\Upsilon, \Psi$ in the center of $T$.
These elements act on each irreducible $T$-module $W$ as $q^{r+t+d-D}I$, $q^{r-t}I$ where $r, t, d$ are the endpoint, dual endpoint, and diameter of $W$, respectively.
We now give our fourth and fifth main results.

\begin{theorem}\label{thm:intro-u1}
There exists a unique $U_q(\mathfrak{sl}_2)$-module structure on $V$ such that on $V$,
\begin{align*}
A & = h1 + \kappa\Upsilon^{-1}\Psi x + \upsilon\Upsilon\Psi^{-1} y,\\
A^* & = h^*1 + \kappa^*\Upsilon^{-1}\Psi^{-1} z.
\end{align*}
\end{theorem}

\begin{theorem}\label{thm:intro-u2}
There exists a unique $U_q(\mathfrak{sl}_2)$-module structure on $V$ such that on $V$,
\begin{align*}
A & = h1 + \kappa\Upsilon^{-1}\Psi y + \upsilon\Upsilon\Psi^{-1} x,\\
A^* & = h^*1 + \kappa^*\Upsilon^{-1}\Psi^{-1} z.
\end{align*}
\end{theorem}

Each of the above $U_q(\mathfrak{sl}_2)$-module structures on $V$ induces a $\C$-algebra homomorphism $U_q(\mathfrak{sl}_2) \rightarrow T$.
We now describe how their images are related to $T$.
The following is our sixth main result.

\begin{theorem}\label{thm:intro-u-gen-t}
For either of our two $\C$-algebra homomorphisms $U_q(\mathfrak{sl}_2) \rightarrow T$, let $U$ denote the image.
Then the algebra $T$ is generated by $U$ together with the elements $\Upsilon^{\pm 1}, \Psi^{\pm 1}$.
\end{theorem}

We now describe some relations in $T$ that we find attractive.
In order to state these relations, it is convenient to decompose
$A = L + F + R$ where
$L$ (resp. $F$) (resp. $R$) is the lowering matrix (resp. flattening matrix) (resp. raising matrix) of $\Gamma$ with respect to $x$.
For vertices $y, z$ of $\Gamma$ the $(y,z)$-entry of $L$ (resp. $F$) (resp. $R$) is 1 whenever $y, z$ are adjacent and $\partial(x,z) - \partial(x,y)$ is 1 (resp. 0) (resp. $-1$).  
Define a diagonal matrix $K \in {\rm Mat}_X(\C)$ with $(y,y)$-entry $q^{-2\partial(x,y)}$ for $y \in X$.
In other words $K = \sum_{i=0}^D q^{-2i}E^*_i$. 
By construction $A^* \in {\rm span}\{I,K\}$.
The algebra $T$ is generated by $L, F, R, K$.
We now describe how $L, F, R, K$ are related.
The following is our seventh main result.

\begin{theorem}\label{thm:intro-lfrk}
The matrices $L, F, R, K$ satisfy
\begin{eqnarray*}
&KL = q^2 LK, \qquad KF = FK, \qquad KR = q^{-2}RK,\\
&LF - q^2 FL = (q^{2e}-1) L, \qquad FR - q^2 RF = (q^{2e}-1) R,\\
&\displaystyle\frac{q^4}{q^2+1} RL^2 - LRL + \frac{q^{-2}}{q^2+1} L^2 R = - q^{2e+2D-2} L,\\
&\displaystyle\frac{q^4}{q^2+1} R^2 L - RLR + \frac{q^{-2}}{q^2+1} L R^2 = - q^{2e+2D-2} R,
\end{eqnarray*}
where $e$ is given in the table below:
\begin{center}
\begin{tabular}[b]{c|cccccc}
{\rm form} & $C_D(b)$ & $B_D(b)$ & $D_D(b)$ & $^2D_{D+1}(b)$ & $^2A_{2D}(q)$ & $^2A_{2D-1}(q)$\\ \hline
e & $1$ & $1$ & $0$ & $2$ & $3/2$ & $1/2$\\
\end{tabular}
\end{center}
\end{theorem}

We will repeatedly use the relations in Theorem \ref{thm:intro-lfrk}.

This paper is organized as follows.
In Sections 2--7 we recall some background concerning Leonard systems. 
In Section 8 we introduce the normalized split basis for a Leonard system.
In Section 9 we discuss the intersection matrix of a Leonard system.
In Section 10 we recall the tridiagonal relations and Askey-Wilson relations of a Leonard system.
In Section 11 we recall the Leonard systems of dual $q$-Krawtchouk type.
In Section 12 we recall $U_q(\mathfrak{sl}_2)$ and describe its finite-dimensional irreducible modules.
In Section 13 we prove Theorem \ref{thm:intro-ls-u1} and Theorem \ref{thm:intro-ls-u2}.
In Section 14 we discuss the subconstituent algebra $T$ of a distance-regular graph.
In Section 15 we recall the notion of a near polygon.
After obtaining some basic facts about near polygons, we focus on a particular type of near polygon called a dual polar graph.
In Sections 16--19 we discuss some basic facts about a dual polar graph and its irreducible $T$-modules.
In Sections 20, 21 we discuss some central elements $\Omega, G, G^*$ of $T$ that come from the Askey-Wilson relations.
We describe the entries of the matrices $\Omega, G, G^*$.
In Section 22 we introduce three central elements $\Upsilon, \Psi, \Lambda$ of $T$ which will be used to relate $T$ to $U_q(\mathfrak{sl}_2)$.
In Section 23 we prove Theorem \ref{thm:intro-ls-on-w}.
In Section 24 we prove Theorem \ref{thm:intro-u1} and Theorem \ref{thm:intro-u2}.
In Section 25 we prove Theorem \ref{thm:intro-u-gen-t}.
In Section 26 we discuss the matrices $L, F, R, K$ and prove Theorem \ref{thm:intro-lfrk}.
In Section 27 we describe $\Omega, G, G^*$ in terms of $L, F, R, K$.
In Section 28 we introduce three central elements $C_0, C_1, C_2$ of $T$ that involve $L, F, R, K$.
We show that $C_0, C_1, C_2$ generate the center of $T$.
In Section 29 we show how $C_0, C_1, C_2$ relate to $\Omega, G, G^*$ and $\Upsilon, \Psi, \Lambda$.
In Section 30 we describe the two $U_q(\mathfrak{sl}_2)$-module structures from Theorem \ref{thm:intro-u1} and Theorem \ref{thm:intro-u2} in terms of $L, F, R, K$.

\section{Leonard pairs}

We now begin our formal argument.
We start by recalling the notion of a Leonard pair.  We will use the following terms.  A square matrix $X$ is said to be {\em tridiagonal} whenever each nonzero entry lies on either the diagonal, the subdiagonal, or the superdiagonal.  Assume $X$ is tridiagonal.  Then $X$ is said to be {\em irreducible} whenever each entry on the subdiagonal is nonzero and each entry on the superdiagonal is nonzero.  We now define a Leonard pair.  For the rest of this paper $\K$ will denote a field.

\begin{definition} \label{def:lp} \rm \cite[Definition~1.1]{LS99}
Let $V$ denote a vector space over $\K$ with finite positive dimension.  By a {\em Leonard pair} on $V$ we mean an ordered pair $A,A^*$ where $A : V \rightarrow V$ and $A^* : V \rightarrow V$ are linear transformations that satisfy (i), (ii) below:
\begin{enumerate}
\item[\rm (i)] There exists a basis for $V$ with respect to which the matrix representing $A$ is irreducible tridiagonal and the matrix representing $A^*$ is diagonal.
\item[\rm (ii)] There exists a basis for $V$ with respect to which the matrix representing $A^*$ is irreducible tridiagonal and the matrix representing $A$ is diagonal.
\end{enumerate}
\end{definition}

\begin{note} \rm
It is a common notational convention to use $A^*$ to represent the conjugate-transpose of $A$.  We are not using this convention.  In a Leonard pair $A,A^*$ the linear transformations $A$ and $A^*$ are arbitrary subject to (i), (ii) above.
\end{note}

\section{Leonard systems}
When working with a Leonard pair, it is convenient to consider a closely related object called a Leonard system.
To prepare for our definition of a Leonard system, we recall a few concepts from linear algebra.  
Throughout the paper an algebra is meant to be associative and have a $1$,
and a subalgebra has the same $1$ as the parent algebra.
Let $d$ denote a nonnegative integer and let ${\rm Mat}_{d+1}(\K)$ denote the $\K$-algebra consisting of all $d+1$ by $d+1$ matrices that have entries in $\K$.  We index the rows and columns by $0,1,\dots ,d$.
Let $V$ denote a vector space over $\K$ with dimension $d+1$.
Let $\{v_i\}_{i=0}^d$ denote a basis for $V$.  For $A \in {\rm End}(V)$ and $B \in {\rm Mat}_{d+1}(\K)$, we say that $B$ {\em represents $A$ with respect to $\{v_i\}_{i=0}^d$} whenever $Av_j = \sum_{i=0}^d B_{ij}v_i$ for $0 \leq j \leq d$.
An element $A \in {\rm End}(V)$ is said to be {\em multiplicity-free} whenever it has $d+1$ mutually distinct eigenvalues in $\K$.
Assume $A$ is multiplicity-free.
Let $\{V_i\}_{i=0}^d$ denote an ordering of the eigenspaces of $A$.
For $0 \leq i \leq d$ let $\theta_i$ denote the eigenvalue of $A$ corresponding to $V_i$.
Define $E_i \in {\rm End}(V)$ such that $(E_i - I)V_i = 0$ and $E_i V_j = 0$ for $j \neq i$ $(0 \leq j \leq d)$.
Here $I$ denotes the identity of ${\rm End}(V)$.
We call $E_i$ the {\em primitive idempotent} of $A$ corresponding to $V_i$ (or $\theta_i$).
Observe that
(i) $AE_i = \theta_i E_i \; (0 \leq i \leq d)$;
(ii) $E_i E_j = \delta_{ij} E_i \; (0 \leq i,j \leq d);$
(iii) I = $\sum_{i=0}^d E_i$;
(iv) $A = \sum_{i=0}^d \theta_i E_i$.
Moreover
\begin{eqnarray*}
E_i = \prod_{\genfrac{}{}{0pt}{}{0 \leq j \leq d}{j \neq i}} \frac{A-\theta_j I}{\theta_i - \theta_j}.
\end{eqnarray*}

\medskip \noindent
We now define a Leonard system.

\begin{definition} \label{def:ls} \rm \cite[Definition~1.4]{LS99}
Let $d$ denote a nonnegative integer
and let $V$ denote a vector space over $\K$ with dimension $d+1$.
By a {\em Leonard system} on $V$ we mean a sequence
\begin{eqnarray*}
\Phi = (A;\{E_i\}_{i=0}^d;A^*;\{E^*_i\}_{i=0}^d)
\end{eqnarray*}
that satisfies (i)--(v) below.
\begin{enumerate}
\item[\rm (i)] Each of $A,A^*$ is a multiplicity-free element in ${\rm End}(V)$.
\item[\rm (ii)] $\{E_i\}_{i=0}^d$ is an ordering of the primitive idempotents of $A$.
\item[\rm (iii)] $\{E^*_i\}_{i=0}^d$ is an ordering of the primitive idempotents of $A^*$.
\item[\rm (iv)] $E_iA^*E_j =
\begin{cases}
0 &\mbox{if } |i-j| > 1\\
\neq 0 &\mbox{if } |i-j| = 1
\end{cases}
\qquad (0 \leq i, j \leq d).$
\item[\rm (v)] $E^*_iAE^*_j =
\begin{cases}
0 &\mbox{if } |i-j| > 1\\
\neq 0 &\mbox{if } |i-j| = 1
\end{cases}
\qquad (0 \leq i, j \leq d).$
\end{enumerate}
We refer to $d$ as the {\em diameter} of $\Phi$, and say $\Phi$ is {\em over} $\K$.
We call $V$ the {\em vector space underlying $\Phi$}.
\end{definition}

\noindent We comment on how 
Leonard pairs and Leonard systems are related.
Fix an integer $d \geq 0$ and let $V$ denote
a vector space over $\K$ with dimension $d+1$. Let 
$(A; \lbrace E_i\rbrace_{i=0}^d; A^*; \lbrace E^*_i\rbrace_{i=0}^d) $
denote a Leonard system on $V$. 
For $0 \leq i \leq d$ let $v_i$ denote a nonzero vector in 
$E_iV$. Then the sequence 
$\{v_i\}_{i=0}^d$ is a basis for $V$ that satisfies
Definition 
\ref{def:lp}(ii).
For $0 \leq i \leq d$ let $v^*_i$ denote a nonzero vector in 
$E^*_iV$. Then the sequence 
$\{v^*_i\}_{i=0}^d$ is a basis for $V$ that satisfies
Definition 
\ref{def:lp}(i). By these comments the pair $A,A^*$ is
a Leonard pair on $V$. Conversely let $A,A^*$ denote
a Leonard pair on $V$.
By \cite[Lemma~3.1]{madrid09} each of $A, A^*$ is multiplicity-free.
 Let 
$\{v_i\}_{i=0}^d$ denote a basis for $V$ that satisfies
Definition 
\ref{def:lp}(ii). For $0 \leq i \leq d$ the vector
$v_i$ is an eigenvector for $A$;
 let
$E_i$ denote the corresponding primitive idempotent.
 Let 
$\{v^*_i\}_{i=0}^d$ denote a basis for $V$ that satisfies
Definition 
\ref{def:lp}(i). For $0 \leq i \leq d$ the vector
$v^*_i$ is an eigenvector for $A^*$;
 let
$E^*_i$ denote the corresponding primitive idempotent.
Then 
$(A; \lbrace E_i\rbrace_{i=0}^d;  A^*; \lbrace E^*_i\rbrace_{i=0}^d)$
is a Leonard system on $V$.

\begin{definition} \label{eigenvalue-seqn} \rm
Referring to the Leonard system $\Phi$ from Definition \ref{def:ls},
for $0 \leq i \leq d$ let $\theta_i$ (resp. $\theta^*_i$) denote the eigenvalue of $A$ (resp. $A^*$) associated with the eigenspace $E_iV$ (resp. $E^*_iV$).
We call $\{\theta_i\}_{i=0}^d$ (resp. $\{\theta^*_i\}_{i=0}^d$) the {\em eigenvalue sequence} (resp. {\em dual eigenvalue sequence}) of $\Phi$.
\end{definition}

\medskip
\noindent The following notation will be useful.  Let $\lambda$ denote an indeterminate and let $\K[\lambda]$ denote the $\K$-algebra consisting of the polynomials in $\lambda$ that have all coefficients in $\K$.

\begin{definition} \label{def:tau-eta} \rm
Referring to the Leonard system $\Phi$ from Definition \ref{def:ls},
let $\{\theta_i\}_{i=0}^d$ (resp. $\{\theta^*_i\}_{i=0}^d$) denote the eigenvalue sequence (resp. dual eigenvalue sequence) of $\Phi$.  For $0 \leq i \leq d$ define polynomials $\tau_i, \eta_i, \tau^*_i, \eta^*_i$ in $\K[\lambda]$ as follows.
\begin{align*}
\tau_i &= (\lambda - \theta_0)  (\lambda - \theta_1) \cdots  (\lambda - \theta_{i-1}),\\
\eta_i &= (\lambda - \theta_d)  (\lambda - \theta_{d-1}) \cdots  (\lambda - \theta_{d-i+1}),\\
\tau^*_i &= (\lambda - \theta^*_0)  (\lambda - \theta^*_1) \cdots  (\lambda - \theta^*_{i-1}),\\
\eta^*_i &= (\lambda - \theta^*_d)  (\lambda - \theta^*_{d-1}) \cdots  (\lambda - \theta^*_{d-i+1}).
\end{align*}
Observe that each of $\tau_i, \eta_i, \tau^*_i, \eta^*_i$ is monic of degree $i$.
\end{definition}

\section{The $D_4$ action}

\noindent Let $\Phi$ denote the Leonard system on $V$ from Definition \ref{def:ls}.
Then each of the following three sequences is a Leonard system
on $V$.
\begin{align*}
 \;\Phi^*&:= (A^*; \lbrace E^*_i\rbrace_{i=0}^d; 
 A; \lbrace E_i\rbrace_{i=0}^d),
\\
\Phi^{\downarrow}&:= (A; \lbrace E_i\rbrace_{i=0}^d ;
A^*; \lbrace E^*_{d-i} \rbrace_{i=0}^d),
\\
\Phi^{\Downarrow} 
&:= (A; \lbrace E_{d-i}\rbrace_{i=0}^d;
A^*; \lbrace E^*_i\rbrace_{i=0}^d).
\end{align*}
Viewing $*, \downarrow, \Downarrow$
as permutations on the set of all Leonard systems,
\begin{eqnarray}
&&\qquad \qquad \qquad *^2 = \; \downarrow^2 \; = \; \Downarrow^2 \; = 1,\label{eq:deightrelationsAS99}\\
&&\Downarrow * = * \downarrow,
\qquad \qquad
\downarrow * = * \Downarrow,
\qquad \qquad
\downarrow \Downarrow \; = \; \Downarrow \downarrow.
\label{eq:deightrelationsBS99}
\end{eqnarray}
The group generated by symbols 
$*, \downarrow, \Downarrow $ subject to the relations
(\ref{eq:deightrelationsAS99}),
(\ref{eq:deightrelationsBS99})
is the dihedral group $D_4$.  
We recall $D_4$ is the group of symmetries of a square,
and has 8 elements.
Apparently $*, \downarrow, \Downarrow $ induce an action of 
 $D_4$ on the set of all Leonard systems.
Two Leonard systems will be called {\it relatives} whenever they
are in the same orbit of this $D_4$ action.

\medskip
\noindent For the rest of this paper we will use the following convention.

\begin{definition}
\label{def:notconv}
\rm
Referring to Leonard system $\Phi$ from Definition \ref{def:ls},
for any element $g$ in
the group $D_4$ and for any object $f$ associated with
$\Phi$, let $f^g$ denote the corresponding object for the
Leonard system 
$\Phi^{g^{-1}}$.
\end{definition}

\section{The standard decomposition and the standard basis}

\noindent Throughout this section fix an integer $d \geq 0$ and let $V$ denote a vector space over $\K$ with dimension $d+1$.
By a {\em decomposition} of $V$ we mean a sequence $\{V_i\}_{i=0}^d$ of subspaces of $V$ such that $V_i$ has dimension $1$ for $0 \leq i \leq d$ and $V = \sum_{i=0}^d V_i$ (direct sum).
Let $\{V_i\}_{i=0}^d$ denote a decomposition of $V$.  By the {\em inversion} of this decomposition we mean the decomposition $\{V_{d-i}\}_{i=0}^d$.

\begin{definition} \label{def:stand-decomp} \rm
Referring to the Leonard system $\Phi$ on $V$ from Definition \ref{def:ls},
observe that $\{E^*_iV\}_{i=0}^d$ is a decomposition of $V$.
We say that this decomposition is {\em $\Phi$-standard}.
\end{definition}

\begin{lemma} \label{lem:stand-basis} {\rm \cite[Lemma~5.1]{terw-lp24}}
Referring to the Leonard system $\Phi$ on $V$ from Definition \ref{def:ls},
let $v$ denote a nonzero vector in $E_0V$.
Then for $0 \leq i \leq d$ the element $E^*_iv$ is nonzero and hence a basis for $E^*_iV$.
Moreover the sequence $\{E^*_iv\}_{i=0}^d$ is a basis for $V$.
\end{lemma}

\begin{definition} \label{def:stand-basis} \rm \cite[Definition~5.2]{terw-lp24}
Referring to the Leonard system $\Phi$ on $V$ from Definition \ref{def:ls},
by a {\em $\Phi$-standard basis} for $V$ we mean the sequence $\{E^*_iv\}_{i=0}^d$ where $v$ denotes a nonzero vector in $E_0V$.
\end{definition}

\begin{lemma}
\label{lem:aastar-stand-basis}
Referring to the Leonard system $\Phi$ on $V$ from Definition \ref{def:ls},
with respect to a $\Phi$-standard basis
the matrix representing $A$ is irreducible tridiagonal and the matrix representing $A^*$ is ${\rm diag}(\theta^*_0, \theta^*_1, \ldots, \theta^*_d)$
where $\{\theta^*_i\}_{i=0}^d$ is the dual eigenvalue sequence of $\Phi$.
\end{lemma}

\noindent {\it Proof:} \rm
Immediate from Definition \ref{def:ls}.
\hfill $\Box$\\

\section{The split decomposition and the split basis}
Throughout this section let $\Phi$ denote the Leonard system on $V$ from Definition \ref{def:ls}.
For $0 \leq i \leq d$ define
\begin{eqnarray}
U_i = (E^*_0V + E^*_1V + \cdots + E^*_iV) \cap (E_iV + E_{i+1}V + \cdots + E_dV). \label{eq:split_decomp_def}
\end{eqnarray}
By \cite[Theorem~20.7]{madrid09} the sequence $\{U_i\}_{i=0}^d$ is a decomposition of $V$.
This decomposition is said to be {\em $\Phi$-split} \cite[Definition~20.2]{madrid09}.
By \cite[Theorem~20.7]{madrid09} for $0 \leq i \leq d$ both
\begin{align}
U_0 + U_1 + \cdots + U_i &= E^*_0V + E^*_1V + \cdots + E^*_iV, \nonumber\\%\label{eq:u_oi_sum}\\
U_i + U_{i+1} + \cdots + U_d &= E_iV + E_{i+1}V + \cdots + E_dV. \label{eq:u_id_sum}
\end{align}
By \cite[Lemma~20.9]{madrid09},
\begin{eqnarray}
(A-\theta_iI)U_i = U_{i+1} \qquad (0 \leq i \leq d-1), \qquad (A-\theta_dI)U_d = 0, \label{eq:a_raise_ui}\\
(A^*-\theta^*_iI)U_i = U_{i-1} \qquad (1 \leq i \leq d), \qquad (A^*-\theta^*_0I)U_0 = 0. \label{eq:astar_lower_ui}
\end{eqnarray}
By (\ref{eq:a_raise_ui}), (\ref{eq:astar_lower_ui}) for $1 \leq i \leq d$, $U_i$ is invariant under the action of 
$(A-\theta_{i-1}I)(A^*-\theta^*_iI)$,
and the corresponding eigenvalue
is a nonzero scalar in $\K$.
We denote this eigenvalue by $\varphi_i$.
We display a basis for $V$ that illuminates the 
significance of $\varphi_i$.
Setting $i=0$ in (\ref{eq:split_decomp_def}) we find $U_0=E^*_0V$.
Combining this with 
(\ref{eq:a_raise_ui}) we find
\begin{eqnarray}
U_i = (A-\theta_{i-1}I)\cdots (A-\theta_1I)(A-\theta_0I)E^*_0V
\qquad (0 \leq i \leq d).
\label{eq:uiinterp}
\end{eqnarray}
Let $v$ denote a nonzero vector in $E^*_0V$. From
(\ref{eq:uiinterp}) we find that for $0 \leq i \leq d$
the vector
$(A-\theta_{i-1}I)\cdots (A-\theta_0I)v$
is a basis for $U_i$. By this and since $\{U_i\}_{i=0}^d$
is a decomposition of $V$ we find the sequence
\begin{eqnarray}
(A-\theta_{i-1}I)\cdots (A-\theta_1I)(A-\theta_0I)v
\qquad (0 \leq i \leq d)
\label{eq:splitbasisdef}
\end{eqnarray}
is a basis for $V$. With respect to this basis the
matrices representing $A$ and $A^*$ are
\begin{equation}
A:
\left(
\begin{array}{c c c c c c}
\theta_0 & & & & & {\bf 0} \\
1 & \theta_1 &  & & & \\
& 1 & \theta_2 &  & & \\
& & \cdot & \cdot &  &  \\
& & & \cdot & \cdot &  \\
{\bf 0}& & & & 1 & \theta_d
\end{array}
\right),
\qquad
A^*:
\left(
\begin{array}{c c c c c c}
\theta^*_0 &\varphi_1 & & & & {\bf 0} \\
 & \theta^*_1 & \varphi_2 & & & \\
&  & \theta^*_2 & \cdot & & \\
& &  & \cdot & \cdot &  \\
& & &  & \cdot & \varphi_d \\
{\bf 0}& & & &  & \theta^*_d
\end{array}
\right).
\label{eq:matrepaastar}
\end{equation}
By a {\it $\Phi$-split basis} for $V$ we mean
a sequence
 of the form 
(\ref{eq:splitbasisdef}), where $v$ is a nonzero vector
in $E^*_0V$.
We call
$\{\varphi_i\}_{i=1}^d$
the {\it first split sequence} of $\Phi$.
We let 
$\{\phi_i\}_{i=1}^d$
denote the first split sequence
of $\Phi^\Downarrow$ and call
this the 
 {\it second split sequence} of $\Phi$.
For notational convenience define
$\varphi_0=0$, 
$\varphi_{d+1}=0$, 
$\phi_0=0$, 
$\phi_{d+1}=0$.

\noindent
We now define the parameter array of $\Phi$.

\begin{definition} \rm 
\cite[Definition~10.1]{tdd-lbub}
\label{def:parameter-array}
By the {\em parameter array} of $\Phi$ we mean the sequence 
$(\{\theta_i\}_{i=0}^d;\{\theta_i^*\}_{i=0}^d;\{\varphi_i\}_{i=1}^d;\{\phi_i\}_{i=1}^d)$
where $\{\theta_i\}_{i=0}^d$ (resp. $\{\theta^*_i\}_{i=0}^d$) denotes the eigenvalue sequence (resp. dual eigenvalue sequence) of $\Phi$ and $\{\varphi_i\}_{i=1}^d$ (resp. $\{\phi_i\}_{i=1}^d$) denotes the first split sequence (resp. second split sequence) of $\Phi$.
\end{definition}

\noindent We finish this section with a few characterizations of the $\Phi$-split basis.

\begin{lemma}
\label{lem:lp-split-basis-3}
{\rm \cite[Lemma~13.2]{tdd-lbub}}
Let $\{v_i\}_{i=0}^d$ denote a sequence of vectors in $V$, not all zero.
Then $\{v_i\}_{i=0}^d$ is a $\Phi$-split basis for $V$ if and only if both
{\rm (i)} $v_0 \in E^*_0V$;
{\rm (ii)} $Av_i = \theta_i v_i + v_{i+1}$ for $0 \leq i \leq d-1$.
\end{lemma}

\begin{lemma}
\label{lem:lp-split-basis-4}
Let $\{v_i\}_{i=0}^d$ denote a sequence of vectors in $V$, not all zero.
Then $\{v_i\}_{i=0}^d$ is a $\Phi$-split basis for $V$ if and only if both
{\rm (i)} $v_d \in E_dV$;
{\rm (ii)} $A^*v_i = \theta^*_i v_i + \varphi_i v_{i-1}$ for $1 \leq i \leq d$.
\end{lemma}

\noindent {\it Proof:} \rm
First assume $\{v_i\}_{i=0}^d$ is a $\Phi$-split basis for $V$.
With respect to $\{v_i\}_{i=0}^d$ the matrices representing $A, A^*$ satisfy (\ref{eq:matrepaastar}).
Therefore the basis $\{v_i\}_{i=0}^d$ satisfies (i), (ii), and we are done in one direction.
To prove the other direction assume $\{v_i\}_{i=0}^d$ satisfies (i), (ii).
We will invoke Lemma \ref{lem:lp-split-basis-3}.
To do this we need to verify that $\{v_i\}_{i=0}^d$ satisfies Lemma \ref{lem:lp-split-basis-3}(i), (ii).
By (\ref{eq:split_decomp_def}) we have $E_dV = U_d$.
By this and (i) we have $v_d \in U_d$.
By this, (ii) and (\ref{eq:astar_lower_ui}) we have $v_i \in U_i$ for $0 \leq i \leq d$.
In particular $v_0 \in U_0$.
By (\ref{eq:split_decomp_def}) we have $U_0 = E^*_0V$ so $v_0 \in E^*_0V$.
Therefore $\{v_i\}_{i=0}^d$ satisfies Lemma \ref{lem:lp-split-basis-3}(i).
By (ii) we have $(A-\theta_iI)v_i = \varphi_{i+1}^{-1} (A-\theta_iI)(A^*-\theta^*_{i+1}I)v_{i+1}$ for $0 \leq i \leq d-1$.
By this and the discussion below (\ref{eq:astar_lower_ui}) we have $(A-\theta_iI)v_i = v_{i+1}$.
Therefore $\{v_i\}_{i=0}^d$ satisfies Lemma \ref{lem:lp-split-basis-3}(ii).
By Lemma \ref{lem:lp-split-basis-3} the sequence $\{v_i\}_{i=0}^d$ is a $\Phi$-split basis for $V$.
\hfill $\Box$\\

\begin{lemma} \label{lem:lp-split-basis-2}
Let $\{v_i\}_{i=0}^d$ denote a sequence of vectors in $V$.  Then the following are equivalent:
\begin{enumerate}
\item[\rm (i)] The sequence $\{v_i\}_{i=0}^d$ is a $\Phi$-split basis for $V$.
\item[\rm (ii)] There exists a nonzero $w \in E_dV$ such that
\begin{eqnarray*}
v_i = \frac{(A^*-\theta^*_{i+1}I) \cdots (A^*-\theta^*_{d-1}I)(A^*-\theta^*_{d}I)w}{\varphi_{i+1} \cdots \varphi_{d-1} \varphi_d} \qquad (0 \leq i \leq d). \label{eq:v_i_split_2}
\end{eqnarray*}
\end{enumerate}
\end{lemma}

\noindent {\it Proof:} \rm
Immediate from Lemma \ref{lem:lp-split-basis-4}.
\hfill $\Box$\\

\section{A classification of Leonard systems}

In \cite[Theorem~1.9]{LS99} Leonard systems are classified up to isomorphism.
We now recall this classification.

\begin{theorem} 
{\rm \cite[Theorem~1.9]{LS99}}
\label{thm:ls}
Let $d$ denote a nonnegative integer and let
\begin{eqnarray}
(\{\theta_i\}_{i=0}^d; \{\theta_i^*\}_{i=0}^d; \{\varphi_i\}_{i=1}^d; \{\phi_i\}_{i=1}^d) \label{eq:pa_class}
\end{eqnarray}
denote a sequence of scalars taken from $\K$.
There exists a Leonard system $\Phi$ over $\K$ with parameter array {\rm (\ref{eq:pa_class})}
if and only if the following conditions {\rm (PA1)--(PA5)} hold.
\begin{description}
\item[(PA1)]  
$
\theta_i \not=\theta_j, \qquad 
\theta^*_i \not=\theta^*_j \qquad \mbox{if}\quad i\not=j, \qquad
(0 \leq i,j\leq d)
$.
\\
\item[(PA2)]
$
\varphi_i \not=0, \qquad \phi_i \not=0 \qquad (1 \leq i \leq d).
$
\\
\item[(PA3)]
$\displaystyle \varphi_i = \phi_1 \sum_{h=0}^{i-1} 
\frac{\theta_h-\theta_{d-h}}{\theta_0-\theta_d}
+ (\theta^*_i-\theta^*_0)(\theta_{i-1}-\theta_d)
\qquad (1 \leq i \leq d)$.
\\
\item[(PA4)] 
$\displaystyle \phi_i = \varphi_1 \sum_{h=0}^{i-1} 
\frac{\theta_h-\theta_{d-h}}{\theta_0-\theta_d}
+ (\theta^*_i-\theta^*_0)(\theta_{d-i+1}-\theta_0)
\qquad (1 \leq i \leq d)$.
\\
\item[(PA5)]
The expressions
\begin{eqnarray} 
\frac{\theta_{i-2}-\theta_{i+1}}{\theta_{i-1}-\theta_i},
\qquad 
\frac{\theta^*_{i-2}-\theta^*_{i+1}}{\theta^*_{i-1}-\theta^*_i}
\label{eq:betaplusone}
\end{eqnarray}
are equal and independent of $i$ for $2 \leq i \leq d-1$.
\end{description}
Moreover, suppose {\rm (PA1)--(PA5)} hold.
Then $\Phi$ is unique up to isomorphism of Leonard systems.
\end{theorem}

\begin{theorem}
{\rm \cite[Theorem~1.11]{LS99}}
\label{thm:param_d4}
Let $\Phi$ denote a Leonard system with parameter array
\begin{eqnarray*}
(\{\theta_i\}_{i=0}^d; \{\theta_i^*\}_{i=0}^d; \{\varphi_i\}_{i=1}^d; \{\phi_i\}_{i=1}^d).
\end{eqnarray*}
Then {\rm (i)--(iii)} hold below.
\begin{enumerate}
\item[\rm (i)] The parameter array of $\Phi^*$ is $(\{\theta^*_i\}_{i=0}^d; \{\theta_i\}_{i=0}^d; \{\varphi_i\}_{i=1}^d; \{\phi_{d-i+1}\}_{i=1}^d)$.
\item[\rm (ii)] The parameter array of $\Phi^\downarrow$ is $(\{\theta_i\}_{i=0}^d; \{\theta^*_{d-i}\}_{i=0}^d; \{\phi_{d-i+1}\}_{i=1}^d; \{\varphi_{d-i+1}\}_{i=1}^d)$.
\item[\rm (iii)] The parameter array of $\Phi^\Downarrow$ is $(\{\theta_{d-i}\}_{i=0}^d; \{\theta^*_i\}_{i=0}^d; \{\phi_i\}_{i=1}^d; \{\varphi_i\}_{i=1}^d)$.
\end{enumerate}
\end{theorem}

\section{The normalized split basis}

Throughout this section let $\Phi$ denote the Leonard system on $V$ from Definition \ref{def:ls}.
In an earlier section we discussed the $\Phi$-split basis.
For our purpose it is convenient to modify the $\Phi$-split basis by adjusting the normalization.

\begin{lemma} \label{lem:lp-norm-split-basis}
Let $v$ denote a nonzero vector in $E^*_0V$.
For $0 \leq i \leq d$ define
\begin{eqnarray}
u_i = \frac{\tau^*_i(\theta^*_d) (A-\theta_{i-1}I) \cdots (A-\theta_{1}I) (A-\theta_{0}I)v}{\varphi_1 \varphi_2 \cdots \varphi_i}. \label{eq:u_i_norm_split}
\end{eqnarray}
Then $u_i$ is a basis for the subspace $U_i$ from line {\rm (\ref{eq:split_decomp_def})}.
Moreover the sequence $\{u_i\}_{i=0}^d$ is a basis for $V$.
\end{lemma}

\noindent {\it Proof:} \rm
Since $\{\theta^*_j\}_{j=0}^d$ are mutually distinct, $\tau^*_i(\theta^*_d) \neq 0$.
The first assertion follows from this and the comment below (\ref{eq:uiinterp}).
The second assertion follows from this and the fact that $\{U_i\}_{i=0}^d$ is a decomposition of $V$.
\hfill $\Box$\\

\begin{definition}\rm
By a {\em normalized $\Phi$-split basis} for $V$ we mean a sequence $\{u_i\}_{i=0}^d$ of the form (\ref{eq:u_i_norm_split}),
where $v$ is a nonzero vector in $E^*_0V$.
\end{definition}

\noindent
For the rest of this section we describe the normalized $\Phi$-split basis from various points of view.

\begin{lemma} \label{lem:lp-connection-norm-split-basis}
The following {\rm (i), (ii)} hold.
\begin{enumerate}
\item[\rm (i)] Let $\{v_i\}_{i=0}^d$ denote a $\Phi$-split basis for $V$.  Then the sequence
\begin{eqnarray*}
\frac{\tau^*_i(\theta^*_d) v_i}{\varphi_1 \varphi_2 \cdots \varphi_i} \qquad (0 \leq i \leq d)
\end{eqnarray*}
is a normalized $\Phi$-split basis for $V$.
\item[\rm (ii)] Let $\{u_i\}_{i=0}^d$ denote a normalized $\Phi$-split basis for $V$.  Then the sequence
\begin{eqnarray*}
\frac{\varphi_1 \varphi_2 \cdots \varphi_i u_i }{\tau^*_i(\theta^*_d)} \qquad (0 \leq i \leq d)
\end{eqnarray*}\
is a $\Phi$-split basis for $V$.
\end{enumerate}
\end{lemma}

\noindent {\it Proof:} \rm
Compare (\ref{eq:splitbasisdef}) and (\ref{eq:u_i_norm_split}).
\hfill $\Box$\\

\begin{lemma} \label{lem:lp-norm-split-basis-2}
Let $\{u_i\}_{i=0}^d$ denote a sequence of vectors in $V$.  Then the following are equivalent:
\begin{enumerate}
\item[\rm (i)] The sequence $\{u_i\}_{i=0}^d$ is a normalized $\Phi$-split basis for $V$.
\item[\rm (ii)] There exists a nonzero $w \in E_dV$ such that
\begin{eqnarray*}
u_i = \tau^*_i(\theta^*_d) (A^*-\theta^*_{i+1}I) \cdots (A^*-\theta^*_{d-1}I) (A^*-\theta^*_{d}I)w \qquad (0 \leq i \leq d). \label{eq:u_i_norm_split_2}
\end{eqnarray*}
\end{enumerate}
\end{lemma}

\noindent {\it Proof:}
Immediate from Lemma \ref{lem:lp-split-basis-2} and Lemma \ref{lem:lp-connection-norm-split-basis}.
\hfill $\Box$\\

\begin{lemma} \label{lem:char_norm_split_basis_3}
Let $\{u_i\}_{i=0}^d$ denote a sequence of vectors in $V$, not all zero.  Then $\{u_i\}_{i=0}^d$ is a normalized $\Phi$-split basis for $V$ if and only if both {\rm (i)} $u_0 \in E^*_0 V$; {\rm (ii)} $Au_i = \theta_i u_i + \varphi_{i+1}(\theta^*_d - \theta^*_i)^{-1} u_{i+1}$ for $0 \leq i \leq d-1$.
\end{lemma}

\noindent {\it Proof:} \rm
Immediate from Lemma \ref{lem:lp-norm-split-basis}.
\hfill $\Box$\\

\begin{lemma} \label{lem:char_norm_split_basis_4}
Let $\{u_i\}_{i=0}^d$ denote a sequence of vectors in $V$, not all zero.  Then $\{u_i\}_{i=0}^d$ is a normalized $\Phi$-split basis for $V$ if and only if both {\rm (i)} $u_d \in E_d V$; {\rm (ii)} $A^*u_i = \theta_i^* u_i + (\theta_d^* - \theta_{i-1}^*) u_{i-1}$ for $1 \leq i \leq d$.
\end{lemma}

\noindent {\it Proof:} \rm
Immediate from Lemma \ref{lem:lp-norm-split-basis-2}.
\hfill $\Box$\\

\noindent
Using Lemma \ref{lem:char_norm_split_basis_3} and Lemma \ref{lem:char_norm_split_basis_4} we now describe the matrices  representing $A, A^*$ with respect to a normalized $\Phi$-split basis for $V$.

\begin{lemma} \label{lem:lp-norm-split-basis-aastar}
With respect to a normalized $\Phi$-split basis for $V$, the matrices in ${\rm Mat}_{d+1}(\fld)$ that represent $A, A^*$ are described as follows.  The matrix representing $A$ is lower bidiagonal with $(i,i)$-entry $\theta_i$ for $0 \leq i \leq d$ and $(i,i-1)$-entry $\varphi_i/(\theta_d^* - \theta_{i-1}^*)$ for $1\leq i \leq d$.
The matrix representing $A^*$ is upper bidiagonal with $(i,i)$-entry $\theta^*_i$ for $0 \leq i \leq d$ and $(i-1,i)$-entry $\theta^*_d - \theta^*_{i-1}$ for $1 \leq i \leq d$.
Moreover the matrix representing $A^*$ has constant row sum $\theta^*_d$.
\end{lemma}

\begin{example} \label{ex:lp-norm-split-basis-aastar}
With reference to Definition \ref{def:ls} assume $d = 4$.
With respect to a normalized $\Phi$-split basis for $V$, the matrices representing $A, A^*$ are given below.
\begin{eqnarray*}
A :
\left(
\begin{array}{cccccc}
\theta_0 & 0 & 0 & 0 & 0\\
\frac{\varphi_1}{\theta^*_4 - \theta^*_0} & \theta_1 & 0 & 0 & 0\\
0 & \frac{\varphi_2}{\theta^*_4 - \theta^*_1} & \theta_2 & 0 & 0\\
0 & 0 & \frac{\varphi_3}{\theta^*_4 - \theta^*_2} & \theta_3 & 0\\
0 & 0 & 0 & \frac{\varphi_4}{\theta^*_4 - \theta^*_3} & \theta_4
\end{array}
\right),
\quad
A^* :
\left(
\begin{array}{cccccc}
\theta_0^* & \theta_4^* - \theta_{0}^* & 0 & 0 & 0\\
0 & \theta_1^* & \theta_4^* - \theta_{1}^* & 0 & 0\\
0 & 0 & \theta_2^* & \theta_4^* - \theta_{2}^* & 0\\
0 & 0 & 0 & \theta_3^* & \theta_4^* - \theta_{3}^*\\
0 & 0 & 0 & 0 & \theta_4^*
\end{array}
\right).
\end{eqnarray*}
Observe that the matrix representing $A^*$ has constant row sum $\theta_4^*$.
\end{example}

\begin{lemma} \label{lem:char_norm_split_basis_5}
Let $\{u_i\}_{i=0}^d$ denote a sequence of vectors in $V$, not all zero.
Then $\{u_i\}_{i=0}^d$ is a normalized $\Phi$-split basis for $V$ if and only if both
{\rm(i)} $u_i \in U_i$ for $0 \leq i \leq d$;
{\rm(ii)} $\sum_{i=0}^d u_i \in E^*_dV$.
\end{lemma}

\noindent {\it Proof:} \rm
First assume that $\{u_i\}_{i=0}^d$ is a normalized $\Phi$-split basis for $V$.
The sequence $\{u_i\}_{i=0}^d$ satisfies (i) by Lemma \ref{lem:lp-norm-split-basis}.
By Lemma \ref{lem:lp-norm-split-basis-aastar} the matrix representing $A^*$ with respect to $\{u_i\}_{i=0}^d$ has constant row sum $\theta^*_d$.
Therefore $(A^*-\theta^*_dI) \sum_{i=0}^d u_i = 0$,
and thus the sequence $\{u_i\}_{i=0}^d$ satisfies (ii).
We have shown that $\{u_i\}_{i=0}^d$ satisfies (i), (ii),
and we are done in one direction.
To prove the other direction assume $\{u_i\}_{i=0}^d$ satisfies (i), (ii).
We will invoke Lemma \ref{lem:char_norm_split_basis_4}.
To do this it suffices to verify that $\{u_i\}_{i=0}^d$ satisfies Lemma \ref{lem:char_norm_split_basis_4}(i), (ii).
By assumption $u_d \in U_d$.
By (\ref{eq:split_decomp_def}) we have $U_d = E_dV$ so $u_d \in E_dV$.
Therefore $\{u_i\}_{i=0}^d$ satisfies Lemma \ref{lem:char_norm_split_basis_4}(i).
By (ii) we have
\begin{align*}
0 &= \sum_{i=0}^d (A^*-\theta^*_d I)u_i\\
&= \sum_{i=0}^d (A^*-\theta^*_i I)u_i + \sum_{i=0}^d (\theta^*_i-\theta^*_d)u_i\\
&= \sum_{i=0}^d (A^*-\theta^*_i I)u_i + \sum_{i=1}^d (\theta^*_{i-1}-\theta^*_d)u_{i-1}.
\end{align*}
By (i) and (\ref{eq:astar_lower_ui}) we have $(A^*-\theta^*_i I)u_i \in U_{i-1}$ for $1 \leq i \leq d$ and $(A^*-\theta^*_0 I)u_0 = 0$.
By (i), the comments above  and since $\{U_i\}_{i=0}^d$ is a decomposition of $V$, the sequence $\{u_i\}_{i=0}^d$ satisfies Lemma \ref{lem:char_norm_split_basis_4}(ii).
Therefore by Lemma \ref{lem:char_norm_split_basis_4} the sequence $\{u_i\}_{i=0}^d$ is a normalized $\Phi$-split basis for $V$.
\hfill $\Box$\\

\section{The intersection matrix}

Throughout this section let $\Phi$ denote the Leonard system from Definition \ref{def:ls}.
In this section we recall the intersection matrix of $\Phi$ and the dual intersection matrix of $\Phi$.

\begin{definition}\rm \label{def:int_mat}
Consider the matrix in ${\rm Mat}_{d+1}(\fld)$ that represents $A$ with respect to a $\Phi$-standard basis.
This matrix is irreducible tridiagonal by Lemma \ref{lem:aastar-stand-basis}.  This matrix will be written as
\begin{eqnarray*}
\left(
\begin{array}{cccccc}
a_0 & b_0 & & & & {\bf 0}\\
c_1 & a_1 & b_1 & & &\\
& c_2 & \cdot & \cdot & &\\
& & \cdot & \cdot & \cdot &\\
& & & \cdot & \cdot & b_{d-1}\\
{\bf 0} & & & & c_d & a_d\\
\end{array}
\right).
\end{eqnarray*}
We call this matrix the {\em intersection matrix} of $\Phi$.
For notational convenience define $b_d = 0$ and $c_0 = 0$.
We call $a_i, b_i, c_i$ $(0 \leq i \leq d)$ the {\em intersection numbers} of $\Phi$.
\end{definition}

\begin{definition}\rm \label{def:dual_int_mat}
By the {\em dual intersection matrix} of $\Phi$ we mean the intersection matrix for $\Phi^*$.
We call $a^*_i, b^*_i, c^*_i$ $(0 \leq i \leq d)$ the {\em dual intersection numbers} of $\Phi$.
\end{definition}

\begin{lemma} 
{\rm \cite[Lemma~11.2]{madrid09}}
\label{lem:sum_abc}
We have $a_i+b_i+c_i = \theta_0$ for $0 \leq i \leq d$.
\end{lemma}

\noindent
We now give explicit formulas for the intersection numbers and the dual intersection numbers.
To avoid trivialities assume $d \geq 1$.

\begin{lemma} {\rm \cite[Theorem~23.5]{madrid09}}
\label{lem:int_mat_entries}
The following {\rm (i), (ii)} hold.
\begin{enumerate}
\item[\rm (i)] $\displaystyle b_i = \varphi_{i+1} \frac{\tau^*_i(\theta^*_i)}{\tau^*_{i+1}(\theta^*_{i+1})} \qquad (0 \leq i \leq d-1)$.
\item[\rm (ii)] $\displaystyle c_i = \phi_{i} \frac{\eta^*_{d-i}(\theta^*_i)}{\eta^*_{d-i+1}(\theta^*_{i-1})} \qquad (1 \leq i \leq d)$.
\end{enumerate}
\end{lemma}

\begin{lemma} {\rm \cite[Theorem~23.6]{madrid09}}
\label{lem:int_mat_entries2}
We have
\begin{align*}
a_0 &= \theta_0 + \frac{\varphi_1}{\theta^*_0-\theta^*_1},\\
a_i &= \theta_i + \frac{\varphi_i}{\theta^*_i-\theta^*_{i-1}} +  \frac{\varphi_{i+1}}{\theta^*_i-\theta^*_{i+1}} \qquad (1 \leq i \leq d-1),\\
a_d &= \theta_d + \frac{\varphi_d}{\theta^*_d - \theta^*_{d-1}}.
\end{align*}
\end{lemma}

\begin{lemma} \label{lem:dual_int_mat_entries}
The following {\rm (i), (ii)} hold.
\begin{enumerate}
\item[\rm (i)] $\displaystyle b^*_i = \varphi_{i+1} \frac{\tau_i(\theta_i)}{\tau_{i+1}(\theta_{i+1})} \qquad (0 \leq i \leq d-1)$.
\item[\rm (ii)] $\displaystyle c^*_i = \phi_{d-i+1} \frac{\eta_{d-i}(\theta_i)}{\eta_{d-i+1}(\theta_{i-1})} \qquad (1 \leq i \leq d)$.
\end{enumerate}
\end{lemma}

\noindent {\it Proof:} \rm
Apply Lemma \ref{lem:int_mat_entries} to the Leonard system $\Phi^*$
and use Lemma \ref{thm:param_d4}(i) to obtain the result.
\hfill $\Box$\\

\begin{lemma}
\label{lem:dual_int_mat_entries2}
We have
\begin{align*}
a^*_0 &= \theta^*_0 + \frac{\varphi_1}{\theta_0-\theta_1},\\
a^*_i &= \theta^*_i + \frac{\varphi_i}{\theta_i-\theta_{i-1}} +  \frac{\varphi_{i+1}}{\theta_i-\theta_{i+1}} \qquad (1 \leq i \leq d-1),\\
a^*_d &= \theta^*_d + \frac{\varphi_d}{\theta_d-\theta_{d-1}}.
\end{align*}
\end{lemma}

\noindent {\it Proof:} \rm
Apply Lemma \ref{lem:int_mat_entries2} to the Leonard system $\Phi^*$
and use Lemma \ref{thm:param_d4}(i) to obtain the result.
\hfill $\Box$\\

\section{The tridiagonal relations and the Askey-Wilson relations}

Throughout this section let $\Phi$ denote the Leonard system from Definition \ref{def:ls}.
We recall the corresponding tridiagonal relations and Askey-Wilson relations.

\begin{lemma}
\label{lem:td-relns}
{\rm \cite[Theorem~10.1]{somealg}}
There exists a sequence of scalars $\beta, \gamma, \gamma^*, \varrho, \varrho^*$ taken from $\fld$ such that both
\begin{align}
\left[ A, A^2A^* - \beta AA^*A + A^*A^2 - \gamma (AA^* + A^*A) - \varrho A^* \right] & = 0, \label{eq:td-reln1}\\
\left[ A^*, A^{*2}A - \beta A^*AA^* + AA^{*2} - \gamma^* (A^*A + AA^*) - \varrho^* A \right] & = 0. \label{eq:td-reln2}
\end{align}
The notation $\left[ r,s \right]$ means $rs - sr$.
The sequence is uniquely determined by the Leonard system $\Phi$ provided $d \geq 3$.
\end{lemma}

\noindent
We call (\ref{eq:td-reln1}), (\ref{eq:td-reln2}) the {\em tridiagonal relations}.

\begin{lemma} {\rm \cite[Corollary~4.4,~Theorem~4.5]{lpaw}}
\label{lem:lp-td-scalars}
Let $\beta,\gamma,\gamma^*,\varrho,\varrho^*$ denote scalars in $\fld$.
Then these scalars satisfy {\rm (\ref{eq:td-reln1}), (\ref{eq:td-reln2})} if and only if the following {\rm (i)--(v)} hold.
\begin{enumerate}
\item[\rm (i)] $\beta + 1$ is the common value of {\rm (\ref{eq:betaplusone})} for $2 \leq i \leq d-1$.
\item[\rm (ii)] $\gamma = \theta_{i-1} - \beta \theta_{i} + \theta_{i+1} \qquad (1 \leq i \leq d-1)$.
\item[\rm (iii)] $\gamma^* = \theta_{i-1}^* - \beta \theta_{i}^* + \theta_{i+1}^* \qquad (1 \leq i \leq d-1)$.
\item[\rm (iv)] $\varrho = \theta_{i-1}^2 - \beta \theta_{i-1} \theta_{i} + \theta_{i}^2 - \gamma(\theta_{i-1}+\theta_{i}) \qquad (1 \leq i \leq d)$.
\item[\rm (v)] $\varrho^* = \theta_{i-1}^{*2} - \beta \theta_{i-1}^* \theta_{i}^* + \theta_{i}^{*2} - \gamma^*(\theta_{i-1}^*+\theta_{i}^*) \qquad (1 \leq i \leq d)$.
\end{enumerate}
\end{lemma}

\begin{lemma} {\rm \cite[Theorem~1.5]{lpaw}}
\label{lem:lpaw}
Let $\beta,\gamma,\gamma^*,\varrho,\varrho^*$ denote the scalars from Lemma \ref{lem:td-relns}.  Then there exists a sequence of scalars $\omega, \eta, \eta^*$ taken from $\K$ such that both
\begin{align}
A^2A^* - \beta AA^*A + A^*A^2 - \gamma (AA^*+A^*A) - \varrho A^* &= \gamma^*A^2 + \omega A + \eta I \label{eq:aw_rel_1},\\
A^{*2}A - \beta A^*AA^* + AA^{*2} - \gamma^* (A^*A+AA^*) - \varrho^* A &= \gamma A^{*2} + \omega A^* + \eta^* I \label{eq:aw_rel_2}.
\end{align}
The sequence is uniquely determined by the Leonard system $\Phi$ provided $d \geq 3$.
\end{lemma}

\noindent
We call (\ref{eq:aw_rel_1}), (\ref{eq:aw_rel_2}) the {\em Askey-Wilson relations}.

\medskip
\noindent
For notational convenience let $\theta_{-1}$ and $\theta_{d+1}$ (resp. $\theta^*_{-1}$ and $\theta^*_{d+1}$) denote the scalars in $\fld$ which satisfy Lemma \ref{lem:lp-td-scalars}{\rm (ii)} (resp. Lemma \ref{lem:lp-td-scalars}{\rm (iii)}) for $i = 0$ and $i = d$.

\begin{lemma} {\rm \cite[Corollary~5.2,~Theorem~5.3]{lpaw}}
\label{lem:lpawconsts}
Let $\beta,\gamma,\gamma^*,\varrho,\varrho^*$ denote the scalars from Lemma \ref{lem:td-relns}.
Let $\omega, \eta, \eta^*$ denote scalars in $\fld$.
Then $\omega, \eta, \eta^*$ satisfy {\rm (\ref{eq:aw_rel_1}), (\ref{eq:aw_rel_2})} if and only if the following {\rm (i)--(iv)} hold.
\begin{enumerate}
\item[\rm (i)] $\omega = a_i(\theta_i^* - \theta_{i+1}^*) + a_{i-1} (\theta_{i-1}^* - \theta_{i-2}^*) - \gamma (\theta_{i}^* + \theta_{i-1}^*) \qquad (1 \leq i \leq d).$
\item[\rm (ii)] $\omega = a_i^*(\theta_i - \theta_{i+1}) + a_{i-1}^* (\theta_{i-1} - \theta_{i-2}) - \gamma^* (\theta_{i} + \theta_{i-1}) \qquad (1 \leq i \leq d).$
\item[\rm (iii)] $\eta = a_i^* (\theta_{i} - \theta_{i-1})(\theta_{i} - \theta_{i+1}) - \gamma^* \theta_i^2 - \omega \theta_{i} \qquad (0 \leq i \leq d).$
\item[\rm (iv)] $\eta^* = a_i (\theta_{i}^* - \theta_{i-1}^*)(\theta_{i}^* - \theta_{i+1}^*) - \gamma \theta_i^{*2} - \omega \theta_{i}^* \qquad (0 \leq i \leq d).$
\end{enumerate}
\end{lemma}

\section{Leonard systems of dual $q$-Krawtchouk type}\label{sec:ls-dqk}

For the past few sections we have been discussing general Leonard systems.
We now consider a family of Leonard systems said to have dual $q$-Krawtchouk type.
\begin{eqnarray*}
\mbox{For the rest of the paper let $q$ denote a nonzero scalar in $\K$ such that $q^2 \neq 1$.}
\end{eqnarray*}
We will be discussing the Leonard system $\Phi$ from Definition \ref{def:ls}.
Let $\overline{\fld}$ denote the algebraic closure of $\fld$.

\begin{definition} \rm
\cite[Example~35.8]{madrid09}
\label{def:dualqkrawtchouk2}
The Leonard system $\Phi$ is said to have {\it dual $q$-Krawtchouk} type whenever there exist scalars $h,h^*,\kappa,\kappa^*,\upsilon$ in $\overline{\fld}$ such that $\kappa,\kappa^*,\upsilon$ are nonzero and both
\begin{align}
\theta_i & = h + \kappa q^{d-2i} + \upsilon q^{2i-d}, \label{eq:lp-eval}\\
\theta_i^* & = h^* + \kappa^*q^{d-2i} \label{eq:lp-dualeval}
\end{align}
for $0 \leq i \leq d$ and both
\begin{align}
\varphi_i & = \kappa\kappa^*q^{d+1-2i}(q^i-q^{-i})(q^{i-d-1}-q^{d+1-i}), \label{eq:lp-split1}\\
\phi_i & = \kappa^*\upsilon q^{d+1-2i}(q^i-q^{-i})(q^{i-d-1}-q^{d+1-i}) \label{eq:lp-split2}
\end{align}
for $1 \leq i \leq d$.
\end{definition}

\noindent
For the rest of this section assume $\Phi$ has dual $q$-Krawtchouk type with the scalars $h,h^*,\kappa,\kappa^*,\upsilon$ as in Definition \ref{def:dualqkrawtchouk2}.

\begin{lemma}
\label{lem:diffofevals}
The following {\rm (i), (ii)} hold for $0 \leq i,j \leq d$.
\begin{enumerate}
\item[\rm (i)] $\theta_i - \theta_j = q^d(q^{-i}-q^{-j})(q^{-i}+q^{-j})(\kappa-\upsilon q^{2i+2j-2d})$.
\item[\rm (ii)] $\theta_i^* - \theta_j^* = \kappa^*q^d(q^{-i}-q^{-j})(q^{-i}+q^{-j})$.
\end{enumerate}
\end{lemma}

\noindent {\it Proof:} \rm
Immediate from (\ref{eq:lp-eval}) and (\ref{eq:lp-dualeval}).
\hfill $\Box$\\

\begin{lemma}
\label{lem:q_cond}
The following {\rm (i), (ii)} hold.
\begin{enumerate}
\item[\rm (i)] $q^{2i} \neq 1 \qquad (1 \leq i \leq d)$.
\item[\rm (ii)] $\kappa \neq \upsilon q^{2i-2d} \qquad (1 \leq i \leq 2d-1)$.
\end{enumerate}
\end{lemma}

\noindent {\it Proof:} \rm
By Lemma \ref{lem:diffofevals}(i) and since $\{\theta_i\}_{i=0}^d$ are mutually distinct.
\hfill $\Box$\\

\begin{lemma}
\label{lem:dqk-bca}
The intersection numbers of $\Phi$ satisfy the following {\rm (i)--(iii)}.
\begin{enumerate}
\item[\rm (i)] $b_i = \kappa q^i(q^{d-i}-q^{i-d}) \qquad (0 \leq i \leq d-1)$.
\item[\rm (ii)] $c_i = \upsilon q^{i-d}(q^{-i}-q^i) \qquad (1 \leq i \leq d)$.
\item[\rm (iii)] $a_i = h + (\kappa+\upsilon)q^{2i-d} \qquad (0 \leq i \leq d)$.
\end{enumerate}
\end{lemma}

\noindent {\it Proof:} \rm
In the equations of Lemma \ref{lem:int_mat_entries} and Lemma \ref{lem:int_mat_entries2}, evaluate the right-hand sides using Definition \ref{def:dualqkrawtchouk2} and Lemma \ref{lem:diffofevals}(ii), and then simplify the result.
\hfill $\Box$\\

\begin{lemma}
\label{lem:dqk-bcastar}
The dual intersection numbers of $\Phi$ satisfy the following {\rm (i)--(iii)}.
\begin{enumerate}
\item[\rm (i)] $\displaystyle b_i^* = \frac{\kappa\kappa^*q^i(q^{d-i}-q^{i-d})(\kappa-\upsilon q^{2i-2d})}{(\kappa-\upsilon q^{4i-2d})(\kappa-\upsilon q^{4i-2d+2})} \qquad (1 \leq i \leq d-1),\qquad
b^*_0 = \frac{\kappa \kappa^*(q^d-q^{-d})}{\kappa - \upsilon q^{2-2d}}.$
\item[\rm (ii)] $\displaystyle c_i^* = \frac{\upsilon\kappa^*q^{5i-3d-2}(q^{-i}-q^{i})(\kappa-\upsilon q^{2i})}{(\kappa-\upsilon q^{4i-2d})(\kappa-\upsilon q^{4i-2d-2})} \qquad (1 \leq i \leq d-1), \qquad
c^*_d = \frac{\upsilon \kappa^* q^{2d-2}(q^{-d}-q^d)}{\kappa - \upsilon q^{2d-2}}$.
\item[\rm (iii)] $a_i^* = \theta_0^* - b_i^* - c_i^* \qquad (0 \leq i \leq d)$.
\end{enumerate}
\end{lemma}

\noindent {\it Proof:} \rm
(i), (ii) In the equations of Lemma \ref{lem:dual_int_mat_entries}, evaluate the right-hand sides using Definition \ref{def:dualqkrawtchouk2} and Lemma \ref{lem:diffofevals}(i), and then simplify the result.\\
(iii) Apply Lemma \ref{lem:sum_abc} to the Leonard system $\Phi^*$.
\hfill $\Box$\\

\begin{lemma}
\label{lem:lp-dkrawtchouk-consts}
Let $\beta,\gamma,\gamma^*,\varrho,\varrho^*$ denote the scalars from $\fld$ which satisfy the following {\rm (i)--(v)}.
\begin{enumerate}
\item[\rm (i)] $\beta = q^2 + q^{-2}$.
\item[\rm (ii)] $\gamma = h(2-\beta)$.
\item[\rm (iii)] $\gamma^* = h^*(2-\beta)$.
\item[\rm (iv)] $\varrho = h^2(\beta-2)-\kappa \upsilon(\beta^2-4)$.
\item[\rm (v)] $\varrho^* = h^{*2}(\beta-2)$.
\end{enumerate}
Then $\beta,\gamma,\gamma^*,\varrho,\varrho^*$ satisfy {\rm (\ref{eq:td-reln1}), (\ref{eq:td-reln2})}.
\end{lemma}

\noindent {\it Proof:} \rm
Routine verification using Lemma \ref{lem:lp-td-scalars} and Definition \ref{def:dualqkrawtchouk2}.
\hfill $\Box$\\

\begin{lemma}
\label{lem:lp-dkrawtchouk-awconsts}
Let $\beta,\gamma,\gamma^*,\varrho,\varrho^*$ denote the scalars from Lemma \ref{lem:lp-dkrawtchouk-consts}.
Let $\omega, \eta, \eta^*$ denote the scalars from $\fld$ which satisfy the following {\rm (i)--(iii)}.
\begin{enumerate}
\item[\rm (i)] $\omega = (\beta-2)(2hh^*-(\kappa+\upsilon) \kappa^*).$
\item[\rm (ii)] $\eta = \kappa \upsilon h^*(\beta^2-4) + \kappa \upsilon \kappa^*(q+q^{-1})(q-q^{-1})^2(q^{d+1}+q^{-d-1}) - h(\beta-2)(hh^*-(\kappa+\upsilon) \kappa^*).$
\item[\rm (iii)] $\eta^* = h^*(\beta-2)((\kappa+\upsilon) \kappa^*-hh^*).$
\end{enumerate}
Then $\beta,\gamma,\gamma^*,\varrho,\varrho^*,\omega,\eta,\eta^*$ satisfy {\rm (\ref{eq:aw_rel_1}), (\ref{eq:aw_rel_2})}.
\end{lemma}

\noindent {\it Proof:} \rm
In the equations of Lemma \ref{lem:lpawconsts}, evaluate the right-hand sides using Definition \ref{def:dualqkrawtchouk2}, Lemma \ref{lem:dqk-bca}, Lemma \ref{lem:dqk-bcastar} and Lemma \ref{lem:lp-dkrawtchouk-consts}, and then simplify the result.
\hfill $\Box$\\

\begin{note} \rm
\label{note:swap-uv}
Among the relatives of $\Phi$ we find that the Leonard system $\Phi^\Downarrow$ also has dual $q$-Krawtchouk type with $h^\Downarrow=h, (h^*)^\Downarrow=h^*, \kappa^\Downarrow=\upsilon, (\kappa^*)^\Downarrow=\kappa^*, \upsilon^\Downarrow=\kappa$.
\end{note}

\section{The algebra $U_q(\mathfrak{sl}_2)$}

In the previous section we discussed Leonard systems of dual $q$-Krawtchouk type.  We now turn our attention to the algebra $U_q(\mathfrak{sl}_2)$.  Later we will relate the Leonard systems of dual $q$-Krawtchouk type and the algebra $U_q(\mathfrak{sl}_2)$.

\begin{definition} \rm
	\label{def:uqsl2chevalley}
	Let $U_q(\mathfrak{sl}_2)$ denote the $\fld$-algebra with generators $k^{\pm 1},e,f$ and relations
	\begin{eqnarray}
		&&\qquad kk^{-1} = k^{-1}k = 1,\nonumber\\%\label{eq:uqsl2chevalley_1}\\
		&&ke = q^2ek, \quad \quad
		kf = q^{-2}fk,\nonumber\\%\label{eq:uqsl2chevalley_3}\\
		&&\qquad ef - fe = \frac{k-k^{-1}}{q-q^{-1}}.\nonumber%\label{eq:uqsl2chevalley_4}
	\end{eqnarray}
\end{definition}

\begin{lemma} {\rm \cite[Theorem~2.1]{equit1}}
	\label{lem:uqsl2equitable}
	The algebra $U_q(\mathfrak{sl}_2)$ is isomorphic to the $\fld$-algebra with generators $x,y,z^{\pm 1}$ and relations
	\begin{eqnarray}
		&&zz^{-1} = z^{-1}z = 1,	\label{eq:uqsl2zzinverse}\\
		&&\frac{qxy-q^{-1}yx}{q-q^{-1}} = 1,	\label{eq:uqsl2xy}\\
		&&\frac{qyz-q^{-1}zy}{q-q^{-1}} = 1,	\label{eq:uqsl2yz}\\
		&&\frac{qzx-q^{-1}xz}{q-q^{-1}} = 1.	\label{eq:uqsl2zx}
	\end{eqnarray}
	An isomorphism with the presentation in Definition \ref{def:uqsl2chevalley} is given by
	\begin{eqnarray*}
		z^{\pm 1} & \rightarrow & k^{\pm 1},\\
		y & \rightarrow & k^{-1} - q(q-q^{-1})k^{-1}e,\\
		x & \rightarrow & k^{-1} + (q-q^{-1})f.
	\end{eqnarray*}
	The inverse of this isomorphism is given by
	\begin{eqnarray*}
		k^{\pm 1} & \rightarrow & z^{\pm 1},\\
		f & \rightarrow & (q-q^{-1})^{-1}(x-z^{-1}),\\
		e & \rightarrow & q^{-1}(q-q^{-1})^{-1}(1-zy).
	\end{eqnarray*}
\end{lemma}

\noindent
We call $x,y,z^{\pm 1}$ the {\it equitable generators} for $U_q(\mathfrak{sl}_2)$.
From now on we identify the versions of $U_q(\mathfrak{sl}_2)$ in Definition \ref{def:uqsl2chevalley} and Lemma \ref{lem:uqsl2equitable} via the isomorphism in Lemma \ref{lem:uqsl2equitable}.\\

\noindent
We now discuss finite-dimensional $U_q(\mathfrak{sl}_2)$-modules.
For an integer $n$ we define
\begin{eqnarray*}
[n]_q = \frac{q^n - q^{-n}}{q-q^{-1}}.
\end{eqnarray*}

\begin{lemma} {\rm\cite[Theorem~2.6]{jantzen}}
\label{lem:fd-u-module-kfe}
For an integer $d \geq 0$ and for $\epsilon \in \{-1,1\}$ there exists a $U_q(\mathfrak{sl}_2)$-module $L(d,\epsilon)$ with the following properties.
$L(d,\epsilon)$ has a basis $\{v_i\}_{i=0}^d$ such that
\begin{align*}
kv_i &= \epsilon q^{d-2i}v_i \qquad (0 \leq i \leq d),\\
fv_i &= [i+1]_q v_{i+1} \qquad (0 \leq i \leq d-1), \qquad fv_d = 0,\\
ev_i &= \epsilon [d-i+1]_q v_{i-1} \qquad (1 \leq i \leq d), \qquad ev_0 = 0.
\end{align*}
The $U_q(\mathfrak{sl}_2)$-module $L(d,\epsilon)$ is irreducible provided that $q^{2i} \neq 1$ for $1 \leq i \leq d$.
Assume $q^{2i} \neq 1$ for $1 \leq i \leq d$.
Then every irreducible $U_q(\mathfrak{sl}_2)$-module of dimension $d+1$ is isomorphic to either $L(d,1)$ or $L(d,-1)$.
\end{lemma}

\begin{lemma}
	{\rm \cite[Lemma~4.2]{equit1}}
	\label{lem:lnep-xyz}
	For an integer $d \geq 0$ and for $\epsilon \in \{-1,1\}$, the $U_q(\mathfrak{sl}_2)$-module $L(d,\epsilon)$ has a basis $\{u_i\}_{i=0}^d$ such that
	\begin{align*}
		\epsilon xu_i &= q^{d-2i}u_i \qquad (0 \leq i \leq d),\\%\label{eq:lnep-xyz-x}\\
		(\epsilon y - q^{2i-d})u_i &= (q^{-d} - q^{2i+2-d})u_{i+1} \qquad (0 \leq i \leq d-1), \qquad (\epsilon y - q^d)u_d = 0,\\%\label{eq:lnep-xyz-y}\\
		(\epsilon z - q^{2i-d})u_i &= (q^d - q^{2i-2-d})u_{i-1} \qquad (1 \leq i \leq d), \qquad (\epsilon z - q^{-d})u_0 = 0.%\label{eq:lnep-xyz-z}
	\end{align*}
\end{lemma}

\begin{definition} \label{def:norm-x-e-basis} \rm
We call the basis $\{u_i\}_{i=0}^d$ from Lemma \ref{lem:lnep-xyz} a {\em normalized $x$-eigenbasis} for $L(d,\epsilon)$.
Observe that this basis satisfies $\epsilon yu = q^{-d}u$ and $\epsilon zu = q^du$ where $u = \sum_{i=0}^d u_i$.
\end{definition}

\begin{lemma}
	\label{lem:lnep-yzx}
	{\rm \cite[Lemma~4.2]{equit1}}
	For an integer $d \geq 0$ and for $\epsilon \in \{-1,1\}$, the $U_q(\mathfrak{sl}_2)$-module $L(d,\epsilon)$ has a basis $\{u_i\}_{i=0}^d$ such that
	\begin{align}
		\epsilon yu_i &= q^{d-2i}u_i \qquad (0 \leq i \leq d), \label{eq:lnep-yzx-y}\\
		(\epsilon z - q^{2i-d})u_i &= (q^{-d} - q^{2i+2-d})u_{i+1} \qquad (0 \leq i \leq d-1), \qquad (\epsilon z - q^d)u_d = 0, \label{eq:lnep-yzx-z} \\
		(\epsilon x - q^{2i-d})u_i &= (q^d - q^{2i-2-d})u_{i-1} \qquad (1 \leq i \leq d), \qquad (\epsilon x - q^{-d})u_0 = 0. \label{eq:lnep-yzx-x}
	\end{align}
\end{lemma}

\begin{definition} \label{def:norm-y-e-basis} \rm
We call the basis $\{u_i\}_{i=0}^d$ from Lemma \ref{lem:lnep-yzx} a {\em normalized $y$-eigenbasis} for $L(d,\epsilon)$.
Observe that this basis satisfies $\epsilon zu = q^{-d}u$ and $\epsilon xu = q^du$ where $u = \sum_{i=0}^d u_i$.
\end{definition}

\begin{lemma}
	\label{lem:lnep-zxy}
	{\rm \cite[Lemma~4.2]{equit1}}
	For an integer $d \geq 0$ and for $\epsilon \in \{-1,1\}$, the $U_q(\mathfrak{sl}_2)$-module $L(d,\epsilon)$ has a basis $\{u_i\}_{i=0}^d$ such that
	\begin{align*}
		\epsilon zu_i & = q^{d-2i}u_i \qquad (0 \leq i \leq d),\\%\label{eq:lnep-zxy-z}\\
		(\epsilon x - q^{2i-d})u_i & = (q^{-d} - q^{2i+2-d})u_{i+1} \qquad (0 \leq i \leq d-1), \qquad (\epsilon x - q^d)u_d = 0,\\%\label{eq:lnep-zxy-x}\\
		(\epsilon y - q^{2i-d})u_i & = (q^d - q^{2i-2-d})u_{i-1} \qquad (1 \leq i \leq d), \qquad (\epsilon y - q^{-d})u_0 = 0.%\label{eq:lnep-zxy-y}
	\end{align*}
\end{lemma}

\begin{definition} \label{def:norm-z-e-basis} \rm
We call the basis $\{u_i\}_{i=0}^d$ from Lemma \ref{lem:lnep-zxy} a {\em normalized $z$-eigenbasis} for $L(d,\epsilon)$.
Observe that this basis satisfies $\epsilon xu = q^{-d}u$ and $\epsilon yu = q^du$ where $u = \sum_{i=0}^d u_i$.
\end{definition}

\begin{definition} \rm
\cite[p.~21]{jantzen}
\label{def:casimir}
Let
\begin{eqnarray*}
\Delta = ef + \frac{q^{-1}k+qk^{-1}}{(q-q^{-1})^2}.
\end{eqnarray*}
We call $\Delta$ the {\em Casimir element} of $U_q(\mathfrak{sl}_2)$.
\end{definition}

\begin{lemma}
{\rm \cite[Lemma~2.7]{jantzen}}
The element $\Delta$ is central in $U_q(\mathfrak{sl}_2)$.
\end{lemma}

\begin{lemma}
{\rm \cite[Lemma~2.7]{jantzen}}
\label{lem:casimir-scalar}
For an integer $d \geq 0$ and for $\epsilon \in \{-1,1\}$, the element $\Delta$ acts on the $U_q(\mathfrak{sl}_2)$-module $L(d,\epsilon)$ as the identity times
\begin{eqnarray}
\epsilon \frac{q^{d+1}+q^{-d-1}}{(q-q^{-1})^2}.\label{eq:casimir-scalar}
\end{eqnarray}
\end{lemma}

\noindent
For the rest of this section suppose $q$ is not a root of unity and ${\rm char}(\K) \neq 2$.

\begin{lemma}
{\rm \cite[Lemma~2.7]{jantzen}}
\label{lem:casimir-scalar-distinct}
For all integers $d \geq 0$ and for all $\epsilon \in \{-1,1\}$, the scalars {\rm (\ref{eq:casimir-scalar})} are mutually distinct.
\end{lemma}

\noindent
Let $M$ denote a finite-dimensional $U_q(\mathfrak{sl}_2)$-module.
For an integer $d \geq 0$ and for $\epsilon \in \{-1,1\}$
let $M_{d,\epsilon}$ denote the subspace of $M$ spanned by the irreducible $U_q(\mathfrak{sl}_2)$-submodules of $M$ which are isomorphic to $L(d,\epsilon)$.
Observe that $M_{d,\epsilon}$ is a $U_q(\mathfrak{sl}_2)$-submodule of $M$.
We call $M_{d,\epsilon}$ the {\em homogeneous component of M associated with $d$ and $\epsilon$.}
By \cite[p.~22]{jantzen} the homogeneous component $M_{d,\epsilon}$ is the eigenspace for $\Delta$ associated with eigenvalue (\ref{eq:casimir-scalar}).
Moreover by \cite[p.~22]{jantzen},
\begin{eqnarray}
M = \sum M_{d,\epsilon} \qquad {\rm (direct \ sum)},\label{eq:u_homeg_decomp}
\end{eqnarray}
where the sum is over all integers $d \geq 0$ and $\epsilon \in \{-1,1\}$.

\medskip
\noindent
We emphasize one point for later use.

\begin{lemma}
{\rm \cite[Theorem~2.9]{jantzen}}
\label{lem:umod-ss}
Every finite-dimensional $U_q(\mathfrak{sl}_2)$-module is semisimple.
\end{lemma}

\section{$U_q(\mathfrak{sl}_2)$ and Leonard systems of dual $q$-Krawtchouk type}

In this section we display two $U_q(\mathfrak{sl}_2)$-module structures associated with a given Leonard system of dual $q$-Krawtchouk type.
Prior to this display we make some comments.
Throughout this section let $V$ denote a vector space over $\K$ with finite positive dimension.  Let $A:V \rightarrow V$ denote a linear transformation.
For $\theta \in \K$ define
$V_A(\theta) = \{ v \in V | Av = \theta v \}$.

\begin{lemma} {\rm \cite[Lemma~6.2]{q-inv}}
\label{lem:q-weyl-linalg1}
Let $A: V \rightarrow V$ and $B: V \rightarrow V$ denote linear transformations.  Then for all nonzero $\theta \in \K$ the following are equivalent:
\begin{enumerate}
	\item[\rm (i)] The expression $qAB-q^{-1}BA-(q-q^{-1})I$ vanishes on $V_A(\theta)$.
	\item[\rm (ii)] $(B-\theta^{-1}I)V_A(\theta) \subseteq V_A(q^{-2}\theta)$.
\end{enumerate}
\end{lemma}

\begin{lemma} {\rm \cite[Lemma~6.3]{q-inv}}
\label{lem:q-weyl-linalg2}
Let $A: V \rightarrow V$ and $B: V \rightarrow V$ denote linear transformations.  Then for all nonzero $\theta \in \K$ the following are equivalent:
\begin{enumerate}
	\item[\rm (i)] The expression $qAB-q^{-1}BA-(q-q^{-1})I$ vanishes on $V_B(\theta)$.
	\item[\rm (ii)] $(A-\theta^{-1}I)V_B(\theta) \subseteq V_B(q^{2}\theta)$.
\end{enumerate}
\end{lemma}

\noindent
For the rest of this section let $\Phi$ denote a Leonard system on $V$ as in Definition \ref{def:ls}.  Assume $\Phi$ has dual $q$-Krawtchouk type. Let $h,h^*,\kappa,\kappa^*,\upsilon$ denote the corresponding parameters from Definition \ref{def:dualqkrawtchouk2}.
We are going to define something that turns out to be unique up to sign.
For notational convenience we fix $\epsilon \in \{1,-1\}$.

\begin{definition} \label{def:xyz-lintranfs} \rm
Let $\{X_i\}_{i=0}^d$ denote the $\Phi^\downarrow$-split decomposition of $V$.
Let $\{Y_i\}_{i=0}^d$ denote the inverted $\Phi$-split decomposition of $V$.
Let $\{Z_i\}_{i=0}^d$ denote the $\Phi$-standard decomposition of $V$.
Define the linear transformations $X,Y,Z$ in ${\rm End}(V)$ as follows.
\begin{align}
(\epsilon X-q^{d-2i}I)X_i &= 0 \qquad (0 \leq i \leq d), \label{eq:x-map}\\
(\epsilon Y-q^{d-2i}I)Y_i &= 0 \qquad (0 \leq i \leq d), \label{eq:y-map}\\
(\epsilon Z-q^{d-2i}I)Z_i &= 0 \qquad (0 \leq i \leq d). \label{eq:z-map}
\end{align} 
\end{definition}

\begin{lemma}
\label{lem:xyz-lintranfs-uqsl2}
There exists a unique $U_q(\mathfrak{sl}_2)$-module structure on $V$ such that the equitable generators $x,y,z$ of $U_q(\mathfrak{sl}_2)$ act on $V$ as $X,Y,Z$, respectively.
\end{lemma}

\noindent {\it Proof:} \rm
By construction $Z$ is invertible.
We show that $X,Y,Z,Z^{-1}$ satisfy the defining relations (\ref{eq:uqsl2zzinverse})--(\ref{eq:uqsl2zx}) of $U_q(\mathfrak{sl}_2)$.
By construction $Z, Z^{-1}$ satisfy (\ref{eq:uqsl2zzinverse}).
Next we show that $Y,Z$ satisfy (\ref{eq:uqsl2yz}).
Since $\{Z_i\}_{i=0}^d$ is the $\Phi$-standard decomposition of $V$, we have $(A^* - \theta^*_iI)Z_i = 0$ for $0 \leq i \leq d$.
Using (\ref{eq:lp-dualeval}) to compare this with (\ref{eq:z-map}) we find
\begin{eqnarray}
\epsilon Z = \kappa^{*-1}(A^*-h^*I).\label{eq:z_to_astar}
\end{eqnarray}
By (\ref{eq:astar_lower_ui}) and since $\{Y_i\}_{i=0}^d$ is the inverted $\Phi$-split decomposition of $V$, we have $(A^*-\theta_{d-i}^*I)Y_i = Y_{i+1}$ for $0 \leq i \leq d-1$ and $(A^*-\theta_0^*I)Y_d = 0$.
By this, (\ref{eq:lp-dualeval}) and (\ref{eq:z_to_astar}) we have $(\epsilon Z-q^{2i-d}I)Y_i = Y_{i+1}$ for $0 \leq i \leq d-1$ and $(\epsilon Z-q^{d}I)Y_d = 0$.
By this, Lemma \ref{lem:q-weyl-linalg1} and since $\{Y_i\}_{i=0}^d$ is a decomposition of $V$, we obtain that $Y,Z$ satisfy (\ref{eq:uqsl2yz}).
Next we show that $Z,X$ satisfy (\ref{eq:uqsl2zx}).
By applying (\ref{eq:astar_lower_ui}) to $\Phi^\downarrow$ and since $\{X_i\}_{i=0}^d$ is the $\Phi^\downarrow$-split decomposition of $V$, we have $(A^*-\theta_{d-i}^*I)X_i = X_{i-1}$ for $1 \leq i \leq d$ and $(A^*-\theta_d^*I)X_0 = 0$.
By this, (\ref{eq:lp-dualeval}) and (\ref{eq:z_to_astar}) we have $(\epsilon Z-q^{2i-d}I)X_i = X_{i-1}$ for $1 \leq i \leq d$ and $(\epsilon Z-q^{-d}I)X_0 = 0$.
By this, Lemma \ref{lem:q-weyl-linalg2} and since $\{X_i\}_{i=0}^d$ is a decomposition of $V$, we obtain that $Z,X$ satisfy (\ref{eq:uqsl2zx}).
Next we show that $X,Y$ satisfy (\ref{eq:uqsl2xy}).
By construction $Z_i = E^*_iV$ for $0 \leq i \leq d$.
By (\ref{eq:u_id_sum}) and since $\{Y_i\}_{i=0}^d$ is the inverted $\Phi$-split decomposition of $V$, we have $E_{i}V + E_{i+1}V + \cdots + E_{d}V = Y_0 + Y_1 + \cdots + Y_{d-i}$ for $0 \leq i \leq d$.
Now for $0 \leq i \leq d$ we have
\begin{eqnarray}
X_i
&=& (E_d^*V + E_{d-1}^*V + \cdots + E_{d-i}^*V) \cap (E_{i}V + E_{i+1}V + \cdots + E_{d}V) \nonumber \\
&=& (Z_d + Z_{d-1} + \cdots + Z_{d-i}) \cap (Y_0 + Y_1 + \cdots + Y_{d-i}). \label{eq:x_i_temp}
\end{eqnarray}
By (\ref{eq:uqsl2yz}) and Lemma \ref{lem:q-weyl-linalg2}, we have $(\epsilon Y-q^{2i-d}I)Z_i \subseteq Z_{i-1}$ for $1 \leq i \leq d$ and $(\epsilon Y-q^{-d}I)Z_0 = 0$.
Combining this with (\ref{eq:y-map}) and (\ref{eq:x_i_temp}), we have $(\epsilon Y-q^{2i-d}I)X_i \subseteq X_{i+1}$ for $0 \leq i \leq d-1$ and $(\epsilon Y-q^{d}I)X_d = 0$.
By this, Lemma \ref{lem:q-weyl-linalg1} and since $\{X_i\}_{i=0}^d$ is a decomposition of $V$, we obtain $X,Y$ satisfy (\ref{eq:uqsl2xy}).
We have now shown that $X,Y,Z,Z^{-1}$ satisfy the relations (\ref{eq:uqsl2zzinverse})--(\ref{eq:uqsl2zx}).
Therefore there exists a $U_q(\mathfrak{sl}_2)$-module structure on $V$ such that $x,y,z$ act on $V$ as $X,Y,Z$, respectively.
This $U_q(\mathfrak{sl}_2)$-module structure is unique since $x,y,z^{\pm 1}$ generate $U_q(\mathfrak{sl}_2)$.
\hfill $\Box$\\

\noindent
We now display the second $U_q(\mathfrak{sl}_2)$-module structure on $V$.

\begin{definition} \label{def:xyzdown-lintranfs} \rm
Let $\{X^\Downarrow_i\}_{i=0}^d$ denote the $\Phi^{\downarrow\Downarrow}$-split decomposition of $V$.
Let $\{Y^\Downarrow_i\}_{i=0}^d$ denote the inverted $\Phi^{\Downarrow}$-split decomposition of $V$.
Define the linear transformations $X^\Downarrow,Y^\Downarrow$ in ${\rm End}(V)$ as follows.
\begin{align*}
(\epsilon X^\Downarrow-q^{d-2i}I)X^\Downarrow_i &= 0 \qquad (0 \leq i \leq d),\\%\label{eq:xprime-map}\\
(\epsilon Y^\Downarrow-q^{d-2i}I)Y^\Downarrow_i &= 0 \qquad (0 \leq i \leq d).%\label{eq:yprime-map}
\end{align*} 
\end{definition}

\begin{lemma}
\label{lem:xyzprime-lintranfs-uqsl2}
There exists a unique $U_q(\mathfrak{sl}_2)$-module structure on $V$ such that the equitable generators $x,y,z$ of $U_q(\mathfrak{sl}_2)$ act on $V$ as $X^\Downarrow,Y^\Downarrow,Z$, respectively.
\end{lemma}

\noindent {\it Proof:} \rm
Recall that the Leonard system $\Phi^\Downarrow$ has dual $q$-Krawtchouk type.
Apply Lemma \ref{lem:xyz-lintranfs-uqsl2} to $\Phi^\Downarrow$.
\hfill $\Box$\\

\begin{lemma}
\label{lem:v_lde_1}
The $U_q(\mathfrak{sl}_2)$-module $V$ from Lemma \ref{lem:xyz-lintranfs-uqsl2} is isomorphic to $L(d,\epsilon)$.
\end{lemma}

\noindent
{\it Proof:} \rm
By construction $Z$ is diagonalizable on $V$.
By this and \cite[Theorem~2.9]{jantzen} we find that $V$ is a direct sum of irreducible $U_q(\mathfrak{sl}_2)$-submodules of $V$.
Observe that $Z$ is diagonalizable on each $U_q(\mathfrak{sl}_2)$-module in the sum.
By construction $\epsilon q^d$ is an eigenvalue of $Z$.
By these comments there exists a $U_q(\mathfrak{sl}_2)$-submodule $W$ of $V$ in the sum such that
$\epsilon q^d$ is an eigenvalue for $Z$ on $W$.
By Lemma \ref{lem:q_cond}(i) and Lemma \ref{lem:lnep-zxy} any irreducible $U_q(\mathfrak{sl}_2)$-module that has $\epsilon q^d$ as an eigenvalue for $Z$ has dimension at least $d+1$.
Therefore $V = W$.
By Lemma \ref{lem:fd-u-module-kfe} the $U_q(\mathfrak{sl}_2)$-module $V$ is isomorphic to $L(d,\epsilon)$.
\hfill $\Box$\\

\begin{lemma}
\label{lem:v_lde_2}
The $U_q(\mathfrak{sl}_2)$-module $V$ from Lemma \ref{lem:xyzprime-lintranfs-uqsl2} is isomorphic to $L(d,\epsilon)$.
\end{lemma}

\noindent {\it Proof:} \rm
Recall that the Leonard system $\Phi^\Downarrow$ has dual $q$-Krawtchouk type.
Apply Lemma \ref{lem:v_lde_1} to $\Phi^\Downarrow$.
\hfill $\Box$\\

\begin{lemma}
\label{lem:norm_phi_split_is_norm_y_eigen}
For the $U_q(\mathfrak{sl}_2)$-module structure on $V$ from Lemma \ref{lem:xyz-lintranfs-uqsl2},
the following coincide:
\begin{enumerate}
\item[{\rm (i)}] The inversion of a normalized $\Phi$-split basis for $V$.
\item[{\rm (ii)}] A normalized $y$-eigenbasis for the $U_q(\mathfrak{sl}_2)$-module $V$.
\end{enumerate}
\end{lemma}

\noindent {\it Proof:} \rm
Let $\{u_i\}_{i=0}^d$ denote a basis for $V$.
Recall the $\Phi$-split decomposition $\{U_i\}_{i=0}^d$ of $V$ from (\ref{eq:split_decomp_def}).
By construction $Y_i = U_{d-i}$ $(0 \leq i \leq d)$ and $Z_d = E^*_dV$.
By this and Lemma \ref{lem:char_norm_split_basis_5} the sequence $\{u_i\}_{i=0}^d$ is the inversion of a normalized $\Phi$-split basis for $V$ if and only if
$u_i \in Y_i$ $(0 \leq i \leq d)$ and $(\epsilon Z - q^{-d}I) \sum_{i=0}^d u_i = 0$.
By comparing this to (\ref{eq:lnep-yzx-y}) and the comment in Definition \ref{def:norm-y-e-basis},
the sequence $\{u_i\}_{i=0}^d$ is the inversion of a normalized $\Phi$-split basis for $V$
 if and only if
it is a normalized $Y$-eigenbasis for the $U_q(\mathfrak{sl}_2)$-module $V$.
\hfill $\Box$\\

\begin{lemma}
\label{lem:norm_phi_split_is_norm_y_eigen2}
For the $U_q(\mathfrak{sl}_2)$-module structure on $V$ from Lemma \ref{lem:xyzprime-lintranfs-uqsl2},
the following coincide:
\begin{enumerate}
\item[{\rm (i)}] The inversion of a normalized $\Phi^\Downarrow$-split basis for $V$.
\item[{\rm (ii)}] A normalized $y$-eigenbasis for the $U_q(\mathfrak{sl}_2)$-module $V$.
\end{enumerate}
\end{lemma}

\noindent {\it Proof:} \rm
Recall that the Leonard system $\Phi^\Downarrow$ has dual $q$-Krawtchouk type.
Apply Lemma \ref{lem:norm_phi_split_is_norm_y_eigen} to $\Phi^\Downarrow$.
\hfill $\Box$\\

\noindent
Recall our Leonard system $\Phi = (A;\{E_i\}_{i=0}^d;A^*;\{E^*_i\}_{i=0}^d)$.
We now show how $A,A^*$ are related to the maps $X,Y,Z$ from Definition \ref{def:xyz-lintranfs}.

\begin{lemma} \label{cor:aastar-in-xyz}
The maps $A,A^*$ can be expressed in terms of $X,Y,Z$ as follows.
\begin{align}
A & = h1 + \epsilon \kappa X + \epsilon \upsilon Y, \label{eq:a_temp_dualqkrawtchouk}\\
A^* & = h^*1 + \epsilon \kappa^* Z. \label{eq:astar_temp_dualqkrawtchouk}
\end{align}
\end{lemma}

\noindent {\it Proof:} \rm
Line (\ref{eq:astar_temp_dualqkrawtchouk}) follows from (\ref{eq:z_to_astar}).
It remains to show that (\ref{eq:a_temp_dualqkrawtchouk}) holds.
Let $\{u_i\}_{i=0}^d$ denote the inversion of a normalized $\Phi$-split basis for $V$.
By Lemma \ref{lem:norm_phi_split_is_norm_y_eigen} this is also a normalized $Y$-eigenbasis for the $U_q(\mathfrak{sl}_2)$-module $V$ from Lemma \ref{lem:xyz-lintranfs-uqsl2}.
Now compare the action of each side of (\ref{eq:a_temp_dualqkrawtchouk}) on the basis $\{u_i\}_{i=0}^d$ as follows.
Let $i$ be given.
On the left-hand side evaluate $Au_i$ using Lemma \ref{lem:lp-norm-split-basis-aastar}.
On the right-hand side evaluate $(h1 + \epsilon \kappa X + \epsilon \upsilon Y)u_i$ using (\ref{eq:lnep-yzx-y}), (\ref{eq:lnep-yzx-x}).
By comparing the result using (\ref{eq:lp-eval}), (\ref{eq:lp-split1}) and Lemma \ref{lem:diffofevals}(ii), we find that both sides are equal.
Therefore (\ref{eq:a_temp_dualqkrawtchouk}) holds.
\hfill $\Box$\\

\noindent
Next we show how $A,A^*$ are related to the maps $X^\Downarrow,Y^\Downarrow$ from Definition \ref{def:xyzdown-lintranfs} and the map $Z$ from Definition \ref{def:xyz-lintranfs}.

\begin{lemma} \label{cor:aastar-in-xyzprime}
The maps $A,A^*$ can be expressed in terms of $X^\Downarrow,Y^\Downarrow,Z$ as follows.
\begin{align*}
A & = h1 + \epsilon \kappa Y^\Downarrow + \epsilon \upsilon X^\Downarrow,\\
A^* & = h^*1 + \epsilon \kappa^* Z.
\end{align*}
\end{lemma}

\noindent {\it Proof:} \rm
Recall that the Leonard system $\Phi^\Downarrow$ has dual $q$-Krawtchouk type.
Apply Lemma \ref{cor:aastar-in-xyz} to $\Phi^\Downarrow$ and use Note \ref{note:swap-uv} to obtain the result.
\hfill $\Box$\\

\noindent
We now express $X,Y,Z$ in terms of $A,A^*$.  The expressions will involve the inverse of $A^* - h^*I$.  We take a moment to verify that the inverse exists.

\begin{lemma}
The eigenvalues of $A^*-h^*I$ are $\kappa^*q^{d-2i}$ $(0 \leq i \leq d)$.  Moreover $A^*-h^*I$ is invertible.
\end{lemma}

\noindent {\it Proof:} \rm
By construction $A^*-h^*I$ is diagonalizable with eigenvalues $\theta^*_i-h^*$ $(0 \leq i \leq d)$.
By (\ref{eq:lp-dualeval}) we have $\theta^*_i-h^* = \kappa^*q^{d-2i}$ $(0 \leq i \leq d)$
so these scalars are all nonzero.
By these comments $A^*-h^*I$ is invertible.
\hfill $\Box$

\begin{lemma}
\label{lem:xyz_aastar}
The maps $X,Y,Z$ can be expressed in terms of $A,A^*$ as follows.
\begin{align}
\epsilon X & = \frac{qB-q^{-1}B^*BB^{*-1}}{\kappa q^{-1}(q^2-q^{-2})} + \frac{\kappa^*(\kappa q^{-1}-\upsilon q)B^{*-1}}{\kappa (q+q^{-1})} \label{eq:uqsl2_x1}\\
& = \frac{qB^{*-1}BB^*-q^{-1}B}{\kappa q(q^2-q^{-2})} + \frac{\kappa^*(\kappa q-\upsilon q^{-1})B^{*-1}}{\kappa (q+q^{-1})}, \label{eq:uqsl2_x2}\\
\epsilon Y & = \frac{qB-q^{-1}B^{*-1}BB^*}{\upsilon q^{-1}(q^2-q^{-2})} + \frac{\kappa^*(\upsilon q^{-1}-\kappa q)B^{*-1}}{\upsilon (q+q^{-1})} \label{eq:uqsl2_y1}\\
& = \frac{qB^*BB^{*-1}-q^{-1}B}{\upsilon q(q^2-q^{-2})} + \frac{\kappa^*(\upsilon q-\kappa q^{-1})B^{*-1}}{\upsilon (q+q^{-1})}, \label{eq:uqsl2_y2}\\
\epsilon Z & = \kappa^{*-1}B^*, \label{eq:uqsl2_z}
\end{align}
where $B = A - hI$ and $B^* = A^* - h^*I$.
\end{lemma}

\noindent {\it Proof:} \rm
By Lemma \ref{cor:aastar-in-xyz},
\begin{eqnarray}
B = \epsilon \kappa X+ \epsilon \upsilon Y, \qquad B^* = \epsilon \kappa^* Z.	\label{eq:bbstar}
\end{eqnarray}
Thus (\ref{eq:uqsl2_z}) holds.
Next we show that (\ref{eq:uqsl2_x1}), (\ref{eq:uqsl2_x2}) hold.
To do this we claim
\begin{eqnarray}
 \kappa\kappa^* \frac{q X Z - q^{-1}Z X}{q-q^{-1}} + \upsilon \kappa^*I = \frac{qBB^*-q^{-1}B^*B}{q-q^{-1}}.
\label{eq:bbstar_x}
\end{eqnarray}
To prove the claim we evaluate $B,B^*$ using (\ref{eq:bbstar}) and simplify the result using (\ref{eq:uqsl2yz}).
By (\ref{eq:uqsl2zx}),
\begin{eqnarray}
\frac{q ZX - q^{-1}XZ}{q-q^{-1}} = I.\label{eq:temp-uq-zx}
\end{eqnarray}
Now view (\ref{eq:bbstar_x}) and (\ref{eq:temp-uq-zx}) as a system of linear equations in $X Z$ and $Z X$.  Solving the system we have
\begin{align}
X Z & = \frac{q^2BB^*-B^*B}{\kappa\kappa^*(q^2-q^{-2})} + \frac{\kappa q^{-1}- \upsilon q}{\kappa (q+q^{-1})}I, \label{eq:xzbbstar}\\
Z X & = \frac{BB^*-q^{-2}B^*B}{\kappa\kappa^*(q^2-q^{-2})} + \frac{\kappa q - \upsilon q^{-1}}{\kappa(q+q^{-1})}I. \label{eq:zxbbstar}
\end{align}
In (\ref{eq:xzbbstar}) multiply each side on the right by $Z^{-1}$ and evaluate $Z$ using (\ref{eq:uqsl2_z}) to obtain (\ref{eq:uqsl2_x1}).
In (\ref{eq:zxbbstar}) multiply each side on the left by $Z^{-1}$ and evaluate $Z$ using (\ref{eq:uqsl2_z}) to obtain (\ref{eq:uqsl2_x2}).
Next we show that (\ref{eq:uqsl2_y1}), (\ref{eq:uqsl2_y2}) hold. To do this we claim
\begin{eqnarray}
\kappa^* \upsilon  \frac{q Z Y - q^{-1}Y Z}{q-q^{-1}} + \kappa\kappa^*I = \frac{qB^*B-q^{-1}BB^*}{q-q^{-1}}.
\label{eq:bbstar_y}
\end{eqnarray}
To prove the claim we evaluate $B,B^*$ using (\ref{eq:bbstar}) and simplify the result using (\ref{eq:uqsl2zx}).
By (\ref{eq:uqsl2yz}),
\begin{eqnarray}
\frac{q YZ - q^{-1}ZY}{q-q^{-1}} = I.\label{eq:temp-uq-yz}
\end{eqnarray}
Now view (\ref{eq:bbstar_y}) and (\ref{eq:temp-uq-yz}) as a system of linear equations in $Z Y$ and $Y Z$.  Solving the system we have
\begin{align}
Z Y & = \frac{q^{2}B^*B-BB^*}{\upsilon\kappa^*(q^2-q^{-2})} + \frac{\upsilon q^{-1} - \kappa q}{\upsilon (q+q^{-1})}I, \label{eq:zybbstar}\\
Y Z & = \frac{B^*B-q^{-2}BB^*}{\upsilon \kappa^*(q^2-q^{-2})} + \frac{\upsilon q -  \kappa q^{-1}}{\upsilon (q+q^{-1})}I. \label{eq:yzbbstar}
\end{align}
In (\ref{eq:zybbstar}) multiply each side on the left by $Z^{-1}$ and evaluate $Z$ using (\ref{eq:uqsl2_z}) to obtain (\ref{eq:uqsl2_y1}).
In (\ref{eq:yzbbstar}) multiply each side on the right by $Z^{-1}$ and evaluate $Z$ using (\ref{eq:uqsl2_z}) to obtain (\ref{eq:uqsl2_y2}).
\hfill $\Box$\\

\noindent
Next we express $X^\Downarrow,Y^\Downarrow$ in terms of $A,A^*$

\begin{lemma}
\label{lem:xyz_aastar2}
The maps $X^\Downarrow,Y^\Downarrow$ can be expressed in terms of $A,A^*$ as follows.
\begin{align}
\epsilon X^\Downarrow & = \frac{qB-q^{-1}B^*BB^{*-1}}{\upsilon q^{-1}(q^2-q^{-2})} + \frac{\kappa^*(\upsilon q^{-1}-\kappa q)B^{*-1}}{\upsilon (q+q^{-1})} \label{eq:uqsl2_x1prime}\\
& = \frac{qB^{*-1}BB^*-q^{-1}B}{\upsilon q(q^2-q^{-2})} + \frac{\kappa^*(\upsilon q - \kappa q^{-1})B^{*-1}}{\upsilon (q+q^{-1})}, \nonumber\\%\label{eq:uqsl2_x2prime}\\
\epsilon Y^\Downarrow & = \frac{qB-q^{-1}B^{*-1}BB^*}{\kappa q^{-1}(q^2-q^{-2})} + \frac{\kappa^*(\kappa q^{-1}-\upsilon q)B^{*-1}}{\kappa (q+q^{-1})} \label{eq:uqsl2_y1prime}\\
& = \frac{qB^*BB^{*-1}-q^{-1}B}{\kappa q(q^2-q^{-2})} + \frac{\kappa^*(\kappa q - \upsilon q^{-1})B^{*-1}}{\kappa (q+q^{-1})}, \nonumber%\label{eq:uqsl2_y2prime}
\end{align}
where $B = A - hI$ and $B^* = A^* - h^*I$.
\end{lemma}

\noindent {\it Proof:} \rm
Recall that the Leonard system $\Phi^\Downarrow$ has dual $q$-Krawtchouk type.
Apply Lemma \ref{lem:xyz_aastar} to $\Phi^\Downarrow$ and use Note \ref{note:swap-uv} to obtain the result.
\hfill $\Box$\\

\noindent
We comment on how the two $U_q(\mathfrak{sl}_2)$-module structures are related.

\begin{lemma}
The $U_q(\mathfrak{sl}_2)$-module structure from Lemma \ref{lem:xyz-lintranfs-uqsl2} and the $U_q(\mathfrak{sl}_2)$-module structure from Lemma \ref{lem:xyzprime-lintranfs-uqsl2} are related as follows.
\begin{align}
X^\Downarrow &= \kappa \upsilon^{-1}X + (1-\kappa \upsilon^{-1})Z^{-1},\label{eq:xxddown}\\
Y^\Downarrow &= \upsilon \kappa^{-1}Y + (1-\upsilon \kappa^{-1})Z^{-1}.\label{eq:yyddown}
\end{align}
\end{lemma}

\noindent {\it Proof:} \rm
To verify (\ref{eq:xxddown}) evaluate $X, Z, X^\Downarrow$ using (\ref{eq:uqsl2_x1}), (\ref{eq:uqsl2_z}) and (\ref{eq:uqsl2_x1prime}), respectively, and simplify the result.
To verify (\ref{eq:yyddown}) evaluate $Y, Z, Y^\Downarrow$ using (\ref{eq:uqsl2_y1}), (\ref{eq:uqsl2_z}) and (\ref{eq:uqsl2_y1prime}), respectively, and simplify the result.
\hfill $\Box$\\

\noindent
We summarize this section with the following theorems.

\begin{theorem}
\label{thm:uqsl2_on_v1}
Let $\Phi$ denote a Leonard system on $V$ as in Definition \ref{def:ls}.  Assume $\Phi$ has dual $q$-Krawtchouk type.
Let $X,Y,Z$ denote the corresponding elements from Definition \ref{def:xyz-lintranfs}.
Then the following {\rm (i)--(iii)} hold.
\begin{enumerate}
\item[\rm (i)] There exists a unique $U_q(\mathfrak{sl}_2)$-module structure on $V$ such that the equitable generators $x,y,z$ of $U_q(\mathfrak{sl}_2)$ act on $V$ as the maps $X,Y,Z$, respectively.
\item[\rm (ii)] The $U_q(\mathfrak{sl}_2)$-module $V$ is isomorphic to $L(d,\epsilon)$.
\item[\rm (iii)] The inversion of a normalized $\Phi$-split basis for $V$ coincides with a normalized $y$-eigenbasis for the $U_q(\mathfrak{sl}_2)$-module $V$.
\end{enumerate}
\end{theorem}

\noindent {\it Proof:} \rm
Combine Lemma \ref{lem:xyz-lintranfs-uqsl2} and Lemma \ref{lem:norm_phi_split_is_norm_y_eigen}.
\hfill $\Box$\\

\begin{theorem}
\label{thm:uqsl2_on_v2}
Let $\Phi$ denote a Leonard system on $V$ as in Definition \ref{def:ls}.  Assume $\Phi$ has dual $q$-Krawtchouk type.
Let $X^\Downarrow,Y^\Downarrow$ denote the corresponding elements from Definition \ref{def:xyzdown-lintranfs} and let $Z$ denote the corresponding element from Definition \ref{def:xyz-lintranfs}.
Then the following {\rm (i)--(iii)} hold.
\begin{enumerate}
\item[\rm (i)] There exists a unique $U_q(\mathfrak{sl}_2)$-module structure on $V$ such that the equitable generators $x,y,z$ of $U_q(\mathfrak{sl}_2)$ act on $V$ as the maps $X^\Downarrow,Y^\Downarrow,Z$, respectively.
\item[\rm (ii)] The $U_q(\mathfrak{sl}_2)$-module $V$ is isomorphic to $L(d,\epsilon)$.
\item[\rm (iii)] The inversion of a normalized $\Phi^{\Downarrow}$-split basis for $V$ coincides with a normalized $y$-eigenbasis for the $U_q(\mathfrak{sl}_2)$-module $V$.
\end{enumerate}
\end{theorem}

\noindent {\it Proof:} \rm
Combine Lemma \ref{lem:xyzprime-lintranfs-uqsl2} and Lemma \ref{lem:norm_phi_split_is_norm_y_eigen2}.
\hfill $\Box$\\

\begin{theorem}
\label{thm:uqsl2ondualqkrawtchouk}
Let $\Phi$ denote a Leonard system on $V$ as in Definition \ref{def:ls}.  Assume $\Phi$ has dual $q$-Krawtchouk type.
Let $h,h^*,\kappa,\kappa^*,\upsilon$ denote the corresponding parameters from Definition \ref{def:dualqkrawtchouk2}.
Then there exists a unique $U_q(\mathfrak{sl}_2)$-module structure on $V$ such that on $V$,
\begin{align*}
A & = h1 + \epsilon \kappa x + \epsilon \upsilon y,\\%\label{eq:a_dualqkrawtchouk}\\
A^* & = h^*1 + \epsilon \kappa^*z,%\label{eq:astar_dualqkrawtchouk}
\end{align*}
where $x,y,z$ are the equitable generators for $U_q(\mathfrak{sl}_2)$.
This $U_q(\mathfrak{sl}_2)$-module structure coincides with the $U_q(\mathfrak{sl}_2)$-module structure from Theorem \ref{thm:uqsl2_on_v1}.
\end{theorem}

\noindent {\it Proof:} \rm
The existence follows from Lemma \ref{lem:xyz-lintranfs-uqsl2} and Lemma \ref{cor:aastar-in-xyz}.
The uniqueness follows from Lemma \ref{lem:xyz_aastar} and since $x,y,z^{\pm 1}$ generate $U_q(\mathfrak{sl}_2)$.
\hfill $\Box$\\

\begin{theorem}
\label{thm:uqsl2ondualqkrawtchouk2}
Let $\Phi$ denote a Leonard system on $V$ as in Definition \ref{def:ls}.  Assume $\Phi$ has dual $q$-Krawtchouk type.
Let $h,h^*,\kappa,\kappa^*,\upsilon$ denote the corresponding parameters from Definition \ref{def:dualqkrawtchouk2}.
Then there exists a unique $U_q(\mathfrak{sl}_2)$-module structure on $V$ such that on $V$,
\begin{align*}
A & = h1 + \epsilon \kappa y + \epsilon \upsilon x,\\%\label{eq:a_dualqkrawtchouk2}\\
A^* & = h^*1 + \epsilon \kappa^*z,%\label{eq:astar_dualqkrawtchouk2}
\end{align*}
where $x,y,z$ are the equitable generators for $U_q(\mathfrak{sl}_2)$.
This $U_q(\mathfrak{sl}_2)$-module structure coincides with the $U_q(\mathfrak{sl}_2)$-module structure from Theorem \ref{thm:uqsl2_on_v2}.
\end{theorem}

\noindent {\it Proof:} \rm
Recall that the Leonard system $\Phi^\Downarrow$ has dual $q$-Krawtchouk type.
Apply Theorem \ref{thm:uqsl2ondualqkrawtchouk} to $\Phi^\Downarrow$ and use Note \ref{note:swap-uv} to obtain the result.
\hfill $\Box$\\

\section{Distance-regular graphs; preliminaries}\label{sec:drg}

We now turn our attention to distance-regular graphs.
After a brief review of the basic definitions
we recall the subconstituent algebra and the $Q$-polynomial structure.
For more information we refer the reader to \cite{bannai,bcn,godsil,terwSub1}.

\medskip
\noindent
Let $X$ denote a nonempty  finite  set.
Let ${\rm Mat}_X(\C)$ 
denote the $\C$-algebra
consisting of the matrices with entries in $\C$,
and rows and columns indexed by $X$.
Let $V=\C^X$ denote the vector space over $\C$
consisting of the column vectors with entries in $\C$
and rows indexed by $X$.
Observe that
${\rm Mat}_X(\C)$ 
acts on $V$ by left multiplication.
We call $V$ the {\it standard module} of ${\rm Mat}_X(\C)$.
We endow $V$ with the Hermitean inner product $\langle \, , \, \rangle$ 
that satisfies
$\langle u,v \rangle = u^t\overline{v}$ for 
$u,v \in V$,
where $t$ denotes transpose and $\overline{\phantom{v}}$
denotes complex conjugation.
For all $y \in X,$ let $\hat{y}$ denote the element
of $V$ with a 1 in the $y$ coordinate and 0 in all other coordinates.
Observe that $\{\hat{y}\;|\;y \in X\}$ is an orthonormal basis for $V.$

\medskip
\noindent
Let $\Gamma = (X,R)$ denote a finite, undirected, connected graph,
without loops or multiple edges, with vertex set $X$,
edge set $R$,
path-length distance function $\partial$,
and diameter $D := {\rm max}\{\partial(x,y)|x,y \in X\}$.  
For $x \in X$ and an integer $i$ $(0 \leq i \leq D)$ let $\Gamma_i(x) = \{y \in X | \partial(x,y) = i\}$.
We abbreviate $\Gamma(x) = \Gamma_1(x)$.
For an integer $k\geq 0$ we say $\Gamma$ is {\it regular with
valency $k$} whenever $|\Gamma(x)| = k$ for every $x \in X$.
We say $\Gamma$ is {\it distance-regular}
whenever for all integers $h,i,j\;(0 \le h,i,j \le D)$ 
and for all
vertices $x,y \in X$ with $\partial(x,y)=h,$ the number
\begin{eqnarray*}
p_{ij}^h = |\Gamma_i(x) \cap \Gamma_j(y)|
\end{eqnarray*}
is independent of $x$ and $y$.
The constants $p_{ij}^h$ are called
the {\it intersection numbers} of $\Gamma.$ 
We abbreviate $c_i=p^i_{1,i-1}$ $(1 \leq i \leq D)$,
$b_i=p^i_{1,i+1}$ $(0 \leq i \leq D-1)$,
$a_i=p^i_{1i}$ $(0 \leq i \leq D)$,
and $c_0 = 0$, $b_D = 0$.

\medskip
\noindent
For the rest of the paper assume  $\Gamma$  
is  distance-regular with $D\geq 3$.
Observe that $\Gamma$ is regular with valency $k=b_0$. Moreover
 $k=c_i+a_i+b_i$ for $0 \leq i \leq D$.
By the triangle inequality, for $0 \leq h,i,j\leq D$ we have
$p^h_{ij}= 0$
(resp. 
$p^h_{ij}\not= 0$) whenever one of $h,i,j$ is greater than
(resp. equal to) the sum of the other two.
In particular $c_i \neq 0$ for $1 \leq i \leq D$ and $b_i \neq 0$ for $0 \leq i \leq D-1$.

\medskip
\noindent 
We recall the Bose-Mesner algebra of $\Gamma.$ 
For 
$0 \le i \le D$ let $A_i$ denote the matrix in ${\rm Mat}_X(\C)$ with
$(x,y)$-entry
\begin{eqnarray*}
(A_i)_{xy} =
\begin{cases}
1, & \mbox{if } \partial(x,y)=i\\
0, & \mbox{if } \partial(x,y) \ne i
\end{cases}
\qquad (x,y \in X).
\end{eqnarray*}
We call $A_i$ the $i$th {\it distance matrix} of $\Gamma.$
We abbreviate $A=A_1$ and call this  the {\it adjacency
matrix} of $\Gamma.$
Observe that
(ai) $A_0 = I$;
 (aii)
$J = \sum_{i=0}^D A_i$;
(aiii)
$\overline{A_i} = A_i \;(0 \le i \le D)$;
(aiv) $A_i^t = A_i  \;(0 \le i \le D)$;
(av) $A_iA_j = \sum_{h=0}^D p_{ij}^h A_h \;( 0 \le i,j \le D)
$,
where $I$ (resp. $J$) denotes the identity matrix 
(resp. all 1's matrix) in 
 ${\rm Mat}_X(\C)$.
 Using these facts  we find
 $\lbrace A_i\rbrace_{i=0}^D$
is a basis for a commutative subalgebra $M$ of 
${\rm Mat}_X(\C)$.
We call $M$ the
{\it Bose-Mesner algebra} of $\Gamma$.
By \cite[p.~190]{bannai} $A$ generates $M$.
By \cite[p.~45]{bcn} $M$ has a second basis 
$\lbrace E_i\rbrace_{i=0}^D$ such that
(ei) $E_0 = |X|^{-1}J$;
(eii) $I = \sum_{i=0}^D E_i$;
(eiii) $\overline{E_i} = E_i \;(0 \le i \le D)$;
(eiv) $E_i^t =E_i  \;(0 \le i \le D)$;
(ev) $E_iE_j =\delta_{ij}E_i  \;(0 \le i,j \le D)$.
We call $\lbrace E_i\rbrace_{i=0}^D$  the {\it primitive idempotents}
of $\Gamma$.  

\medskip
\noindent
We  recall the eigenvalues
of  $\Gamma $.
Since $\lbrace E_i\rbrace_{i=0}^D$ form a basis for  
$M$, there exist complex scalars 
$\lbrace\theta_i\rbrace_{i=0}^D$
such that
$A = \sum_{i=0}^D \theta_iE_i$.
Combining this with (ev) we find
$AE_i = E_iA =  \theta_iE_i$ for $0 \leq i \leq D$.
We call $\theta_i$  the {\it eigenvalue}
of $\Gamma$ associated with $E_i$.
By \cite[p.~197]{bannai} the 
scalars $\lbrace \theta_i\rbrace_{i=0}^D$ are
in $\R.$
The
$\lbrace\theta_i\rbrace_{i=0}^D$ are mutually distinct 
since $A$ generates $M$.
By (ei) we have $\theta_0 = k$.
By (eii)--(ev),
\begin{eqnarray*}
V = E_0V+E_1V+ \cdots +E_DV \qquad {\rm (orthogonal\ direct\ sum}).
\end{eqnarray*}
For $0 \le i \le D$ the space $E_iV$ is the  eigenspace of $A$ associated 
with $\theta_i$.
Let $m_i$ denote the rank of $E_i$
and note that $m_i$ is the dimension of $E_iV$.
We call $m_i$ the {\em multiplicity} of $E_i$ (or $\theta_i$).

\medskip
\noindent 
We recall the Krein parameters of $\Gamma$.
Let $\circ $ denote the entrywise product in
${\rm Mat}_X(\C)$.
Observe that
$A_i\circ A_j= \delta_{ij}A_i$ for $0 \leq i,j\leq D$,
so
$M$ is closed under
$\circ$. Thus there exist complex scalars
$q^h_{ij}$  $(0 \leq h,i,j\leq D)$ such
that
$$
E_i\circ E_j = |X|^{-1}\sum_{h=0}^D q^h_{ij}E_h
\qquad (0 \leq i,j\leq D).
$$
By \cite[p.~170]{biggs}, 
$q^h_{ij}$ is real and nonnegative  for $0 \leq h,i,j\leq D$.
The $q^h_{ij}$ are called the {\it Krein parameters} of $\Gamma$.

\medskip
\noindent
We  recall the dual Bose-Mesner algebra of $\Gamma.$
For the rest of this section we fix
a vertex $x \in X.$ We view $x$ as a ``base vertex''.
For 
$ 0 \le i \le D$ let $E_i^*=E_i^*(x)$ denote the diagonal
matrix in ${\rm Mat}_X(\C)$ with $(y,y)$-entry
\begin{eqnarray}\label{DEFDEI}
(E_i^*)_{yy} =
\begin{cases}
1, & \mbox{if } \partial(x,y)=i\\
0, & \mbox{if } \partial(x,y) \ne i
\end{cases}
\qquad (y \in X).
\end{eqnarray}
We call $E_i^*$ the  $i$th {\it dual idempotent} of $\Gamma$
 with respect to $x$ \cite[p.~378]{terwSub1}.
Observe that
(esi) $I = \sum_{i=0}^D E_i^*$;
(esii) $\overline{E_i^*} = E_i^*$ $(0 \le i \le D)$;
(esiii) $E_i^{*t} = E_i^*$ $(0 \le i \le D)$;
(esiv) $E_i^*E_j^* = \delta_{ij}E_i^* $ $(0 \le i,j \le D)$.
By these facts 
$\lbrace E_i^*\rbrace_{i=0}^D$ form a 
basis for a commutative subalgebra
$M^*=M^*(x)$ of 
${\rm Mat}_X(\C).$ 
We call 
$M^*$ the {\it dual Bose-Mesner algebra} of
$\Gamma$ with respect to $x$ \cite[p.~378]{terwSub1}.
For $0 \leq i \leq D$ let $A^*_i = A^*_i(x)$ denote the diagonal
matrix in 
 ${\rm Mat}_X(\C)$
with $(y,y)$-entry
$(A^*_i)_{yy}=\vert X \vert (E_i)_{xy}$ for $y \in X$.
Then $\lbrace A^*_i\rbrace_{i=0}^D$ is a basis for $M^*$ 
\cite[p.~379]{terwSub1}.
Moreover
(asi) $A^*_0 = I$;
(asii)
$\overline{A^*_i} = A^*_i \;(0 \le i \le D)$;
(asiii) $A^{*t}_i = A^*_i  \;(0 \le i \le D)$;
(asiv) $A^*_iA^*_j = \sum_{h=0}^D q_{ij}^h A^*_h \;( 0 \le i,j \le D)
$
\cite[p.~379]{terwSub1}.
We call 
 $\lbrace A^*_i\rbrace_{i=0}^D$
the {\it dual distance matrices} of $\Gamma$ with respect to $x$.

\medskip
\noindent 
We recall the subconstituents of $\Gamma $.
From
(\ref{DEFDEI}) we find
\begin{equation}\label{DEIV}
E_i^*V = \mbox{span}\{\hat{y}|y \in X, \partial(x,y)=i\}
\qquad (0 \le i \le D).
\end{equation}
We call $E_i^*V$ the $i$th {\it subconstituent} of $\Gamma$
with respect to $x$.
By 
(\ref{DEIV})  and since
 $\{\hat{y}\;|\;y \in X\}$ is an orthonormal basis for $V$
 we find
\begin{eqnarray*}
\label{vsub}
V = E_0^*V+E_1^*V+ \cdots +E_D^*V \qquad 
{\rm (orthogonal\ direct\ sum}).
\end{eqnarray*}

\medskip
\noindent
We recall the subconstituent algebra of $\Gamma $.
Let $T=T(x)$ denote the subalgebra of ${\rm Mat}_X(\C)$ generated by 
$M$ and $M^*$. 
We call $T$ the {\it subconstituent algebra} 
(or {\it Terwilliger algebra}) of $\Gamma$ 
 with respect to $x$ \cite[Definition~3.3]{terwSub1}.
Observe that $T$ has finite dimension. Moreover $T$ is 
semisimple since it
is closed under the conjugate transpose map
\cite[p.~157]{CR}.
%We note that $A,A^*$ together generate $T$.
By \cite[Lemma~3.2]{terwSub1}
the following are relations in $T$.
\begin{eqnarray}
&&E^*_hA_iE^*_j = 0 \quad \mbox{iff} \quad p^h_{ij}=0,
\qquad (0 \leq h,i,j \leq D),
\label{eq:triple2}
\\
&&E_hA^*_iE_j = 0 \quad \mbox{iff} \quad q^h_{ij}=0, \qquad
(0 \leq h,i,j \leq D).
\label{eq:triple1}
\end{eqnarray}

\medskip
\noindent
We recall the $Q$-polynomial property.
The graph $\Gamma$ is said to be {\it $Q$-polynomial}
(with respect to the given ordering
$\lbrace E_i\rbrace_{i=0}^D$
of the primitive idempotents)
whenever for $0 \leq h,i,j\leq D$, 
$q^h_{ij}= 0$
(resp. 
$q^h_{ij}\not= 0$) whenever one of $h,i,j$ is greater than
(resp. equal to) the sum of the other two
\cite[p.~235]{bcn}.

\medskip
\noindent
For the rest of this section assume $\Gamma$ is $Q$-polynomial
with respect to $\lbrace E_i\rbrace_{i=0}^D$.
We abbreviate $c^*_i=q^i_{1,i-1}$ $(1 \leq i \leq D)$,
$b^*_i=q^i_{1,i+1}$ $(0 \leq i \leq D-1)$,
$a^*_i=q^i_{1i}$ $(0 \leq i \leq D)$,
and $c^*_0 = 0$, $b^*_D = 0$.
We call the sequence
$\lbrace \theta_i\rbrace_{i=0}^D$ the {\it eigenvalue sequence}
for this $Q$-polynomial structure.
We abbreviate 
$A^*=A^*_1$ 
and call this the {\it dual adjacency matrix} of $\Gamma$ with
respect to $x$.
The matrix $A^*$ generates $M^*$ \cite[Lemma~3.11]{terwSub1}.
Therefore $A, A^*$ together generate $T$.

\medskip
\noindent By (\ref{eq:triple2}), (\ref{eq:triple1}) we have
\begin{eqnarray}
&&E^*_iAE^*_j = 0 \quad \mbox{if} \quad |i-j| > 1,
\label{eq:q-triple2}
\\
&&E_iA^*E_j = 0 \quad \mbox{if} \quad |i-j| > 1
\nonumber%\label{eq:q-triple1}
\end{eqnarray}
for $0 \leq i, j \leq D$.

\medskip
\noindent We recall the dual eigenvalues of $\Gamma$.
Since $\lbrace E^*_i\rbrace_{i=0}^D$ form a basis for  
$M^*$ there exist complex scalars $\lbrace \theta^*_i\rbrace_{i=0}^D$
such that
$A^* = \sum_{i=0}^D \theta^*_iE^*_i$.
Combining this with (esiv) we find
$A^*E^*_i = E^*_iA^* =  \theta^*_iE^*_i$ for $0 \leq i \leq D$.
By \cite[Lemma~3.11]{terwSub1} the 
scalars $\lbrace \theta^*_i\rbrace_{i=0}^D$ are in $\R$. 
The scalars $\lbrace \theta^*_i\rbrace_{i=0}^D$ are mutually
distinct 
since $A^*$ generates $M^*$. We call $\theta^*_i$ the {\it dual eigenvalue}
of $\Gamma$ associated with $E^*_i$ $(0 \leq i\leq D)$.
We call the sequence $\lbrace \theta^*_i\rbrace_{i=0}^D$ the
{\it dual eigenvalue sequence} for the given $Q$-polynomial structure.
By \cite[Lemma~2.21(ii)]{bcn}, $\theta^*_0 = m_1$.
For $0 \leq i \leq D$ the space $E^*_iV$ is the eigenspace
of $A^*$ associated with $\theta^*_i$.

\medskip
\noindent
Recall the tridiagonal relations from Lemma \ref{lem:td-relns}.

\begin{lemma} {\rm\cite[Lemma~5.4]{terwSub3}}
\label{lem:drg-td-relns}
There exist unique scalars $\beta, \gamma, \gamma^*, \varrho, \varrho^*$ in $\C$ such that both
\begin{align}
\left[ A, A^2A^* - \beta AA^*A + A^*A^2 - \gamma (AA^* + A^*A) - \varrho A^* \right] & = 0, \label{eq:drg-td-reln1}\\
\left[ A^*, A^{*2}A - \beta A^*AA^* + AA^{*2} - \gamma^* (A^*A + AA^*) - \varrho^* A \right] & = 0. \nonumber%\label{eq:drg-td-reln2}
\end{align}
\end{lemma}

\begin{lemma} {\rm\cite[Lemma~5.4]{terwSub3}}
\label{lem:drgscalars}
The scalars $\beta,\gamma,\gamma^*,\varrho,\varrho^*$ from Lemma \ref{lem:drg-td-relns} satisfy the following {\rm (i)--(v)}.
\begin{enumerate}
\item[\rm (i)] The expressions
\begin{eqnarray*}
\frac{\theta_{i-2}-\theta_{i+1}}{\theta_{i-1}-\theta_{i}}, \qquad \frac{\theta_{i-2}^*-\theta_{i+1}^*}{\theta_{i-1}^*-\theta_{i}^*}
\end{eqnarray*}
are both equal to $\beta + 1$ for $2 \leq i \leq D-1$.
\item[\rm (ii)] $\gamma = \theta_{i-1} - \beta \theta_{i} + \theta_{i+1} \qquad (1 \leq i \leq D-1)$.
\item[\rm (iii)] $\gamma^* = \theta_{i-1}^* - \beta \theta_{i}^* + \theta_{i+1}^* \qquad (1 \leq i \leq D-1)$.
\item[\rm (iv)] $\varrho = \theta_{i-1}^2 - \beta \theta_{i-1} \theta_{i} + \theta_{i}^2 - \gamma(\theta_{i-1}+\theta_{i}) \qquad (1 \leq i \leq D)$.
\item[\rm (v)] $\varrho^* = \theta_{i-1}^{*2} - \beta \theta_{i-1}^* \theta_{i}^* + \theta_{i}^{*2} - \gamma^*(\theta_{i-1}^*+\theta_{i}^*) \qquad (1 \leq i \leq D)$.
\end{enumerate}
\end{lemma}

\medskip
\noindent We recall the $T$-modules.
By a {\it T-module}
we mean a subspace $W \subseteq V$ such that $BW \subseteq W$
for all $B \in T.$ 
\noindent
Let $W$ denote a $T$-module and let 
$W'$ denote a  
$T$-module contained in $W$.
Then the orthogonal complement of $W'$ in $W$ is a $T$-module 
\cite[p.~802]{go2}.
It follows that each $T$-module
is an orthogonal direct sum of irreducible $T$-modules.
In particular $V$ is an orthogonal direct sum of irreducible $T$-modules.

\medskip
\noindent 
Let $W$ denote an irreducible $T$-module.
Observe that $W$ is the direct sum of the nonzero spaces among
$E^*_0W,\ldots, E^*_DW$. Similarly
$W$ is the direct sum 
 of the nonzero spaces among
$E_0W,\ldots,$ $ E_DW$.
By the {\it endpoint} of $W$ we mean
$\mbox{min}\lbrace i |0\leq i \leq D, \; E^*_iW\not=0\rbrace $.
By the {\it diameter} of $W$ we mean
$ |\lbrace i | 0 \leq i \leq D,\; E^*_iW\not=0 \rbrace |-1 $.
By the {\it dual endpoint} of $W$ we mean
$\mbox{min}\lbrace i |0\leq i \leq D, \; E_iW\not=0\rbrace $.
By
the {\it dual diameter} of $W$ we mean
$ |\lbrace i | 0 \leq i \leq D,\; E_iW\not=0 \rbrace |-1 $.
It turns out that the
diameter of $W$ is  equal to the dual diameter of
$W$
\cite[Corollary~3.3]{aap1}.
By \cite[Lemma~3.9, Lemma~3.12]{terwSub1}
${\rm dim} \,E^*_iW \leq 1$ for $0 \leq i \leq D$ if and only if
${\rm dim} \,E_iW \leq 1$ for $0 \leq i \leq D$; in this case
$W$ is called {\it thin}.
$\Gamma$ is called {\it thin} (with respect to $x$) whenever all of its irreducible $T$-modules are thin.

\begin{lemma}
{\rm \cite[Lemma~3.4, Lemma~3.9, Lemma~3.12]{terwSub1}}
\label{lem:basic}
Let $W$ denote an irreducible $T$-module.
Let $r, t, d$ denote the endpoint, dual endpoint, and diameter of $W$, respectively.
Then $r,t,d$ are nonnegative integers such that $r+d\leq D$ and
$t+d\leq D$. Moreover the following {\rm (i)--(iv)} hold.
\begin{enumerate}
\item[\rm (i)]
$E^*_iW \not=0$ if and only if $r \leq i \leq r+d$, 
$ \quad (0 \leq i \leq D)$.
\item[\rm (ii)]
$W = \sum_{h=0}^{d} E^*_{r+h}W \qquad ({\rm orthogonal \ direct \ sum}). $
\item[\rm (iii)]
$E_iW \not=0$ if and only if $t \leq i \leq t+d$,
$ \quad (0 \leq i \leq D)$.
\item[\rm (iv)]
$W = \sum_{h=0}^{d} E_{t+h}W \qquad ({\rm orthogonal \ direct \ sum}). $
\end{enumerate}
\end{lemma}

\begin{lemma}
\label{lem:stand-basis}
{\rm \cite[Lemma~4.1, Lemma~8.7]{cerzo}}
Let $W$ denote a thin irreducible $T$-module.
Let $r, t, d$ denote the endpoint, dual endpoint, and diameter of $W$, respectively.
For any nonzero $u \in E_tW$, the sequence $\{E_{r+i}^*u\}_{i=0}^d$ is a basis for $W$.
Consider the matrices in ${\rm Mat}_{d+1}(\C)$ that represent $A$ and $A^*$ with respect to this basis.
The matrix representing $A$ is irreducible tridiagonal with constant row sum $\theta_t$, and the matrix representing $A^*$ is ${\rm diag}(\theta^*_{r}, \theta^*_{r+1}, \ldots, \theta^*_{r+d})$.
\end{lemma}

\begin{definition} \rm
\label{def:stand-basis}
\cite[Definition~8.2]{cerzo}
Let $W$ denote a thin irreducible $T$-module.
Let $r, t, d$ denote the endpoint, dual endpoint, and diameter of $W$, respectively.
By a {\em standard basis} for $W$ we mean a sequence $\{E_{r+i}^*u\}_{i=0}^d$ where $u$ is a nonzero vector in $E_tW$.
The matrix in ${\rm Mat}_{d+1}(\C)$ that represents $A$ with respect to a standard basis will be denoted
\begin{eqnarray}
\left(
\begin{array}{cccccc}
a_0(W) & b_0(W) & & & & {\bf 0}\\
c_1(W) & a_1(W) & b_1(W) & & &\\
& c_2(W) & \cdot & \cdot & &\\
& & \cdot & \cdot & \cdot &\\
& & & \cdot & \cdot & b_{d-1}(W)\\
{\bf 0} & & & & c_d(W) & a_d(W)\\
\end{array}
\right).
\label{eq:int_matrix_a}
\end{eqnarray}
We call the matrix (\ref{eq:int_matrix_a}) the {\em intersection matrix} of $W$.
\end{definition}

\begin{lemma}
{\rm \cite[Lemma~4.2,~Lemma 8.8]{cerzo}}
\label{lem:dual-stand-basis}
Let $W$ denote a thin irreducible $T$-module.
Let $r, t, d$ denote the endpoint, dual endpoint, and diameter of $W$, respectively.
For any nonzero $v \in E_r^*W$, the sequence $\{E_{t+i}v\}_{i=0}^d$ is a basis for $W$.
Consider the matrices in ${\rm Mat}_{d+1}(\C)$ that represent $A$ and $A^*$ with respect to this basis.
The matrix representing $A^*$ is irreducible tridiagonal with constant row sum $\theta^*_r$, and the matrix representing $A$ is ${\rm diag}(\theta_{t}, \theta_{t+1}, \ldots, \theta_{t+d})$.
\end{lemma}

\begin{definition} \rm
\label{def:dual-stand-basis}
\cite[Definition~8.2]{cerzo}
Let $W$ denote a thin irreducible $T$-module.
Let $r, t, d$ denote the endpoint, dual endpoint, and diameter of $W$, respectively.
By a {\em dual standard basis} for $W$ we mean a sequence $\{E_{t+i}v\}_{i=0}^d$ where $v$ is a nonzero vector in $E_r^*W$.
The matrix in ${\rm Mat}_{d+1}(\C)$ that represents $A^*$ with respect to a dual standard basis will be denoted
\begin{eqnarray}
\left(
\begin{array}{cccccc}
a_0^*(W) & b_0^*(W) & & & & {\bf 0}\\
c_1^*(W) & a_1^*(W) & b_1^*(W) & & &\\
& c_2^*(W) & \cdot & \cdot & &\\
& & \cdot & \cdot & \cdot &\\
& & & \cdot & \cdot & b_{d-1}^*(W)\\
{\bf 0} & & & & c_d^*(W) & a_d^*(W)\\
\end{array}
\right).
\label{eq:int_matrix_astar}
\end{eqnarray}
We call the matrix (\ref{eq:int_matrix_astar}) the {\em dual intersection matrix} of $W$.
\end{definition}

\begin{lemma}
\label{lem:irred-t-module-lp}
Let $W$ denote a thin irreducible $T$-module.
Let $r, t, d$ denote the endpoint, dual endpoint, and diameter of $W$, respectively.
Then the sequence
\begin{eqnarray*}
\Phi = (A|_W;\{E_{t+i}|_W\}_{i=0}^d;A^*|_W;\{E^*_{r+i}|_W\}_{i=0}^d)
\end{eqnarray*}
is a Leonard system on $W$.
The intersection matrix of $\Phi$ from Definition \ref{def:int_mat} coincides with the intersection matrix of $W$ from {\rm (\ref{eq:int_matrix_a})}.
The dual intersection matrix of $\Phi$ from Definition \ref{def:dual_int_mat} coincides with the dual intersection matrix of $W$ from {\rm (\ref{eq:int_matrix_astar})}.
\end{lemma}

\noindent {\it Proof:} \rm
The first assertion follows from Lemma \ref{lem:stand-basis} and Lemma \ref{lem:dual-stand-basis}.
The intersection matrix of $\Phi$ coincides with (\ref{eq:int_matrix_a}) because a $\Phi$-standard basis for $W$ is a standard basis for the $T$-module $W$.
The dual intersection matrix of $\Phi$ coincides with (\ref{eq:int_matrix_astar}) because a $\Phi^*$-standard basis for $W$ is a dual standard basis for the $T$-module $W$.
\hfill $\Box$\\

\begin{lemma}
\label{lem:w-aw}
Let $W$ denote a thin irreducible $T$-module.
Then there exist scalars $\omega = \omega(W), \eta = \eta(W), \eta^* = \eta^*(W)$ in $\C$ such that both
\begin{align}
A^2A^* - \beta AA^*A + A^*A^2 - \gamma (AA^*+A^*A) - \varrho A^* &= \gamma^*A^2 + \omega A + \eta I,\label{eq:drg-aw-rel1}\\
A^{*2}A - \beta A^*AA^* + AA^{*2} - \gamma^* (A^*A+AA^*) - \varrho^* A &= \gamma A^{*2} + \omega A^* + \eta^* I\label{eq:drg-aw-rel2}
\end{align}
on $W$.
In the above equations, the scalars $\beta, \gamma, \gamma^*, \varrho, \varrho^*$ are from Lemma \ref{lem:drg-td-relns}.
\end{lemma}

\noindent {\it Proof:}
Apply Lemma \ref{lem:lpaw} to the Leonard system from Lemma \ref{lem:irred-t-module-lp}.
\hfill $\Box$\\

\noindent
For notational convenience let $\theta_{-1}$ and $\theta_{D+1}$ (resp. $\theta^*_{-1}$ and $\theta^*_{D+1}$) denote the scalars in $\fld$ which satisfy Lemma \ref{lem:drgscalars}{\rm (ii)} (resp. Lemma \ref{lem:drgscalars}{\rm (iii)}) for $i = 0$ and $i = D$.

\begin{lemma}
\label{cor:drgaw2}
Let $W$ denote a thin irreducible $T$-module.
Let $r, t, d$ denote the endpoint, dual endpoint, and diameter of $W$, respectively.
Let $\omega, \eta, \eta^*$ denote scalars in $\C$.
Then $\omega,\eta,\eta^*$ satisfy {\rm (\ref{eq:drg-aw-rel1}), (\ref{eq:drg-aw-rel2})} if and only if the following {\rm (i)--(iv)} hold.
\begin{enumerate}
\item[\rm (i)] $\omega = a_i(W)(\theta_{r+i}^* - \theta_{r+i+1}^*) + a_{i-1}(W) (\theta_{r+i-1}^* - \theta_{r+i-2}^*) - \gamma (\theta_{r+i}^* + \theta_{r+i-1}^*) \qquad (1 \leq i \leq d).$
\item[\rm (ii)] $\omega = a_i^*(W)(\theta_{t+i} - \theta_{t+i+1}) + a_{i-1}^*(W) (\theta_{t+i-1} - \theta_{t+i-2}) - \gamma^* (\theta_{t+i} + \theta_{t+i-1}) \qquad (1 \leq i \leq d).$
\item[\rm (iii)] $\eta = a_i^*(W)(\theta_{t+i} - \theta_{t+i-1})(\theta_{t+i} - \theta_{t+i+1}) - \gamma^* \theta_{t+i}^2 - \omega \theta_{t+i} \qquad (0 \leq i \leq d).$
\item[\rm (iv)] $\eta^* = a_i(W) (\theta_{r+i}^* - \theta_{r+i-1}^*)(\theta_{r+i}^* - \theta_{r+i+1}^*) - \gamma \theta_{r+i}^{*2} - \omega \theta_{r+i}^* \qquad (0 \leq i \leq d).$
\end{enumerate}
\end{lemma}

\noindent {\it Proof:} \rm
Apply Lemma \ref{lem:lpawconsts} to the Leonard system from Lemma \ref{lem:irred-t-module-lp}.
\hfill $\Box$\\

\noindent
We will be discussing the center of $T$, denoted $Z(T)$.

\begin{lemma}
\label{lem:omega-z-zstar}
Suppose $\Gamma$ is thin.
Then there exist elements $\Omega, G, G^*$ of $Z(T)$ with the following property.
For every irreducible $T$-module $W$, the elements $\Omega, G, G^*$ act on $W$ as $\omega I, \eta I, \eta^* I$ where $\omega = \omega(W), \eta = \eta(W), \eta^* = \eta^*(W)$ are from Lemma \ref{lem:w-aw}.
\end{lemma}

\noindent {\it Proof:} \rm
Let $\{W_i\}_{i=1}^n$ denote a full set of mutually nonisomorphic irreducible $T$-modules.
For $1 \leq i \leq n$ let $V_i$ denote the subspace of $V$ spanned by the irreducible $T$-modules that are isomorphic to $W_i$.
Observe that each $V_i$ is a $T$-submodule of $V$.
By construction
\begin{eqnarray*}
V = V_1 + V_2 + \cdots + V_n\qquad \mbox{(direct sum)}.
\end{eqnarray*}
There exists matrices $\Omega, G, G^*$ in ${\rm Mat_X(\C)}$ that act on $V_i$ as $\omega(W_i)I, \eta(W_i)I, \eta^*(W_i)I$, respectively for $1 \leq i \leq n$.
Now consider our irreducible $T$-module $W$.
By construction there exists an integer $j$ such that $W$ is contained in $V_j$.
Therefore $\Omega, G, G^*$ act on $W$ as $\omega(W_j)I, \eta(W_j)I, \eta^*(W_j)I$, respectively.
In particular $\Omega, G, G^*$ act on $W$ as scalar multiples of $I$ and leave $W$ invariant.
By this and a similar argument to the proof of \cite[Lemma~12.1]{drgqtet}, each of $\Omega, G, G^*$ is in $Z(T)$.
Since the $T$-module $W$ is irreducible and is contained in $V_j$, the $T$-modules $W, W_j$ are isomorphic.
Therefore $\omega(W) = \omega(W_j), \eta(W) = \eta(W_j), \eta^*(W) = \eta^*(W_j)$.
By these comments the $\Omega, G, G^*$ act on $W$ as $\omega(W)I, \eta(W)I, \eta^*(W)I$, respectively.
\hfill $\Box$\\

\begin{lemma}
\label{lem:drgaw}
Suppose $\Gamma$ is thin.
Let $\beta, \gamma, \gamma^*, \varrho, \varrho^*$ denote the scalars from Lemma \ref{lem:drg-td-relns}.
Let $\Omega, G, G^*$ denote the elements of $T$ from Lemma \ref{lem:omega-z-zstar}.
Then both
\begin{align}
A^2A^* - \beta AA^*A + A^*A^2 - \gamma (AA^*+A^*A) - \varrho A^* &= \gamma^*A^2 + \Omega A + G,\label{drg-aw-rel-m1}\\
A^{*2}A - \beta A^*AA^* + AA^{*2} - \gamma^* (AA^*+A^*A) - \varrho^* A &= \gamma A^{*2} + \Omega A^* + G^*\label{drg-aw-rel-m2}.
\end{align}
\end{lemma}

\noindent {\it Proof:} \rm
Write $V$ as a direct sum of irreducible $T$-modules.
By Lemma \ref{lem:w-aw} and Lemma \ref{lem:omega-z-zstar}, lines (\ref{drg-aw-rel-m1}), (\ref{drg-aw-rel-m2}) hold on each irreducible $T$-module in the sum.
Therefore lines (\ref{drg-aw-rel-m1}), (\ref{drg-aw-rel-m2}) hold on $V$.
\hfill $\Box$\\

\noindent
We mention a result about $Z(T)$ for later use.

\begin{lemma}
\label{lem:centralcombi1}
Let $C$ denote an element in $Z(T)$.  Then for vertices $y,z \in X$
\begin{eqnarray*}
	C_{yz} \neq 0 \qquad \rightarrow \qquad \partial(x,y) = \partial(x,z).
\end{eqnarray*}
\end{lemma}

\noindent {\it Proof:} \rm
Assume $C_{yz} \neq 0$.
Let $i = \partial(x,y)$ and $j = \partial(x,z)$.
We show $i = j$.
Suppose $i \neq j$.
Since $E^*_iC = CE^*_i$ and $E^*_iE^*_j = 0$,
we have $E^*_iCE^*_j = E^*_iE^*_jC = 0$.
However
$(E^*_iCE^*_j)_{yz} = (E^*_i)_{yy} C_{yz} (E^*_j)_{zz} = C_{yz} \neq 0,$
a contradiction.
Therefore $i = j$.
\hfill $\Box$\\

\section{Near polygons}

We continue to discuss the distance-regular graph $\Gamma$ from Section \ref{sec:drg}.
In this section we consider the case in which $\Gamma$ is a near polygon. 
A clique in $\Gamma$ is called {\it maximal} whenever it is not properly contained in another clique.

\begin{definition} \rm \cite[p.~198]{bcn}
The graph $\Gamma$ is called a {\it near polygon} whenever the following two axioms hold.
\begin{description}
\item[(NP1)] There are no induced subgraphs of shape $K_{1,2,1}$.
\item[(NP2)] For a vertex $x$ in $X$ and a maximal clique $M$ of $\Gamma$ with $\partial(x,M) < D$, there exists a unique vertex in $M$ nearest to $x$.
\end{description}
\end{definition}

\begin{definition} \rm \cite[p.~198]{bcn}
Suppose $\Gamma$ is a near polygon.  Then $\Gamma$ is called a {\it near $n$-gon}, where $n = 2D+1$ if there is a vertex at distance $D$ from some maximal clique, and $n = 2D$ otherwise.
\end{definition}

\begin{lemma} {\rm \cite[Theorem~6.4.1]{bcn}}
The graph $\Gamma$ is a near polygon if and only if the axiom {\rm (NP1)} holds and $a_i = a_1 c_i$ for $1 \leq i \leq D-1$.
In this case $\Gamma$ is a near $(2D+1)$-gon if $a_D \neq a_1 c_D$ and a near $2D$-gon if $a_D = a_1 c_D$.
\end{lemma}

\begin{lemma}{\rm \cite[p.~200]{bcn}}
\label{lem:edge-clique}
Assume $\Gamma$ is a near polygon.  Then each edge in $\Gamma$ is contained in a unique maximal clique, and this clique has cardinality $a_1+2$.
\end{lemma}

\begin{lemma}
\label{lem:np_yz1}
Assume $\Gamma$ is a near polygon.
Fix $x \in X$.
Fix adjacent $y,z \in X$ such that $\partial(x,y) = \partial(x,z)-1$.
Define $i = \partial(x,y)$.
Then
\begin{eqnarray}
\Gamma_{i}(x) \cap \Gamma(y) \cap \Gamma(z) = \varnothing, \qquad |\Gamma_{i+1}(x) \cap \Gamma(y) \cap \Gamma(z)| = a_1.\label{eq:xyz1}
\end{eqnarray}
\end{lemma}

\noindent {\it Proof:} \rm
Since $y,z$ are adjacent, we have $|\Gamma(y) \cap \Gamma(z)| = a_1$.
Since $\partial(x,y) = i$ and $\partial(x,z) = i+1$, we have $\Gamma(y) \cap \Gamma(z) \subseteq \Gamma_{i}(x) \cup \Gamma_{i+1}(x)$.
By these comments
\begin{eqnarray*}
|\Gamma_{i}(x) \cap \Gamma(y) \cap \Gamma(z)| + |\Gamma_{i+1}(x) \cap \Gamma(y) \cap \Gamma(z)| = |\Gamma(y) \cap \Gamma(z)| = a_1.
\end{eqnarray*}
By this it suffices to show that $\Gamma_{i}(x) \cap \Gamma(y) \cap \Gamma(z) = \varnothing$.
Suppose there exists a vertex $w \in \Gamma_{i}(x) \cap \Gamma(y) \cap \Gamma(z)$.
Let $C$ denote a maximal clique of $\Gamma$ that contains $y,z$.
Observe that $w \in C$ by (NP1).
By construction $\partial(x,C) = i$.
By (NP2) there exists a unique vertex in $C$ at distance $i$ from $x$.
However $y,w$ are distinct vertices in $C$ at distance $i$ from $x$, a contradiction.
The result follows.
\hfill $\Box$\\

\noindent We recall some definitions for future use.

\begin{definition} \label{def:weak-geod-closed} \rm \cite[Definition~3.1]{weakgeod}
Assume $\Gamma$ is a near polygon.  A subgraph $H$ of $\Gamma$ is called {\em weak-geodetically closed} whenever for all $x,y \in H$ and for all $z \in X$
\begin{eqnarray*}
\partial(x,z) + \partial(z,y) \leq \partial(x,y) + 1 \qquad \rightarrow \qquad z \in H.
\end{eqnarray*}
\end{definition}

\begin{lemma}
Assume $\Gamma$ is a near polygon.
Let $H$ denote a weak-geodetically closed subgraph of $\Gamma$.
For all $x,y \in H$ we have
\begin{eqnarray*}
\partial(x,y) \;(\mbox{in } H) = \partial(x,y) \;(\mbox{in } \Gamma).
\end{eqnarray*}
\end{lemma}

\noindent {\it Proof:} \rm
Clear.
\hfill $\Box$\\

\begin{definition} \rm \cite[p.~145]{quads}
\label{def:quads2}
Assume $\Gamma$ is a near polygon.
A subgraph $Q$ of $\Gamma$ is called a {\it quad} whenever $Q$ has diameter $2$ and $Q$ is weak geodetically-closed.
\end{definition}

\section{Dual polar graphs}\label{sec:dpg}

In this section we discuss a type of near polygon called a dual polar graph.
Let $b$ denote a prime power.
Let $\F_b$ denote the finite field of order $b$.
Let $U$ denote a finite-dimensional vector space over $\F_b$ endowed with one of the following nondegenerate forms.

\medskip
\centerline{
\begin{tabular}[b]{l|l|l|l}
name & ${\rm dim} \,U$ & form & e\\ \hline
$C_D(b)$ & $2D$ & symplectic & 1\\
$B_D(b)$ & $2D+1$ & quadratic & 1\\
$D_D(b)$ & $2D$ & quadratic & 0\\
& & (Witt index $D$) &\\
$^2D_{D+1}(b)$ & $2D+2$ & quadratic & 2\\
& & (Witt index $D$) &\\
$^2A_{2D}(q)$ & $2D+1$ & Hermitean $(b=q^2)$ & $\frac{3}{2}$\\
$^2A_{2D-1}(q)$ & $2D$ & Hermitean $(b=q^2)$ & $\frac{1}{2}$
\end{tabular}
}

\medskip
\noindent A subspace $W$ of $U$ is called {\it isotropic} whenever the form vanishes completely on $W$.
By \cite[Theorem~6.3.1]{cameron}, each maximal isotropic subspace of $U$ has dimension $D$.
Define a graph as follows.
The vertex set consists of the maximal isotropic subspaces of $U$.
Vertices $y,z$ are adjacent whenever ${\rm dim}(y \cap z) = D-1$.
By \cite[p.~274]{bcn} this graph is distance-regular and has diameter $D$.
By \cite[p.~276]{bcn} for vertices $y,z$ we have $\partial(y,z) = D - {\rm dim}(y \cap z)$.
We call this graph the {\em dual polar graph} associated with $U$.

\medskip
\noindent For the rest of the paper assume $\Gamma$ is a dual polar graph.
In the next few lemmas we recall some basic data about $\Gamma$.

\begin{lemma} {\rm \cite[Theorem~9.4.3]{bcn}}
\label{lem:dpgintarray}
The intersection numbers $c_i, b_i$ of $\Gamma$ are given by
\begin{eqnarray*}
c_i = \frac{b^i-1}{b-1}, \qquad
b_i = \frac{b^{i+e}(b^{D-i}-1)}{b-1} \qquad
(0 \leq i \leq D).
\end{eqnarray*}
In particular the valency $k = b_0$ is
\begin{eqnarray*}
k = \frac{b^e(b^D-1)}{b-1}.
\end{eqnarray*}
\end{lemma}

\begin{corollary}
\label{cor:dpg_ai}
The intersection numbers $a_i$ of $\Gamma$ are given by
\begin{eqnarray*}
\label{eq:ai}
a_i = \frac{(b^e-1)(b^i-1)}{b-1} \qquad (0 \leq i \leq D).
\end{eqnarray*}
\end{corollary}

\noindent {\it Proof:} \rm
Use Lemma \ref{lem:dpgintarray} and $a_i = k - c_i - b_i$.
\hfill $\Box$\\

\begin{lemma} {\rm \cite[Proposition~10.4.1]{cameron}}
\label{lem:drgisnp}
The graph $\Gamma$ is a near $2D$-gon.
\end{lemma}

\begin{lemma} {\rm \cite[Theorem~9.4.3]{bcn}}
\label{lem:dpgeval}
Let $\theta_0 > \theta_1 > \cdots > \theta_D$ denote the eigenvalues of $\Gamma$.  Then
\begin{align*}
\theta_i &= \frac{1-b^e+b^{D+e-i}-b^i}{b-1},\\
m_i &= \frac{b^i(b^D-1)(b^{D-1}-1)\cdots(b^{D-i+1}-1)}{(b^i-1)\cdots(b^2-1)(b-1)} \; \frac{1+b^{D+e-2i}}{1+b^{D+e-i}}\prod^{i}_{j=1}\frac{1+b^{D+e-j}}{1+b^{j-e}}
\end{align*}
for $0 \le i \le D$.
\end{lemma}

\begin{lemma}
\label{lem:dpgdualeval}
The graph $\Gamma$ is $Q$-polynomial with respect to the ordering $\theta_0 > \theta_1 > \cdots > \theta_D$ of the eigenvalues.
The corresponding dual eigenvalue sequence is
\begin{eqnarray*}
\theta_i^* = \zeta + \xi b^{-i} \qquad (0 \leq i \leq D),
\end{eqnarray*}
where
\begin{align*}
\zeta &= -\frac{b(b^{D+e-2}+1)}{b-1},\\
\xi &= \frac{b^2(b^{D+e-2}+1)(b^{D+e-1}+1)}{(b-1)(b^e+b)}.
\end{align*}
\end{lemma}

\noindent {\it Proof:}  \rm
By \cite[Table~6.1,~Corollary~8.4.2]{bcn} the graph $\Gamma$ is $Q$-polynomial with respect to $\{\theta_i\}_{i=0}^D$.
By \cite[Theorem~8.4.1]{bcn} the corresponding dual eigenvalues satisfy
\begin{eqnarray}
\label{eq:dualevalratio}
\frac{\theta_i^*}{\theta_0^*} = 1+\frac{\theta_1-k}{k} \; \frac{b-b^{1-i}}{b-1} \qquad (0 \leq i \leq D).
\end{eqnarray}
In (\ref{eq:dualevalratio}) evaluate $k$ using Lemma \ref{lem:dpgintarray} and evaluate $\theta_0^*$ using $\theta_0^* = m_1$ and Lemma \ref{lem:dpgeval}.  The result follows.
\hfill $\Box$\\

\medskip
\noindent
From now on it is understood that the eigenvalues of $\Gamma$ are ordered such that $\theta_0 > \theta_1 > \cdots > \theta_D$.

\begin{lemma}
{\rm \cite[Theorem~17.19]{cerzo}}
The dual intersection numbers of $\Gamma$ are given by
\begin{align*}
c_i^* &= \frac{ \xi b^{-i} (b^i-1)(b^{e-i}+1) }{ (b^{D+e-2i}+1)(b^{D+e-2i+1}+1) },\\
b_i^* &= \frac{ \xi (1-b^{i-D})(b^{-D-e+i}+1) }{ (b^{-D-e+2i}+1)(b^{-D-e+2i+1}+1) },\\
a_i^* &= \theta_0^* - c_i^* - b_i^*
\end{align*}
for $0 \leq i \leq D$.
\end{lemma}

\begin{lemma}
	\label{lem:dpgconsts}
	The scalars $\beta,\gamma,\gamma^*,\varrho,\varrho^*$ from Lemma \ref{lem:drg-td-relns} satisfy the following {\rm (i)--(v)}.
	\begin{enumerate}
		\item[\rm (i)] $\beta = b + b^{-1}$.
		\item[\rm (ii)] $\displaystyle{\gamma = \frac{(b^e-1)(b-1)}{b}}$.
		\item[\rm (iii)] $\gamma^* = (b-1)(b^{D+e-2}+1)$.
		\item[\rm (iv)] $\displaystyle{\varrho = \frac{(b^e-1)^2}{b} + b^{D+e-2}(b+1)^2}$.
		\item[\rm (v)] $\varrho^* = b(b^{D+e-2}+1)^2$.
	\end{enumerate}
\end{lemma}

\noindent {\it Proof:} \rm
Routine verification using Lemma \ref{lem:drgscalars}, Lemma \ref{lem:dpgeval} and Lemma \ref{lem:dpgdualeval}.
\hfill $\Box$\\

\section{The structure of an irreducible $T$-module for a dual polar graph}

We continue to discuss the dual polar graph $\Gamma$ from Section \ref{sec:dpg}.
For the rest of the paper fix a vertex $x \in X$ and write $T = T(x)$ for the subconstituent algebra.
In this section we recall some basic data about irreducible $T$-modules.

\begin{lemma}
\label{lem:w-thin}
{\rm \cite[Example~6.1]{terwSub3}}
The graph $\Gamma$ is thin.
\end{lemma}

\begin{lemma} {\rm \cite[p.~200]{terwSub3}}
\label{lem:dpgintw}
Let $W$ denote an irreducible $T$-module.
Let $r, t, d$ denote the endpoint, dual endpoint, and diameter of $W$, respectively.
Then for $0 \leq i \leq d$
\begin{align*}
c_i(W) & = \frac{b^t(b^{i}-1)}{b-1},\\
b_i(W) & = \frac{b^{D+e-d-t+i}(b^{d-i}-1)}{b-1},\\
a_i(W) & = \frac{b^{D+e-d-t+i}-b^e-b^{t+i}+1}{b-1}.
\end{align*}
\end{lemma}

\begin{lemma} {\rm \cite[Theorem~17.17]{cerzo}}
\label{lem:dpgdintw}
Let $W$ denote an irreducible $T$-module.
Let $r, t, d$ denote the endpoint, dual endpoint, and diameter of $W$, respectively.
Then for $0 \leq i \leq d$
\begin{align*}
c_i^*(W) & = \frac{\xi b^{-r-i}(b^i-1)(b^{D+e-2t-d-i}+1)}{(b^{D+e-2t-2i}+1)(b^{D+e-2t-2i+1}+1)},\\
b_i^*(W) & = \frac{\xi b^{-r}(1-b^{i-d})(b^{-D-e+2t+i}+1)}{(b^{-D-e+2t+2i}+1)(b^{-D-e+2t+2i+1}+1)},\\
a_i^*(W) & = \theta_r^* - b_i^*(W) - c_i^*(W),
\end{align*}
where $\xi$ is the scalar from Lemma \ref{lem:dpgdualeval}.
\end{lemma}

\begin{lemma}
\label{lem:omega_eta_etastar}
Let $W$ denote an irreducible $T$-module.
Let $r, t, d$ denote the endpoint, dual endpoint and diameter of $W$, respectively.
Let $\zeta, \xi$ denote the scalars from Lemma \ref{lem:dpgdualeval}.
Let $\gamma,\varrho$ denote the scalars from Lemma \ref{lem:dpgconsts}.
Let $\omega=\omega(W),\eta=\eta(W),\eta^*=\eta^*(W)$ denote the scalars in $\C$ that satisfy the following {\rm (i)--(iii)}.
\begin{enumerate}
	\item[\rm (i)] $\displaystyle{\omega = \xi (b^{-1}-1) (b^{e+D-d-t-r}-b^{t-r}) - 2\gamma\zeta}$.
	\item[\rm (ii)] $\displaystyle{\eta = \xi b^{-1}(1-b^e)(b^{e+D-d-r-t}-b^{t-r})
	 - \xi b^{e+D-d-r-2}(b+1)(b^{d+1}+1)
	 -\varrho \zeta}$.
	\item[\rm (iii)] $\eta^* = -\gamma \zeta^2 - \zeta \omega$.
\end{enumerate}
Then $\omega,\eta,\eta^*$ satisfy {\rm (\ref{eq:drg-aw-rel1}), (\ref{eq:drg-aw-rel2})}.
\end{lemma}

\noindent {\it Proof:} \rm
We verify Lemma \ref{cor:drgaw2}(i)--(iv) using Lemma \ref{lem:dpgeval}, Lemma \ref{lem:dpgdualeval}, Lemma \ref{lem:dpgconsts}, Lemma \ref{lem:dpgintw} and Lemma \ref{lem:dpgdintw}.
The result follows from Lemma \ref{cor:drgaw2}.
\hfill $\Box$\\

\noindent
We give a comment for future use.

\begin{lemma}
\label{lem:isom-w}
Let $W$ denote an irreducible $T$-module.
Let $r, t, d$ denote the endpoint, dual endpoint and diameter of $W$, respectively.
Then the isomorphism class of $W$ is determined by $r, t, d$.
\end{lemma}

\noindent {\it Proof:}
By Lemma \ref{lem:dpgintw} the action of $A$ on $W$ is determined by $t, d$.
By Lemma \ref{lem:dpgdintw} the action of $A^*$ on $W$ is determined by $r, t, d$.
By these comments and since $A, A^*$ generate $T$, the action of $T$ on $W$ is determined by $r, t, d$.
\hfill $\Box$\\

\section{Some combinatorial aspects of a dual polar graph}

We continue to discuss the dual polar graph $\Gamma$ from Section \ref{sec:dpg}.
In this section we discuss some combinatorial aspects of $\Gamma$.
Recall the quads of $\Gamma$ from Definition \ref{def:quads2}.

\begin{lemma}
{\rm \cite[p.~132]{cameron}}
\label{lem:quad_2v}
For all $y, z \in X$ such that $\partial(y,z) = 2$, there exists a unique quad containing $y, z$.
\end{lemma}

\begin{lemma} {\rm \cite[p.~132]{cameron}}
\label{lem:quad}
Let $Q$ denote a quad in $\Gamma$.
For all $u \in X$ there exists a unique $v \in Q$ nearest to $u$.
Moreover for any $w \in Q$, we have $\partial(u,w) = \partial(u,v) + \partial(v,w)$.
\end{lemma}

\begin{lemma}
\label{lem:dpg_yz2}
Fix $y, z \in X$ such that $\partial(x,y) = \partial(x,z) - 1$ and $\partial(y,z) = 2$.
Define $i = \partial(x,y)$.
Then
\begin{eqnarray}
|\Gamma_i(x) \cap \Gamma(y) \cap \Gamma(z)| = 1, \qquad |\Gamma_{i+1}(x) \cap \Gamma(y) \cap \Gamma(z)| = b.\label{eq:xyz2}
\end{eqnarray}
\end{lemma}

\noindent {\it Proof:} \rm
By Lemma \ref{lem:quad_2v} there exists a unique quad $Q$ containing $y, z$.
Since $\partial(y,z) = 2$, we have $|\Gamma(y) \cap \Gamma(z)| = c_2$.
By Lemma \ref{lem:dpgintarray} we have $c_2 = b+1$.
Since $\partial(x,y) = i$ and $\partial(x,z) = i+1$, the intersection $\Gamma(y) \cap \Gamma(z) \subseteq \Gamma_{i}(x) \cup \Gamma_{i+1}(x)$.
By these comments
\begin{eqnarray*}
|\Gamma_i(x) \cap \Gamma(y) \cap \Gamma(z)| + |\Gamma_{i+1}(x) \cap \Gamma(y) \cap \Gamma(z)| = |\Gamma(y) \cap \Gamma(z)| = b+1.
\end{eqnarray*}
Therefore it suffices to show that $|\Gamma_i(x) \cap \Gamma(y) \cap \Gamma(z)| = 1$.
First we show $\Gamma_i(x) \cap \Gamma(y) \cap \Gamma(z) \neq \varnothing$.
Suppose $\Gamma_i(x) \cap \Gamma(y) \cap \Gamma(z) = \varnothing$.
Then $\Gamma_{i+1}(x) \cap \Gamma(y) \cap \Gamma(z) \neq \varnothing$.
Let $u \in \Gamma_{i+1}(x) \cap \Gamma(y) \cap \Gamma(z)$.
Observe that $u \in Q$.
Moreover the maximal clique containing the edge $uz$ lies in $Q$ and, by (NP2), contains a unique vertex $v$ at distance $i$ from $x$.
Observe that $y, v$ are distinct vertices in $Q$ each at distance $i$ from $x$.
By these comments, Lemma \ref{lem:quad} and since $Q$ has diameter $2$,
there exists a unique vertex $p$ in $Q$ at distance $i-1$ from $x$.
By Lemma \ref{lem:quad} the vertices $p, y$ are adjacent.
Let $P$ denote a maximal clique of $\Gamma$ that contains $p, y$.
Observe that $P$ is contained in $Q$.
Since $\partial(y,z) = 2$ and $Q$ has diameter 2, we have $\partial(z,P) = 1$.
Thus $z$ is adjacent to a unique vertex in $P$ and by construction it is different from $p$ and $y$.
This vertex is in $\Gamma_{i}(x) \cap \Gamma(y) \cap \Gamma(z)$, a contradiction.
Hence $\Gamma_i(x) \cap \Gamma(y) \cap \Gamma(z) \neq \varnothing$.
Now suppose $|\Gamma_i(x) \cap \Gamma(y) \cap \Gamma(z)| > 1$.
Let $u,v$ denote distinct vertices in $\Gamma_i(x) \cap \Gamma(y) \cap \Gamma(z)$.
By (NP1) the vertices $u,v$ are not adjacent.
Observe that $Q$ contains $u, v$.
By these comments, Lemma \ref{lem:quad} and since $Q$ has diameter $2$, there exists a unique vertex $p$ in $Q$ at distance $i-1$ from $x$.
By Lemma \ref{lem:quad} the vertex $p$ is adjacent to $u, v, y$.
The vertices $p, u, v, y$ induce a subgraph of shape $K_{1,2,1}$ contradicting (NP1).
Thus $|\Gamma_i(x) \cap \Gamma(y) \cap \Gamma(z)| = 1$.
The result follows.
\hfill $\Box$\\

\section{A basis for $Z(T)$}\label{sec:centeroft}

We continue to discuss the dual polar graph $\Gamma$ from Section \ref{sec:dpg}.
Recall the subconstituent algebra $T$.
In this section we display a basis for $Z(T)$.
For $0 \leq r, t, d \leq D$, the triple  $(r,t,d)$ is called {\em feasible} whenever there exists an irreducible $T$-module with endpoint $r$, dual endpoint $t$, and diameter $d$.
Let ${\rm Feas}$ denote the set of all feasible triples.
For $\lambda = (r,t,d) \in {\rm Feas}$,
let $V_\lambda$ denote the subspace of $V$ spanned by the irreducible $T$-modules with endpoint $r$, dual endpoint $t$ and diameter $d$.
Observe that $V_\lambda$ is a $T$-module.
We call $V_\lambda$ the {\em homogeneous component of $V$ associated with $\lambda$.}
Observe
\begin{eqnarray*}
V = \sum_{\lambda \in {\rm Feas}} V_{\lambda} \qquad ({\rm direct \ sum}).
\end{eqnarray*}

\noindent
For all $\lambda \in {\rm Feas}$, define the linear transformation $E_{\lambda}:V \rightarrow V$ by
\begin{eqnarray*}
&&(E_{\lambda} - I)V_{\lambda} = 0,\\
&&E_{\lambda} V_{\lambda'} = 0 \qquad \mbox{if } \lambda' \neq \lambda \qquad (\lambda' \in {\rm Feas}).
\end{eqnarray*}

\noindent
Observe that
\begin{eqnarray*}
&&I = \sum_{\lambda \in {\rm Feas}} E_{\lambda},\\
&&E_{\lambda} E_{\lambda'} = \delta_{\lambda\lambda'}E_{\lambda} \qquad (\lambda,\lambda' \in {\rm Feas}).
\end{eqnarray*}

\begin{lemma}
\label{lem:centerofT}
{\rm \cite[Theorem~25.15, Theorem~26.4]{CR}}
The elements $\{E_{\lambda}|\lambda \in {\rm Feas}\}$ form a basis for the vector space $Z(T)$.
\end{lemma}

\section{The central elements $\Omega, G^*$}

We continue to discuss the dual polar graph $\Gamma$ from Section \ref{sec:dpg}.
Recall the subconstituent algebra $T$ and its central elements $\Omega, G^*$ from Lemma \ref{lem:omega-z-zstar}.
In this section we show that $G^*$ is a linear combination of $\Omega$ and $I$.
Then we display $\Omega$ in a certain attractive form.
Using this form we obtain a characterization of $\Omega$.

\begin{lemma}
\label{lem:zstar}
The elements $\Omega, G^*$ are related by
\begin{eqnarray}
G^* = -\gamma \zeta^2 I - \zeta \Omega\label{eq:zstar}
\end{eqnarray}
where $\zeta$ is from Lemma \ref{lem:dpgdualeval} and $\gamma$ is from Lemma \ref{lem:dpgconsts}.
\end{lemma}

\noindent {\it Proof:} \rm
By construction $V$ is a direct sum of irreducible $T$-modules.
Let $W$ denote an irreducible $T$-module in the sum.
It suffices to show that the two sides of (\ref{eq:zstar}) agree on $W$.
By Lemma \ref{lem:omega-z-zstar} the elements $\Omega, G^*$ act on $W$ as $\omega(W)I, \eta^*(W)I$, respectively.
By Lemma \ref{lem:omega_eta_etastar}(iii) we have $\eta^*(W) = -\gamma \zeta^2 - \zeta \omega(W)$.
Therefore the two sides of (\ref{eq:zstar}) agree on $W$.
The result follows.
\hfill $\Box$\\

\begin{lemma}
\label{lem:Omega}
We have
\begin{eqnarray}
	\Omega = \sum^{D}_{i=1} \alpha_i E^*_i A E^*_i + \sum^{D}_{i=0} \beta_i E^*_i
	\label{eq:Omega}
\end{eqnarray}
where
\begin{align*}
\alpha_i & = \frac{(b^{D+e-1}+1)(b^{D+e-2}+1)(1-b)b^{-i}}{b^{e-1}+1} \qquad (1 \leq i \leq D),\\
\beta_i & = \frac{(b^{D+e-2}+1)(b^e-1)(2b^{e-1} + 2 - (b^{D+e-1}+1)b^{-i})}{b^{e-1}+1} \qquad (0 \leq i \leq D).
\end{align*}
\end{lemma}

\noindent {\it Proof:} \rm
Let $\tilde{\Omega}$ denote the expression on the right of (\ref{eq:Omega}).
By construction $V$ is a direct sum of irreducible $T$-modules.
Let $W$ denote an irreducible $T$-module in the sum.
To show $\tilde{\Omega} = \Omega$ it suffices to show that $\tilde{\Omega}, \Omega$ agree on $W$.
Let $r, d$ denote the endpoint and diameter of $W$.
By construction for $r \leq i \leq r+d$ the element $\tilde{\Omega}$ acts on $E^*_iW$ as $(\alpha_i a_{i-r}(W) + \beta_i)I$.
By evaluating $a_{i-r}(W)$ using Lemma \ref{lem:dpgintw} we find $\alpha_i a_{i-r}(W) + \beta_i = \omega(W)$ where $\omega(W)$ is the scalar from Lemma \ref{lem:omega_eta_etastar}(i).
Therefore $\tilde{\Omega}$ acts on $E^*_iW$ as $\omega(W)I$.
By this and since $W$ is a direct sum of $\{E^*_iW\}_{i=r}^{r+d}$, the element $\tilde{\Omega}$ acts on $W$ as $\omega(W)I$.
By Lemma \ref{lem:omega-z-zstar}, the element $\Omega$ acts on $W$ as $\omega(W)I$.
Therefore $\Omega, \tilde{\Omega}$ agree on $W$.
The result follows.
\hfill $\Box$\\

\medskip
\noindent
Motivated by Lemma \ref{lem:Omega} we consider an element $C$ of $T$ of the form
\begin{eqnarray}
\label{eq:centralc}
C := \sum^{D}_{i=1} \alpha_i E^*_i A E^*_i + \sum^{D}_{i=0} \beta_i E^*_i
\end{eqnarray}
for some arbitrary scalars $\alpha_i \in \C$ $(1 \leq i \leq D)$ and $\beta_i \in \C$ $(0 \leq i \leq D)$.

\medskip
\noindent Observe that $C$ is invariant under transposition.  We find necessary and sufficient conditions on $\alpha_i, \beta_i$ for $C$ to be central in $T$.  By construction $C$ commutes with $A^*$.  Since $A,A^*$ generate $T$, the element $C$ is central in $T$ if and only if $C$ commutes with $A$.

\begin{lemma}
\label{lem:Cyz}
For vertices $y,z \in X$, the $(y,z)$-entry of $C$ is described as follows.  First assume $\partial(x,y) \neq \partial(x,z)$.
Then the $(y,z)$-entry of $C$ is zero.
Next assume $\partial(x,y) = \partial(x,z)$ and let $s$ denote this common distance.  Then the $(y,z)$-entry of $C$ is given by

\medskip
\centerline{
\begin{tabular}[t]{c|c}
        Case &$(y,z)$-entry of $C$ \\ \hline
        $y = z$ & $\beta_s$ \\
        $\partial(y,z) = 1$ & $\alpha_s$ \\
        $\partial(y,z) \geq 2$ & $0$
\end{tabular}
}
\end{lemma}

\noindent {\it Proof:} \rm
Routine consequence of (\ref{eq:centralc}).
\hfill $\Box$\\

\begin{lemma}
\label{lem:yz-entry}
For integers $i,j$ $(0 \leq i,j \leq D)$ and for vertices $y,z \in X$, the $(y,z)$-entries of $E_i^* A E_j^* A$ and $A E_i^* A E_j^*$ are described as follows.
\begin{enumerate}
\item[\rm(i)] $(E_i^* A E_j^* A)_{yz} = 
\begin{cases} 0, &\mbox{if } \partial(x,y) \neq i;\\
|\Gamma_j(x) \cap \Gamma(y) \cap \Gamma(z)|, &\mbox{if } \partial(x,y) = i. \end{cases}$
\item[\rm(ii)] $(A E_i^* A E_j^*)_{yz} = 
\begin{cases} 0, &\mbox{if } \partial(x,z) \neq j;\\
|\Gamma_i(x) \cap \Gamma(y) \cap \Gamma(z)|, &\mbox{if } \partial(x,z) = j. \end{cases}$
\end{enumerate}
\end{lemma}

\noindent {\it Proof:} \rm
(i) We have $(E^*_i A E^*_j A)_{yz} = \sum_{w \in X} (E^*_i A E^*_j)_{yw} A_{wz}$.
The result follows.\\
(ii) By (i) above and since $A E^*_i A E^*_j = (E_j^* A E_i^* A)^t $.
\hfill $\Box$\\

\begin{lemma}
\label{lem:zyzero}
Let $y,z$ denote vertices in $X$ such that $\partial(y,z) \geq 3$ or $|\partial(x,y)-\partial(x,z)| \geq 2$.
Then the $(y,z)$-entries of $AC$ and $CA$ are both zero.
\end{lemma}

\noindent {\it Proof:} \rm
First we show $(AC)_{yz} = 0$.
By (\ref{eq:centralc}) it suffices to show that each of  $(A E^*_i A E^*_i)_{yz}$ and $(A E^*_i)_{yz}$ is $0$ for $0 \leq i \leq D$.
Let $i$ be given.
By construction $A_{yz} = 0$ so $(AE^*_i)_{yz} = A_{yz}(E^*_i)_{zz} = 0$.
By Lemma \ref{lem:yz-entry}(ii), we have $(A E^*_i A E^*_i)_{yz} = 0$.
Therefore $(AC)_{yz} = 0$.
By swapping the roles of $y$ and $z$ we have $(AC)_{zy} = 0$.
By this and since $CA = (AC)^t$, we have $(CA)_{yz} = (AC)_{zy} = 0$.
\hfill $\Box$\\

\begin{lemma}
\label{lem:ac=ca1}
Let $y, z$ denote vertices in $X$ such that $\partial(x,y) = \partial(x,z) - 1$ and $\partial(y,z) = 2$.
Then the following are equivalent:
\begin{enumerate}
\item[\rm (i)] $(AC)_{yz} = (CA)_{yz}$.
\item[\rm (ii)] $b \alpha_{i+1} = \alpha_i$ where $i = \partial(x,y)$.
\end{enumerate}
\end{lemma}

\noindent {\it Proof:} \rm
Consider the $(y,z)$-entries of $AC$ and $CA$.
\begin{align*}
(AC)_{yz} & = \sum_{j=1}^D \alpha_j (A E_j^* A E_j^*)_{yz} + \sum_{j=0}^D \beta_j (A E_j^*)_{yz}\\
& = \alpha_{i+1} |\Gamma_{i+1}(x) \cap \Gamma(y) \cap \Gamma(z)|
\qquad \mbox{by Lemma \ref{lem:yz-entry}(ii)}\\
& = b \alpha_{i+1} \qquad \mbox{by the equation on the right of (\ref{eq:xyz2})},\\
(CA)_{yz} & = \sum_{j=1}^D \alpha_j (E_j^* A E_j^* A)_{yz} + \sum_{j=0}^D \beta_j (E_j^* A)_{yz}\\
& = \alpha_{i} |\Gamma_{i}(x) \cap \Gamma(y) \cap \Gamma(z)|
\qquad \mbox{by Lemma \ref{lem:yz-entry}(i)}\\
& = \alpha_{i} \qquad \mbox{by the equation on the left of (\ref{eq:xyz2})}.
\end{align*}
The result follows.
\hfill $\Box$\\

\begin{lemma}
\label{lem:ac=ca2}
Let $y, z$ denote vertices in $X$ such that $\partial(x,y) = \partial(x,z) - 1$ and $\partial(y,z) = 1$.
Then the following are equivalent:
\begin{enumerate}
\item[\rm (i)] $(AC)_{yz} = (CA)_{yz}$.
\item[\rm (ii)] $a_1 \alpha_{i+1} + \beta_{i+1} = \beta_i$ where $i = \partial(x,y)$.
\end{enumerate}
\end{lemma}

\noindent {\it Proof:} \rm
Consider the $(y,z)$-entries of $AC$ and $CA$.
\begin{align*}
(AC)_{yz} & = \sum_{j=1}^D \alpha_j (A E_j^* A E_j^*)_{yz} + \sum_{j=0}^D \beta_j (A E_j^*)_{yz}\\
& = \alpha_{i+1} |\Gamma_{i+1}(x) \cap \Gamma(y) \cap \Gamma(z)| + \beta_{i+1}
\qquad \mbox{by Lemma \ref{lem:yz-entry}(ii)}\\
& = a_1 \alpha_{i+1} + \beta_{i+1} \qquad \mbox{by the equation on the right of (\ref{eq:xyz1})},\\
(CA)_{yz} & = \sum_{j=1}^D \alpha_j (E_j^* A E_j^* A)_{yz} + \sum_{j=0}^D \beta_j (E_j^* A)_{yz}\\
& = \alpha_{i} |\Gamma_{i}(x) \cap \Gamma(y) \cap \Gamma(z)| + \beta_{i}
\qquad \mbox{by Lemma \ref{lem:yz-entry}(i)}\\
& = \beta_{i} \qquad \mbox{by the equation on the left of (\ref{eq:xyz1})}.
\end{align*}
The result follows.
\hfill $\Box$\\

\begin{lemma}
\label{lem:ac=ca3}
Let $y, z$ denote vertices in $X$ such that $\partial(x,y) = \partial(x,z)$ and $\partial(y,z) \leq 2$.
Then
\begin{eqnarray*}
(AC)_{yz} = (CA)_{yz}.
\end{eqnarray*}
\end{lemma}

\noindent {\it Proof:} \rm
Let $i$ denote the common value of $\partial(x,y)$ and $\partial(x,z)$.
We have
\begin{eqnarray*}
	&(A E_i^* A E_i^*)_{yz} = (E_i^* A E_i^* A E_i^*)_{yz} = (E_i^* A E_i^* A)_{yz},\\
	&(A E_i^*)_{yz} = (E_i^* A E_i^*)_{yz} = (E_i^* A)_{yz}.
\end{eqnarray*}
The result follows from this and (\ref{eq:centralc}).% and Lemma \ref{lem:yz-entry}.
\hfill $\Box$\\

\noindent
By combining Lemma \ref{lem:zyzero}, Lemma \ref{lem:ac=ca1}, Lemma \ref{lem:ac=ca2} and Lemma \ref{lem:ac=ca3} we have the following result.

\begin{theorem}
\label{lem:centralcond}
Consider the element $C$ from {\rm (\ref{eq:centralc})}.
Then $C$ is central in $T$ if and only if the following {\rm (i), (ii)} hold.
\begin{enumerate}
\item[\rm(i)] $\alpha_i = b^{1-i} \alpha_1 \qquad (1 \leq i \leq D)$.
\item[\rm(ii)] $\displaystyle{\beta_i = \beta_0 - a_1 b^{1-i} \frac{b^i-1}{b-1} \alpha_1 \qquad (0 \leq i \leq D)}$.
\end{enumerate}
\end{theorem}

\begin{corollary} Let $\mathcal{C}$ denote the subspace of $T$ consisting of central elements of the form {\rm (\ref{eq:centralc})}.  Then $\Omega$ and $I$ form a basis for $\mathcal{C}$.
\end{corollary}

\noindent {\it Proof:} \rm
In the formulas of $\alpha_i, \beta_i$ in Theorem \ref{lem:centralcond}, there are two free variables $\alpha_1, \beta_0$.
Therefore ${\rm dim}~\mathcal{C} \leq 2$.
Observe that $I \in \mathcal{C}$.
By Lemma \ref{lem:Omega} the element $\Omega \in \mathcal{C}$.
Since $I, \Omega$ are linearly independent, they form a basis for $\mathcal{C}$.
\hfill $\Box$\\

\section{The central element $G$}

We continue to discuss the dual polar graph $\Gamma$ from Section \ref{sec:dpg}.
Recall the subconstituent algebra $T$.
In this section we investigate the central element $G \in T$ from Lemma \ref{lem:omega-z-zstar}.

\begin{lemma}
	\label{lem:Z}
	We have
	\begin{eqnarray}\label{eq:Z}
	\begin{split}
		G &= \xi (1-b^2) \sum_{i=1}^D b^{-i} E_i^* A E_{i-1}^* A E_i^*
		+ \xi (1-b^{-2}) \sum_{i=0}^{D-1} b^{-i} E_i^* A E_{i+1}^* A E_i^*\\
		&\quad + \xi (b^{-1}-1) (b^e-1) \sum_{i=1}^D b^{-i} E_i^* A E_i^* 
		- \varrho A^*,
	\end{split}
	\end{eqnarray}
where $\xi$ is from Lemma \ref{lem:dpgdualeval} and $\varrho$ is from Lemma \ref{lem:dpgconsts}.
\end{lemma}

\noindent {\it Proof:} \rm
Recall that $I = \sum_{i=0}^D E_i^*$, $A^* = \sum_{i=0}^D \theta^*_i E^*_i$,  and $E_i^* E_j^* = \delta_{ij} E_i^*$ for $0 \leq i, j \leq D$.
By these facts and since $G$ is central in $T$,
\begin{eqnarray*}
	G = \left( \sum_{i=0}^D E_i^* \right) G \left( \sum_{j=0}^D E_j^* \right)
	= \sum_{i=0}^D \sum_{j=0}^D E_i^* G E_j^*
	= \sum_{i=0}^D E_i^* G E_i^*.
\end{eqnarray*}
For $0 \leq i \leq D$ we compute $E_i^*GE_i^*$.
By (\ref{drg-aw-rel-m1}) and since $A^*E_j^* = E_j^*A^* = \theta_j^* E_j^*$ $(0 \leq j \leq D)$,
\begin{eqnarray} \label{eq:eze}
\begin{split}
E_i^* G E_i^* 
&= (2\theta_i^* - \gamma^*) E_i^*A^2E_i^* - \beta E_i^*AA^*AE_i^*\\
&\quad - 2\gamma \theta_i^* E_i^*AE_i^* - \varrho \theta_i^* E_i^* - E_i^* \Omega A E_i^*.
\end{split}
\end{eqnarray}
We now evaluate the right-hand side of (\ref{eq:eze}).
First assume $1 \leq i \leq D-1$.
Using (\ref{eq:q-triple2}),
\begin{align*}
E_i^*A^2E_i^*
& = E_i^*A \left( \sum_{j=0}^D E_j^* \right) AE_i^*\\
& = E_i^*AE_{i-1}^*AE_i^* + E_i^*AE_{i}^*AE_i^* + E_i^*AE_{i+1}^*AE_i^*,\\
E_i^*AA^*AE_i^*
& = E_i^*A \left( \sum_{j=0}^D \theta_j^* E_j^* \right) AE_i^*\\
& = \theta_{i-1}^* E_i^*AE_{i-1}^*AE_i^* + \theta_{i}^* E_i^*AE_{i}^*AE_i^* + \theta_{i+1}^* E_i^*AE_{i+1}^*AE_i^*.
\end{align*}
Using (\ref{eq:Omega}),
\begin{eqnarray*}
E_i^*\Omega AE_i^*
=
\alpha_i E_i^*AE_i^*AE_i^* + \beta_i E_i^*AE_i^*
\end{eqnarray*}
where $\alpha_i, \beta_i$ are from Lemma \ref{lem:Omega}.
Evaluating the right-hand side of (\ref{eq:eze}) using the above comments,
\begin{align*}
E_i^*GE_i^*
&=
(2\theta_i^*-\gamma^* - \beta \theta_{i-1}^*)E_i^*AE_{i-1}^*AE_i^* + (2\theta_i^*-\gamma^* - \beta \theta_{i}^* - \alpha_i)E_i^*AE_{i}^*AE_i^*\\
&\quad + (2\theta_i^*-\gamma^* - \beta \theta_{i+1}^*)E_i^*AE_{i+1}^*AE_i^*
- (2 \gamma \theta_i^* + \beta_i) E_i^* A E_i^* - \varrho \theta_i^* E_i^*.
\end{align*}
By a similar argument,
\begin{align*}
E_0^*GE_0^* & =
 (2\theta_0^*-\gamma^* - \beta \theta_{1}^*)E_0^*AE_{1}^*AE_0^*
 - \varrho \theta_0^* E_0^*,\\
E_D^*GE_D^* & =
(2\theta_D^*-\gamma^* - \beta \theta_{D-1}^*)E_D^*AE_{D-1}^*AE_D^* + (2\theta_D^*-\gamma^* - \beta \theta_{D}^* - \alpha_D)E_D^*AE_{D}^*AE_D^*\\
&\quad - (2 \gamma \theta_D^* + \beta_D) E_D^* A E_D^* - \varrho \theta_D^* E_D^*.
\end{align*}
In the preceding equations we now evaluate the coefficients on the right-hand side.
The $\theta_i^*$ are from Lemma \ref{lem:dpgdualeval},
the $\beta, \gamma, \gamma^*$ are from Lemma \ref{lem:dpgconsts},
and $\alpha_i, \beta_i$ are from Lemma \ref{lem:Omega}.
By these lemmas,
\begin{align*}
	2\theta_i^*-\gamma^* - \beta \theta_{i-1}^* & = \xi (1-b^2)b^{-i} \qquad (1 \leq i \leq D),\\
	2\theta_i^*-\gamma^* - \beta \theta_{i}^* - \alpha_i & = 0 \qquad (1 \leq i \leq D),\\
	2\theta_i^*-\gamma^* - \beta \theta_{i+1}^* & = \xi (1-b^{-2})b^{-i} \qquad (0 \leq i \leq D-1),\\
	2 \gamma \theta_i^* + \beta_i & = \xi (1-b^{-1})(b^e-1)b^{-i} \qquad (0 \leq i \leq D).
\end{align*}
The result follows.
\hfill $\Box$\\

\begin{corollary}
\label{cor:Zyz}
For vertices $y,z \in X$, the $(y,z)$-entry of $G$ is described as follows.  First assume $\partial(x,y) \neq \partial(x,z)$.
Then the $(y,z)$-entry of $G$ is zero.
Next assume $\partial(x,y) = \partial(x,z)$ and let $s$ denote this common distance.
Then the $(y,z)$-entry of $G$ is given in the table below.

\medskip
\centerline{
\begin{tabular}[t]{c|c}
	Case &$(y,z)$-entry of $G$ \\ \hline
	$y = z$ & $\xi b^{-s-1} (1-b^2)(b c_s - b^{-1} b_s) - \varrho \theta_s^*$ \\
	$\partial(y,z) = 1$ & $\xi b^{-s-1} (1-b) (b^2+b+b^e-1)$ \\
	$\partial(y,z) = 2, \Gamma_{s+1}(x) \cap \Gamma(y) \cap \Gamma(z) = \varnothing$ & $\xi b^{-s} (1-b^2)|\Gamma_{s-1}(x) \cap \Gamma(y) \cap \Gamma(z)|$\\
	$\partial(y,z) = 2, \Gamma_{s+1}(x) \cap \Gamma(y) \cap \Gamma(z) \neq \varnothing$ & $-\xi b^{-s-1} (b+1)(b-1)^2 $\\
	$\partial(y,z) \geq 3$ & $0$
\end{tabular}
}
\medskip
In the above table $\xi$ is from Lemma \ref{lem:dpgdualeval} and $\varrho$ is from Lemma \ref{lem:dpgconsts}.
\end{corollary}

\noindent {\it Proof:} \rm
The first assertion follows from Lemma \ref{lem:centralcombi1}.
Now suppose $\partial(x,y) = \partial(x,z) = s$.
We verify the table.
By (\ref{eq:Z}), the $(y,z)$-entry of $G$ is given by
\begin{align*} %\label{eq:Z_yz}
G_{yz} & =
\xi (1-b^2) b^{-s} (E_s^* A E_{s-1}^* A E_s^*)_{yz}
+ \xi (1-b^{-2}) b^{-s} (E_s^* A E_{s+1}^* A E_s^*)_{yz}\\
&\quad + \xi (b^{-1}-1) (b^e-1) b^{-s} (E_s^* A E_s^*)_{yz}
- \varrho (A^*)_{yz}.
\end{align*}
In the above equation the terms on the right-hand side are given by
\begin{align*}
(E_s^* A E_{s-1}^* A E_s^*)_{yz}
& = |\Gamma_{s-1}(x) \cap \Gamma(y) \cap \Gamma(z)|
\qquad \mbox{by Lemma \ref{lem:yz-entry}},\\
(E_s^* A E_{s+1}^* A E_s^*)_{yz}
& = |\Gamma_{s+1}(x) \cap \Gamma(y) \cap \Gamma(z)|
\qquad \mbox{by Lemma \ref{lem:yz-entry}},\\
(E_s^* A E_s^*)_{yz} & = \begin{cases} 0 &\mbox{if } \partial(y,z) \neq 1;\\
1 &\mbox{if } \partial(y,z) = 1,\end{cases}\\
(A^*)_{yz} & = \sum_{i=0}^D \theta_i^* (E_i^*)_{yz}	= \theta_s^* (E_s^*)_{yz} = \delta_{yz} \theta_s^*.
\end{align*}
We now split our argument into cases.\\
Case $y=z$: We have
$|\Gamma_{s-1}(x) \cap \Gamma(y) \cap \Gamma(z)| = |\Gamma_{s-1}(x) \cap \Gamma(y)| = c_s$
and $|\Gamma_{s+1}(x) \cap \Gamma(y) \cap \Gamma(z)| = |\Gamma_{s+1}(x) \cap \Gamma(y)| = b_s$.
By these comments $G_{yz}$ is as shown in the table.\\
Case $\partial(y,z) = 1$: 
By (NP1) the vertices $y, z$ together with $\Gamma(y) \cap \Gamma(z)$ form a maximal clique of $\Gamma$.
By Lemma \ref{lem:drgisnp} this clique is at distance less than $D$ from $x$.
By (NP2) we have
$|\Gamma_{s-1}(x) \cap \Gamma(y) \cap \Gamma(z)| = 1$
and so $|\Gamma_{s+1}(x) \cap \Gamma(y) \cap \Gamma(z)| = 0$.
By this and the preliminary comments we find that $G_{yz}$ is as shown in the table.\\
Case $\partial(y,z) = 2$: Let $Q$ denote a quad containing $y,z$.
Observe that $\Gamma(y) \cap \Gamma(z)$ is contained in $Q$.
First assume $\Gamma_{s+1}(x) \cap \Gamma(y) \cap \Gamma(z) = \varnothing$.
By the preliminary comments we find that $G_{yz}$ is as shown in the table.
Next assume $\Gamma_{s+1}(x) \cap \Gamma(y) \cap \Gamma(z) \neq \varnothing$.
Then $\partial(x,Q) = s-1$ by Lemma \ref{lem:quad} and since $Q$ has diameter 2.
By Lemma \ref{lem:quad} we have $|\Gamma_{s-1}(x) \cap \Gamma(y) \cap \Gamma(z)| = 1$.
Let $u$ denote the unique vertex in $\Gamma_{s-1}(x) \cap \Gamma(y) \cap \Gamma(z)$.
We claim $\Gamma_{s}(x) \cap \Gamma(y) \cap \Gamma(z) = \varnothing$.
Suppose there exists a vertex $v$ in $\Gamma_{s}(x) \cap \Gamma(y) \cap \Gamma(z)$.
Then $v$ is adjacent to $u$ by Lemma \ref{lem:quad}.
Now $u,v,y,z$ induce a subgraph of shape $K_{1,2,1}$ which contradicts (NP1).
Hence the claim holds.
We have $|\Gamma(y) \cap \Gamma(z)| = c_2$ and $c_2 = b+1$ by Lemma \ref{lem:dpgintarray}.
By these comments and since $|\Gamma(y) \cap \Gamma(z)| = |\Gamma_{s-1}(x) \cap \Gamma(y) \cap \Gamma(z)| + |\Gamma_{s}(x) \cap \Gamma(y) \cap \Gamma(z)| + |\Gamma_{s+1}(x) \cap \Gamma(y) \cap \Gamma(z)|$,
we have $|\Gamma_{s+1}(x) \cap \Gamma(y) \cap \Gamma(z)| = b$.
By this and the preliminary comments we find that $G_{yz}$ is as shown in the table.\\
Case $\partial(y,z) \geq 3$: We have $\Gamma(y) \cap \Gamma(z) = \varnothing$ so $|\Gamma_{s-1}(x) \cap \Gamma(y) \cap \Gamma(z)| = 0$ and $|\Gamma_{s+1}(x) \cap \Gamma(y) \cap \Gamma(z)| = 0$.
By this and the preliminary comments we find that $G_{yz}$ is as shown in the table.
\hfill $\Box$\\

\section{The central elements $\Upsilon$, $\Psi$, $\Lambda$}\label{sec:up-psi-lamb}

We continue to discuss the dual polar graph $\Gamma$ from Section \ref{sec:dpg}.
Recall $q$ from Section \ref{sec:ls-dqk} and $b$ from Section \ref{sec:dpg}.
For the rest of the paper $b, q$ are related as follows.
\begin{eqnarray*}
b = q^2.
\end{eqnarray*}
We note that $q$ is nonzero and not a root of unity.

\medskip
\noindent
In Lemma \ref{lem:omega-z-zstar} we discussed the central elements $\Omega, G, G^*$ of the subconstituent algebra $T$.
In this section we introduce three more central elements $\Upsilon$, $\Psi$, $\Lambda$ of $T$.
These central elements will be useful later when we display some $U_q(\mathfrak{sl}_2)$-module structures on the standard module $V$.

\begin{definition} \rm
\cite[Definition~3.1]{kim}
\label{def:displacememnt_scalars}
Let $W$ denote an irreducible $T$-module.
Let $r, t, d$ denote the endpoint, dual endpoint, and diameter of $W$, respectively.
By the {\it displacement of $W$ of the first kind} we mean $r+t+d-D$.
By the {\it displacement of $W$ of the second kind} we mean $r-t$.
\end{definition}

\begin{lemma}
{\rm \cite[Lemma~3.2]{kim}}
Let $W$ denote an irreducible $T$-module. Then the following {\rm (i), (ii)} hold.
\begin{enumerate}
\item[\rm (i)] Let $\mu$ denote the displacement of $W$ of the first kind.  Then $0 \leq \mu \leq D$.
\item[\rm (ii)] Let $\nu$ denote the displacement of $W$ of the second kind.  Then $-D \leq \nu \leq D$.
\end{enumerate}
\end{lemma}

\begin{definition} \rm
\label{def:map_displacement_1}
For an integer $\mu$ $(0 \leq \mu \leq D)$ let $V_\mu$ denote the subspace of $V$ spanned by the irreducible $T$-modules for which $\mu$ is the displacement of the first kind.
Observe that $V_\mu$ is a $T$-module.
By \cite[Lemma~4.4]{ds} we have $V = \sum_{\mu=0}^D V_\mu$ (orthogonal direct sum).
For $0 \leq \mu \leq D$ we define a matrix $\sigma_\mu \in {\rm Mat}_X(\C)$ such that
\begin{eqnarray*}
&&(\sigma_\mu - I)V_\mu = 0,\\
&&\sigma_\mu V_{\mu'} = 0 \quad \mbox{if} \; \mu \neq \mu' \quad (0 \leq \mu' \leq D).
\end{eqnarray*}
In other words $\sigma_\mu$ is the projection from $V$ onto $V_\mu$.
We note that $V_\mu = \sigma_\mu V$.
\end{definition}

\noindent
The following three lemmas are immediate from Definition \ref{def:map_displacement_1}.

\begin{lemma}
\label{lem:upsilon_prop}
The following {\rm (i), (ii)} hold.
\begin{enumerate}
\item[\rm (i)] $I = \sum_{\mu=0}^D \sigma_\mu$.
\item[\rm (ii)] $\sigma_\mu \sigma_{\mu'} = \delta_{\mu \mu'} \sigma_\mu \quad (0 \leq \mu, \mu' \leq D).$
\end{enumerate}
\end{lemma}

\begin{lemma}
\label{lem:displacement1-decomp}
We have
\begin{eqnarray*}
V = \sum_{\mu=0}^D \sigma_\mu V \quad ({\rm orthogonal \ direct \ sum}).
\label{eq:displacement_decomp1}
\end{eqnarray*}
Moreover for $0 \leq \mu \leq D$ the subspace $\sigma_\mu V$ is spanned by the irreducible $T$-modules for which $\mu$ is the displacement of the first kind.
\end{lemma}

\noindent
Recall the set {\rm Feas} from Section \ref{sec:centeroft}.

\begin{lemma}
For $0 \leq \mu \leq D$ we have
\begin{eqnarray*}
\sigma_\mu V = \sum V_{(r,t,d)}
\end{eqnarray*}
where the sum is over all $(r,t,d) \in {\rm Feas}$ such that $r+t+d-D = \mu$.
\end{lemma}

\begin{definition} \rm
\label{def:Upsilon}
Let $\Upsilon$ denote the matrix in ${\rm Mat}_X(\C)$ such that
\begin{eqnarray*}
\Upsilon = \sum_{\mu=0}^D q^\mu \sigma_\mu.
\end{eqnarray*}
\end{definition}

\begin{lemma}
\label{lem:Upsilon-action}
For $0 \leq \mu \leq D$ the matrix $\Upsilon$ acts on $\sigma_\mu V$ as $q^\mu I$.
\end{lemma}

\noindent {\it Proof:} \rm
Immediate from Lemma \ref{lem:upsilon_prop}(ii).
\hfill $\Box$\\

\begin{lemma}
The matrix $\Upsilon$ is invertible and its inverse is
\begin{eqnarray*}
\Upsilon^{-1} = \sum_{\mu=0}^D q^{-\mu} \sigma_\mu.
\end{eqnarray*}
\end{lemma}

\noindent {\it Proof:}
Immediate from Lemma \ref{lem:upsilon_prop}.
\hfill $\Box$\\

\begin{definition} \rm
\label{def:map_displacement_2}
For an integer $\nu$ $(-D \leq \nu \leq D)$ let $V_\nu$ denote the subspace of $V$ spanned by the irreducible $T$-modules for which $\nu$ is the displacement of the second kind.
Observe that $V_\nu$ is a $T$-module.
By \cite[Lemma~4.4]{ds} we have $V = \sum_{\nu=-D}^D V_\nu$ (orthogonal direct sum).
For $-D \leq \nu \leq D$ we define a matrix $\psi_\nu \in {\rm Mat}_X(\C)$ such that
\begin{eqnarray*}
&&(\psi_\nu - I)V_\nu = 0,\\
&&\psi_\nu V_{\nu'} = 0 \quad \mbox{if} \; \nu \neq \nu' \quad (-D \leq \nu' \leq D).
\end{eqnarray*}
In other words $\psi_\nu$ is the projection from $V$ onto $V_\nu$.
We note that $V_\nu = \psi_\nu V$.
\end{definition}

\noindent
The following three lemmas are immediate from Definition \ref{def:map_displacement_2}.

\begin{lemma}
\label{lem:psi_prop}
The following {\rm (i), (ii)} hold.
\begin{enumerate}
\item[\rm (i)] $I = \sum_{\nu=-D}^D \psi_\nu$.
\item[\rm (ii)] $\psi_\nu \psi_{\nu'} = \delta_{\nu \nu'} \psi_\nu \quad (-D \leq \nu, \nu' \leq D).$
\end{enumerate}
\end{lemma}

\begin{lemma}
\label{lem:displacement2-decomp}
We have
\begin{eqnarray*}
V = \sum_{\nu=-D}^D \psi_\nu V \quad  ({\rm orthogonal \ direct \ sum}).
\label{eq:displacement_decomp2}
\end{eqnarray*}
Moreover for $-D \leq \nu \leq D$ the subspace $\psi_\nu V$ is spanned by the irreducible $T$-modules for which $\nu$ is the displacement of the second kind.
\end{lemma}

\begin{lemma}
For $-D \leq \nu \leq D$ we have
\begin{eqnarray*}
\psi_\nu V = \sum V_{(r,t,d)}
\end{eqnarray*}
where the sum is over all $(r,t,d) \in {\rm Feas}$ such that $r-t = \nu$.
\end{lemma}

\begin{definition} \rm
\label{def:Psi}
Let $\Psi$ denote the matrix in ${\rm Mat}_X(\C)$ such that
\begin{eqnarray*}
\Psi = \sum_{\nu=-D}^D q^\nu \psi_\nu.
\end{eqnarray*}
\end{definition}

\begin{lemma}
\label{lem:Psi-action}
For $-D \leq \nu \leq D$ the matrix $\Psi$ acts on $\psi_\nu V$ as $q^\nu I$.
\end{lemma}

\noindent {\it Proof:} \rm
Immediate from Lemma \ref{lem:psi_prop}(ii).
\hfill $\Box$\\

\begin{lemma}
The matrix $\Psi$ is invertible and its inverse is
\begin{eqnarray*}
\Psi^{-1} = \sum_{\nu=-D}^D q^{-\nu} \psi_\nu.
\end{eqnarray*}
\end{lemma}

\noindent {\it Proof:}
Immediate from Lemma \ref{lem:psi_prop}.
\hfill $\Box$\\

\noindent
Earlier we defined the notion of diameter for an irreducible $T$-module.  Using this notion we define some more projections involving $V$.

\begin{definition} \rm
\label{def:map_diameter}
For an integer $d$ $(0 \leq d \leq D)$ let $V_d$ denote the subspace of $V$ spanned by the irreducible $T$-modules of diameter $d$.
Observe that $V_d$ is a $T$-module,
and $V = \sum_{d=0}^D V_d$ (orthogonal direct sum).
For $0 \leq d \leq D$ we define a matrix $\rho_d \in {\rm Mat}_X(\C)$ such that
\begin{eqnarray*}
&&(\rho_d - I)V_d = 0,\\
&&\rho_d V_{d'} = 0 \quad \mbox{if} \; d \neq d' \quad (0 \leq d' \leq D).
\end{eqnarray*}
In other words $\rho_d$ is the projection from $V$ onto $V_d$.
We note that $V_d = \rho_d V$.
\end{definition}

\noindent
The following three lemmas are immediate from Definition \ref{def:map_diameter}.

\begin{lemma}
\label{lem:lambda_prop}
The following {\rm (i), (ii)} hold.
\begin{enumerate}
\item[\rm (i)] $I = \sum_{d=0}^D \rho_d$.
\item[\rm (ii)] $\rho_d \rho_{d'} = \delta_{d d'} \rho_d \quad (0 \leq d, d' \leq D).$
\end{enumerate}
\end{lemma}

\begin{lemma}
\label{lem:t_diam_decomp}
We have
\begin{eqnarray}
V = \sum_{d=0}^D \rho_d V \quad  ({\rm orthogonal \ direct \ sum}).
\label{eq:diameter_decomp}
\end{eqnarray}
Moreover for $0 \leq d \leq D$ the subspace $\rho_d V$ is spanned by the irreducible $T$-modules of diameter $d$.
\end{lemma}

\begin{lemma}
\label{lem:t_homeg_d}
For $0 \leq d \leq D$ we have
\begin{eqnarray*}
\rho_d V = \sum V_{(r,t,d)}
\end{eqnarray*}
where the sum is over all integers $r, t$ such that $(r,t,d) \in {\rm Feas}$.
\end{lemma}

\begin{definition} \rm
\label{def:Lambda}
Let $\Lambda$ denote the matrix in ${\rm Mat}_X(\C)$ such that
\begin{eqnarray*}
\Lambda = \sum_{d=0}^D \frac{q^{d+1}+q^{-d-1}}{(q-q^{-1})^2} \rho_d.
\end{eqnarray*}
\end{definition}

\begin{lemma}
\label{lem:Lambda_scalar}
For $0 \leq d \leq D$ the matrix $\Lambda$ acts on $\rho_d V$ as
\begin{eqnarray}
\frac{q^{d+1}+q^{-d-1}}{(q-q^{-1})^2} I.\label{eq:Lambda_scalar}
\end{eqnarray}
\end{lemma}

\noindent {\it Proof:} \rm
Immediate from Lemma \ref{lem:lambda_prop}(ii).
\hfill $\Box$\\

\begin{lemma} {\rm \cite[Lemma~12.1]{drgqtet}}
\label{lem:upl_central}
The matrices $\Upsilon,\Psi,\Lambda$ are central elements in $T$.
\end{lemma}

\noindent
The central elements $\Omega, G$, $G^*$ of $T$ from Lemma \ref{lem:omega-z-zstar} are related to $\Upsilon, \Psi, \Lambda$ as follows.

\begin{proposition}
\label{prop:omegazzstar-upsilonpsilambda}
Let $\xi, \zeta$ denote the scalars from Lemma \ref{lem:dpgdualeval},
and let $\gamma,\varrho$ denote the scalars from Lemma \ref{lem:dpgconsts}.
Then $\Omega, G, G^*$ can be expressed in terms of $\Upsilon, \Psi, \Lambda$ as follows.
\begin{enumerate}
\item[\rm (i)] $\displaystyle \Omega = \xi (q^{-2}-1)(q^{2e}\Upsilon^{-2} - \Psi^{-2}) - 2\gamma\zeta I$.
\item[\rm (ii)] $\displaystyle G = \xi q^{-2} (1-q^{2e})(q^{2e}\Upsilon^{-2}-\Psi^{-2})
-\xi q^{D+2e-5}(q^2-1)^2(q^2+1)\Upsilon^{-1}\Psi^{-1}\Lambda - \varrho \zeta I$.
\item[\rm (iii)] $\displaystyle G^* = \zeta \xi (1-q^{-2}) (q^{2e}\Upsilon^{-2} - \Psi^{-2}) + \gamma\zeta^2 I$.
\end{enumerate}
\end{proposition}

\noindent {\it Proof:} \rm
(i) By construction $V$ is a direct sum of irreducible $T$-modules.
Let $W$ denote an irreducible $T$-module in the sum.
It suffices to show that for the equation in (i) the two sides agree on $W$.
Let $\mu, \nu$ denote the displacement of $W$ of the first kind and second kind, respectively.
By Lemma \ref{lem:omega-z-zstar} the element $\Omega$ acts on $W$ as $\omega(W)I$.
By Lemma \ref{lem:displacement1-decomp} we have $W \subseteq \sigma_{\mu}V$.
By this and Lemma \ref{lem:Upsilon-action} the matrix $\Upsilon$ acts on $W$ as $q^{\mu}I$.
By Lemma \ref{lem:displacement2-decomp} we have $W \subseteq \psi_{\nu}V$.
By this and Lemma \ref{lem:Psi-action} the matrix $\Psi$ acts on $W$ as $q^{\nu}I$.
By these comments, Lemma \ref{lem:omega_eta_etastar}(i) and Definition \ref{def:displacememnt_scalars}, the two sides agree on $W$.
The result follows.\\
(ii) Similar to the proof of (i).\\
(iii) Combine Lemma \ref{lem:zstar} and (i).
\hfill $\Box$\\

\section{Irreducible $T$-modules and Leonard systems of dual $q$-Krawtchouk type}

We continue to discuss the dual polar graph $\Gamma$ from Section \ref{sec:dpg}.
In Lemma \ref{lem:irred-t-module-lp} we obtained a Leonard system on each irreducible $T$-module.
In this section we show that this Leonard system has dual $q$-Krawtchouk type.

\begin{theorem}
\label{thm:ls-dqkc-on-w}
Let $W$ denote an irreducible $T$-module.
Let $r, t, d$ denote the endpoint, dual endpoint and diameter of $W$, respectively.
Let $\Phi$ denote the corresponding Leonard system on $W$ from Lemma \ref{lem:irred-t-module-lp}.
Then $\Phi$ has dual $q$-Krawtchouk type.
Let $h,h^*,\kappa,\kappa^*,\upsilon$ denote the parameters corresponding to $\Phi$ from Definition \ref{def:dualqkrawtchouk2}.
Let $\zeta,\xi$ denote the scalars from Lemma \ref{lem:dpgdualeval}.
Then
\begin{eqnarray}
h = \frac{1-q^{2e}}{q^2-1},
\qquad h^* = \zeta, 
\qquad \kappa = \frac{q^{2e+2D-2t-d}}{q^2-1},
\qquad \kappa^* = \xi q^{-2r-d},
\qquad \upsilon = -\frac{q^{2t+d}}{q^2-1}.
\label{eq:ls-qkc-params}
\end{eqnarray}
\end{theorem}

\noindent {\it Proof:} \rm
By Lemma \ref{lem:irred-t-module-lp} the Leonard system $\Phi$ and the $T$-module $W$ have the same intersection matrix and dual intersection matrix.
We show that $\Phi$ has dual $q$-Krawtchouk type by verifying that its parameter array satisfies (\ref{eq:lp-eval})--(\ref{eq:lp-split2}).
By Lemma \ref{lem:dpgeval} the eigenvalue sequence $\{\theta_{t+i}\}_{i=0}^d$ has the form (\ref{eq:lp-eval}) with scalars $h,\kappa,\upsilon$ given in (\ref{eq:ls-qkc-params}).
By Lemma \ref{lem:dpgdualeval} the dual eigenvalue sequence $\{\theta^*_{r+i}\}_{i=0}^d$ has the form (\ref{eq:lp-dualeval}) with scalars $h^*,\kappa^*$ given in (\ref{eq:ls-qkc-params}).
By Lemma \ref{lem:int_mat_entries2} we have $\varphi_1 = (\theta_r^*-\theta_{r+1}^*)(a_0(W)-\theta_t)$.
On the right-hand side evaluate the eigenvalue using Lemma \ref{lem:dpgeval}, evaluate the dual eigenvalues using Lemma \ref{lem:dpgdualeval}, evaluate $a_0(W)$ using Lemma \ref{lem:dpgintw}, and simplify the result to get
\begin{eqnarray}
	\varphi_1 = -\xi (q^{2d}-1)q^{2(e+D-d-t-r-1)}. \label{eq:varphi1}
\end{eqnarray}
By (PA4) we have $\phi_1 = \varphi_1 + (\theta_{r+1}^*-\theta_{r}^*)(\theta_{t+d}-\theta_t)$.
On the right-hand side evaluate the eigenvalues using Lemma \ref{lem:dpgeval}, evaluate the dual eigenvalues using Lemma \ref{lem:dpgdualeval}, evaluate $\varphi_1$ using (\ref{eq:varphi1}), and simplify the result to get
\begin{eqnarray}
	\phi_1 = \xi (q^{2d}-1)q^{2(t-r-1)}. \label{eq:phi1}
\end{eqnarray}
On the right-hand side of (PA3) evaluate the eigenvalues using Lemma \ref{lem:dpgeval}, evaluate the dual eigenvalues using Lemma \ref{lem:dpgdualeval}, evaluate $\phi_1$ using (\ref{eq:phi1}), and simplify the result to get
\begin{eqnarray}
\varphi_i = \frac{\xi}{q^2-1}(q^{2(d+1)}-q^2)(q^{-2i}-1)q^{2(e+D-d-t-r-i)} \qquad (1 \leq i \leq d). \label{eq:varphii}
\end{eqnarray}
On the right-hand side of (PA4) evaluate the eigenvalues using Lemma \ref{lem:dpgeval}, evaluate the dual eigenvalues using Lemma \ref{lem:dpgdualeval}, evaluate $\varphi_1$ using (\ref{eq:varphi1}), and simplify the result to get
\begin{eqnarray}
\phi_i = \frac{\xi}{q^2-1}(q^{2(d+1)}-q^2)(q^{-2i}-1)q^{2(e+D-d-t-r-i)} \qquad (1 \leq i \leq d). \label{eq:phii}
\end{eqnarray}
Using (\ref{eq:ls-qkc-params}) it is routine to verify that (\ref{eq:varphii}), (\ref{eq:phii}) satisfy (\ref{eq:lp-split1}), (\ref{eq:lp-split2}), respectively.  Therefore the Leonard system $\Phi$ has dual $q$-Krawtchouk type.
\hfill $\Box$\\

\section{Two $U_q(\mathfrak{sl}_2)$-module structures on the standard module $V$}

We continue to discuss the dual polar graph $\Gamma$ from Section \ref{sec:dpg}.
In this section we display two $U_q(\mathfrak{sl}_2)$-module structures on the standard module $V$.
Then we show how the two $U_q(\mathfrak{sl}_2)$-module structures are related.

\begin{lemma}
\label{lem:uqsl2onw}
Let $W$ denote an irreducible $T$-module.
Let $r, t, d$ denote the endpoint, dual endpoint and diameter of $W$, respectively.
Let $h,h^*,\kappa,\kappa^*,\upsilon$ denote the corresponding parameters from {\rm (\ref{eq:ls-qkc-params})}.
Then there exists a unique $U_q(\mathfrak{sl}_2)$-module structure on $W$ such that on $W$,
\begin{align}
A & = h1 + \kappa x + \upsilon y, \label{eq:a_w}\\
A^* & = h^*1 + \kappa^*z. \label{eq:astar_w}
\end{align}
Moreover the $U_q(\mathfrak{sl}_2)$-module $W$ is isomorphic to $L(d,1)$.
\end{lemma}

\noindent {\it Proof:} \rm
Combine Theorem \ref{thm:ls-dqkc-on-w} and Theorem \ref{thm:uqsl2ondualqkrawtchouk} taking $\epsilon = 1$.
\hfill $\Box$\\

\begin{lemma}
\label{lem:uqsl2onw_2}
Let $W$ denote an irreducible $T$-module.
Let $r, t, d$ denote the endpoint, dual endpoint and diameter of $W$, respectively.
Let $h,h^*,\kappa,\kappa^*,\upsilon$ denote the corresponding parameters from {\rm (\ref{eq:ls-qkc-params})}.
Then there exists a unique $U_q(\mathfrak{sl}_2)$-module structure on $W$ such that on $W$,
\begin{align}
A & = h1 + \kappa y + \upsilon x,\nonumber\\%\label{eq:a_w_2}\\
A^* & = h^*1 + \kappa^*z.\nonumber%\label{eq:astar_w_2}
\end{align}
Moreover the $U_q(\mathfrak{sl}_2)$-module $W$ is isomorphic to $L(d,1)$.
\end{lemma}

\noindent {\it Proof:} \rm
Combine Theorem \ref{thm:ls-dqkc-on-w} and Theorem \ref{thm:uqsl2ondualqkrawtchouk2} taking $\epsilon = 1$.
\hfill $\Box$\\

\begin{theorem}
\label{thm:uqsl2ondpg}
There exists a $U_q(\mathfrak{sl}_2)$-module structure on $V$ such that on $V$,
\begin{align}
A & = h1 + \kappa\Upsilon^{-1}\Psi x + \upsilon\Upsilon\Psi^{-1} y, \label{eq:a_v}\\
A^* & = h^*1 + \kappa^*\Upsilon^{-1}\Psi^{-1} z, \label{eq:astar_v}
\end{align}
where
\begin{eqnarray}
h = \frac{1-q^{2e}}{q^2-1},
\qquad h^* = \zeta,
\qquad \kappa = \frac{q^{2e+D}}{q^2-1},
\qquad \kappa^* = \xi q^{-D},
\qquad \upsilon = -\frac{q^{D}}{q^2-1}.
\label{eq:dpg-params}
\end{eqnarray}
\end{theorem}

\noindent {\it Proof:} \rm
By construction $V$ is a direct sum of irreducible $T$-modules.
Let $W$ denote an irreducible $T$-module in the sum.
It suffices to show that (\ref{eq:a_v}), (\ref{eq:astar_v}) hold on $W$.
Let $\mu, \nu$ denote the displacement of $W$ of the first kind and second kind, respectively.
Consider the $U_q(\mathfrak{sl}_2)$-module structure on $W$ from Lemma \ref{lem:uqsl2onw}.
Writing (\ref{eq:a_w}), (\ref{eq:astar_w}) in terms of (\ref{eq:dpg-params}) and $\mu, \nu$ we get
\begin{align}
A & = h1 + \kappa q^{\nu-\mu} x + \upsilon q^{\mu-\nu} y,\label{eq:awv-temp}\\
A^* & = h^*1 + \kappa^* q^{-\mu-\nu} z,\label{eq:astarwv-temp}
\end{align}
where $h,h^*,\kappa,\kappa^*,\upsilon$ are from (\ref{eq:dpg-params}).
By Lemma \ref{lem:displacement1-decomp} we have $W \subseteq \sigma_{\mu}V$.
By this and Lemma \ref{lem:Upsilon-action} the matrix $\Upsilon$ acts on $W$ as $q^{\mu}I$.
By Lemma \ref{lem:displacement2-decomp} we have $W \subseteq \psi_{\nu}V$.
By this and Lemma \ref{lem:Psi-action} the matrix $\Psi$ acts on $W$ as $q^{\nu}I$.
By these comments and (\ref{eq:awv-temp}), (\ref{eq:astarwv-temp}) we find that (\ref{eq:a_v}), (\ref{eq:astar_v}) hold on $W$.
Therefore (\ref{eq:a_v}), (\ref{eq:astar_v}) hold on $V$.
\hfill $\Box$\\

\begin{theorem}
\label{thm:uqsl2ondpg_2}
There exists a $U_q(\mathfrak{sl}_2)$-module structure on $V$ such that on $V$,
\begin{align*}
A & = h1 + \kappa\Upsilon^{-1}\Psi y + \upsilon\Upsilon\Psi^{-1} x,\\%\label{eq:a_v_2}\\
A^* & = h^*1 + \kappa^*\Upsilon^{-1}\Psi^{-1} z,%\label{eq:astar_v_2}
\end{align*}
where the scalars $h,h^*,\kappa,\kappa^*,\upsilon$ are from {\rm (\ref{eq:dpg-params})}.
\end{theorem}

\noindent {\it Proof:} \rm
Similar to the proof of Theorem \ref{thm:uqsl2ondpg} but use Lemma \ref{lem:uqsl2onw_2} instead of Lemma \ref{lem:uqsl2onw}.
\hfill $\Box$\\

\begin{lemma}
\label{lem:xyz_in_bbstar}
For the $U_q(\mathfrak{sl}_2)$-module structure on $V$ from Theorem \ref{thm:uqsl2ondpg}, the actions of $x,y,z$ on $V$ are given by
\begin{align}
x & = \frac{\Upsilon\Psi^{-1}(qB-q^{-1}B^*BB^{*-1})}{\kappa q^{-1}(q^2-q^{-2})}
+ \frac{q^{-1}\kappa^*\Upsilon^{-1}\Psi^{-1}B^{*-1}}{q+q^{-1}}
- \frac{q\kappa^*\upsilon\Upsilon\Psi^{-3}B^{*-1}}{\kappa(q+q^{-1})}\label{eq:x1_bbstar}\\
& = \frac{\Upsilon\Psi^{-1}(qB^{*-1}BB^*-q^{-1}B)}{\kappa q(q^2-q^{-2})}
+ \frac{q\kappa^*\Upsilon^{-1}\Psi^{-1}B^{*-1}}{q+q^{-1}}
- \frac{q^{-1}\kappa^*\upsilon\Upsilon\Psi^{-3}B^{*-1}}{\kappa(q+q^{-1})},\nonumber\\%\label{eq:x2_bbstar}\\
y & = \frac{\Upsilon^{-1}\Psi(qB-q^{-1}B^{*-1}BB^*)}{\upsilon q^{-1}(q^2-q^{-2})}
+ \frac{q^{-1}\kappa^*\Upsilon^{-1}\Psi^{-1}B^{*-1}}{q+q^{-1}}
- \frac{q\kappa^*\kappa\Upsilon^{-3}\Psi B^{*-1}}{\upsilon(q+q^{-1})}\label{eq:y1_bbstar}\\
& = \frac{\Upsilon^{-1}\Psi(qB^*BB^{*-1}-q^{-1}B)}{\upsilon q(q^2-q^{-2})}
+ \frac{q\kappa^*\Upsilon^{-1}\Psi^{-1}B^{*-1}}{q+q^{-1}}
- \frac{q^{-1}\kappa^*\kappa\Upsilon^{-3}\Psi B^{*-1}}{\upsilon(q+q^{-1})},\nonumber\\%\label{eq:y2_bbstar}\\
z & = \kappa^{*-1}\Upsilon\Psi B^*,\label{eq:z_bbstar}
	\end{align}
where $B = A - hI$, $B^* = A^* - h^*I$ and the scalars $h,h^*,\kappa,\kappa^*,\upsilon$ are from {\rm (\ref{eq:dpg-params})}.
\end{lemma}

\noindent {\it Proof:} \rm
Similar to the proof of Lemma \ref{lem:xyz_aastar} but use Theorem \ref{thm:uqsl2ondpg} instead of Lemma \ref{cor:aastar-in-xyz}.
\hfill $\Box$\\

\begin{lemma}
\label{lem:xyzprime_in_bbstar}
For the $U_q(\mathfrak{sl}_2)$-module structure on $V$ from Theorem \ref{thm:uqsl2ondpg_2}, the actions of $x,y,z$ on $V$ are given by
\begin{align}
x & = \frac{\Upsilon^{-1}\Psi(qB-q^{-1}B^*BB^{*-1})}{\upsilon q^{-1}(q^2-q^{-2})}
+ \frac{q^{-1}\kappa^*\Upsilon^{-1}\Psi^{-1}B^{*-1}}{q+q^{-1}}
- \frac{q\kappa^*\kappa\Upsilon^{-3}\Psi B^{*-1}}{\upsilon(q+q^{-1})}\label{eq:xprime1_bbstar}\\
& = \frac{\Upsilon^{-1}\Psi(qB^{*-1}BB^*-q^{-1}B)}{\upsilon q(q^2-q^{-2})}
+ \frac{q\kappa^*\Upsilon^{-1}\Psi^{-1}B^{*-1}}{q+q^{-1}}
- \frac{q^{-1}\kappa^*\kappa\Upsilon^{-3}\Psi B^{*-1}}{\upsilon (q+q^{-1})},\nonumber\\%\label{eq:xprime2_bbstar}\\
y & = \frac{\Upsilon\Psi^{-1}(qB-q^{-1}B^{*-1}BB^*)}{\kappa q^{-1}(q^2-q^{-2})}
+ \frac{q^{-1}\kappa^*\Upsilon^{-1}\Psi^{-1}B^{*-1}}{q+q^{-1}}
- \frac{q\kappa^*\upsilon\Upsilon\Psi^{-3} B^{*-1}}{\kappa(q+q^{-1})}\label{eq:yprime1_bbstar}\\
& = \frac{\Upsilon\Psi^{-1}(qB^*BB^{*-1}-q^{-1}B)}{\kappa q(q^2-q^{-2})}
+ \frac{q\kappa^*\Upsilon^{-1}\Psi^{-1}B^{*-1}}{q+q^{-1}}
- \frac{q^{-1}\kappa^*\upsilon\Upsilon\Psi^{-3} B^{*-1}}{\kappa(q+q^{-1})},\nonumber\\%\label{eq:yprime2_bbstar}\\
z & = \kappa^{*-1}\Upsilon\Psi B^*,\label{eq:zprime_bbstar}
	\end{align}
where $B = A - hI$, $B^* = A^* - h^*I$ and the scalars $h,h^*,\kappa,\kappa^*,\upsilon$ are from {\rm (\ref{eq:dpg-params})}.
\end{lemma}

\noindent {\it Proof:} \rm
Similar to the proof of Lemma \ref{lem:xyz_aastar} but use Theorem \ref{thm:uqsl2ondpg_2} instead of Lemma \ref{cor:aastar-in-xyz}.
\hfill $\Box$\\

\noindent
We draw several corollaries from Lemma \ref{lem:xyz_in_bbstar} and Lemma \ref{lem:xyzprime_in_bbstar}.

\begin{corollary}
The $U_q(\mathfrak{sl}_2)$-module structure on $V$ from Theorem \ref{thm:uqsl2ondpg} is unique.
\end{corollary}

\noindent {\it Proof:} \rm
By Lemma \ref{lem:xyz_in_bbstar} and since $x, y, z^{\pm 1}$ generate $U_q(\mathfrak{sl}_2)$.
\hfill $\Box$\\

\begin{corollary}
The $U_q(\mathfrak{sl}_2)$-module structure on $V$ from Theorem \ref{thm:uqsl2ondpg_2} is unique.
\end{corollary}

\noindent {\it Proof:} \rm
By Lemma \ref{lem:xyzprime_in_bbstar} and since $x, y, z^{\pm 1}$ generate $U_q(\mathfrak{sl}_2)$.
\hfill $\Box$\\

\begin{corollary}
\label{cor:uqsl2_t}
For the $\C$-algebra homomorphism $U_q(\mathfrak{sl}_2) \rightarrow {\rm End}(V)$ induced by the $U_q(\mathfrak{sl}_2)$-module structure on $V$ from Theorem \ref{thm:uqsl2ondpg},
the image is contained in $T$.
\end{corollary}

\noindent {\it Proof:} \rm
By Lemma \ref{lem:upl_central} and Lemma \ref{lem:xyz_in_bbstar}.
\hfill $\Box$\\

\begin{corollary}
\label{cor:uqsl2_t_2}
For the $\C$-algebra homomorphism $U_q(\mathfrak{sl}_2) \rightarrow {\rm End}(V)$ induced by the $U_q(\mathfrak{sl}_2)$-module structure on $V$ from Theorem \ref{thm:uqsl2ondpg_2},
the image is contained in $T$.
\end{corollary}

\noindent {\it Proof:} \rm
By Lemma \ref{lem:upl_central} and Lemma \ref{lem:xyzprime_in_bbstar}.
\hfill $\Box$\\

\noindent
In Theorem \ref{thm:uqsl2ondpg} and Theorem \ref{thm:uqsl2ondpg_2} we gave $U_q(\mathfrak{sl}_2)$-actions on the standard module $V$.
We now describe how these actions look on each irreducible $T$-module.

\begin{lemma}
\label{lem:uqsl2_v_to_w}
Let $W$ denote an irreducible $T$-module.
Then $W$ is invariant under the $U_q(\mathfrak{sl}_2)$-action from Theorem \ref{thm:uqsl2ondpg}.
Moreover the action of $U_q(\mathfrak{sl}_2)$ on $W$ coincides with the $U_q(\mathfrak{sl}_2)$-action from Lemma \ref{lem:uqsl2onw}.
\end{lemma}

\noindent {\it Proof:} \rm
The first assertion follows from Corollary \ref{cor:uqsl2_t} .
We now verify the second assertion.
By Lemma \ref{lem:uqsl2onw} it suffices to show that (\ref{eq:a_w}), (\ref{eq:astar_w}) hold on $W$.
Let $\mu, \nu$ denote the displacement of $W$ of the first kind and second kind, respectively.
By Lemma \ref{lem:displacement1-decomp} we have $W \subseteq \sigma_{\mu}V$.
By this and Lemma \ref{lem:Upsilon-action} the matrix $\Upsilon$ acts on $W$ as $q^{\mu}I$.
By Lemma \ref{lem:displacement2-decomp} we have $W \subseteq \psi_{\nu}V$.
By this and Lemma \ref{lem:Psi-action} the matrix $\Psi$ acts on $W$ as $q^{\nu}I$.
By applying these comments, Definition \ref{def:displacememnt_scalars} and (\ref{eq:dpg-params}) to (\ref{eq:a_v}), (\ref{eq:astar_v}) we find that (\ref{eq:a_w}), (\ref{eq:astar_w}) hold on $W$.
\hfill $\Box$\\

\begin{lemma}
\label{lem:uqsl2_v_to_w_2}
Let $W$ denote an irreducible $T$-module.
Then $W$ is invariant under the $U_q(\mathfrak{sl}_2)$-action from Theorem \ref{thm:uqsl2ondpg_2}.
Moreover the action of $U_q(\mathfrak{sl}_2)$ on $W$ coincides with the $U_q(\mathfrak{sl}_2)$-action from Lemma \ref{lem:uqsl2onw_2}.
\end{lemma}

\noindent {\it Proof:} \rm
Similar to the proof of Lemma \ref{lem:uqsl2_v_to_w} but use Lemma \ref{lem:uqsl2onw_2} and Corollary \ref{cor:uqsl2_t_2} instead of Lemma \ref{lem:uqsl2onw} and Corollary \ref{cor:uqsl2_t}.
\hfill $\Box$\\

\noindent
We finish this section with a comment describing how the two $U_q(\mathfrak{sl}_2)$-module structures on $V$ from Theorem \ref{thm:uqsl2ondpg} and Theorem \ref{thm:uqsl2ondpg_2} are related.

\begin{theorem}
Consider the table below.
In the first column the three displayed elements each induces an element in ${\rm End}(V)$ using the $U_q(\mathfrak{sl}_2)$-module structure from Theorem \ref{thm:uqsl2ondpg_2}.
In the second column the three displayed elements each induces an element in ${\rm End}(V)$ using the $U_q(\mathfrak{sl}_2)$-module structure from Theorem \ref{thm:uqsl2ondpg}.
For each row the two elements induce the same element of ${\rm End}(V)$.\\
\begin{center}
\begin{tabular}{c|c}
$U_q(\mathfrak{sl}_2)$-module structure & $U_q(\mathfrak{sl}_2)$-module structure\\
from Theorem \ref{thm:uqsl2ondpg_2} & from Theorem \ref{thm:uqsl2ondpg}\\
\hline
$z$ & $z$ \\
$x$ & $-q^{2e} \Upsilon^{-2} \Psi^2 x + (1 + q^{2e} \Upsilon^{-2} \Psi^2) z^{-1}$ \\
$y$ & $-q^{-2e} \Upsilon^{2} \Psi^{-2} y + (1 + q^{-2e} \Upsilon^{2} \Psi^{-2}) z^{-1}$ \\
\end{tabular}
\end{center}
\end{theorem}

\noindent {\it Proof:} \rm
The first row is immediate from (\ref{eq:z_bbstar}), (\ref{eq:zprime_bbstar}).  To prove the second row evaluate $x$ on the left using (\ref{eq:xprime1_bbstar}), and evaluate $x, z$ on the right using  (\ref{eq:x1_bbstar}), (\ref{eq:z_bbstar}).  To prove the last row evaluate $y$ on the left using (\ref{eq:yprime1_bbstar}), and evaluate $y, z$ on the right using  (\ref{eq:y1_bbstar}), (\ref{eq:z_bbstar}).
\hfill $\Box$\\

\section{Two homomorphisms $U_q(\mathfrak{sl}_2) \rightarrow T$}

We continue to discuss the dual polar graph $\Gamma$ from Section \ref{sec:dpg}.
In Theorem \ref{thm:uqsl2ondpg} and Theorem \ref{thm:uqsl2ondpg_2} we displayed two $U_q(\mathfrak{sl}_2)$-module structures on the standard module $V$.
In this section we show how these two $U_q(\mathfrak{sl}_2)$-module structures are related to $T$.
By Corollary \ref{cor:uqsl2_t} the $U_q(\mathfrak{sl}_2)$-module structure from Theorem \ref{thm:uqsl2ondpg} induces a $\C$-algebra homomorphism $U_q(\mathfrak{sl}_2) \rightarrow T$.
By Corollary \ref{cor:uqsl2_t_2} the $U_q(\mathfrak{sl}_2)$-module structure from Theorem \ref{thm:uqsl2ondpg_2} induces a $\C$-algebra homomorphism $U_q(\mathfrak{sl}_2) \rightarrow T$.
For either of the two $\C$-algebra homomorphisms, let $U$ denote the image.
We show that $T$ is generated by $U$ together with the elements $\Upsilon^{\pm 1}, \Psi^{\pm 1}$ where $\Upsilon$ is from Definition \ref{def:Upsilon} and $\Psi$ is from Definition \ref{def:Psi}.

\medskip
\noindent
Recall below Lemma \ref{lem:casimir-scalar-distinct} we discussed the homogeneous components for a $U_q(\mathfrak{sl}_2)$-module.
We now consider the homogeneous components for the $U_q(\mathfrak{sl}_2)$-module structure on $V$ from either Theorem \ref{thm:uqsl2ondpg} or Theorem \ref{thm:uqsl2ondpg_2}.

\begin{lemma}
\label{lem:ut_homog_comp_d}
Consider the $U_q(\mathfrak{sl}_2)$-module structure on $V$ from either Theorem \ref{thm:uqsl2ondpg} or Theorem \ref{thm:uqsl2ondpg_2}.
Then
\begin{eqnarray}
V_{d,-1} = 0, \qquad V_{d,1} = \rho_d V \qquad (0 \leq d \leq D).\label{eq:vd1-rhod}
\end{eqnarray}
In other words, in the sum {\rm (\ref{eq:diameter_decomp})} the summands are the homogeneous components of the $U_q(\mathfrak{sl}_2)$-module $V$.
\end{lemma}

\noindent {\it Proof:}
First assume the $U_q(\mathfrak{sl}_2)$-module structure on $V$ is from Theorem \ref{thm:uqsl2ondpg}.
Recall from Lemma \ref{lem:t_diam_decomp} that $\rho_d V$ is spanned by the irreducible $T$-modules of diameter $d$.
By Lemma \ref{lem:uqsl2_v_to_w} and Lemma \ref{lem:uqsl2onw} each irreducible $T$-module of diameter $d$ is a $U_q(\mathfrak{sl}_2)$-module isomorphic to $L(d,1)$.
By these comments,
\begin{eqnarray}
\rho_d V \subseteq V_{d,1}.\label{eq:rhodv_vd1}
\end{eqnarray}
By summing (\ref{eq:rhodv_vd1}) over $0 \leq d \leq D$ and comparing the result to (\ref{eq:diameter_decomp}) we have
\begin{eqnarray}
V = \sum_{d=0}^D V_{d,1} \qquad ({\rm direct \ sum}).\label{eq:v_decomp_vd1}
\end{eqnarray}
By comparing (\ref{eq:v_decomp_vd1}) to (\ref{eq:u_homeg_decomp}) with $M = V$ we have (\ref{eq:vd1-rhod}).
We are now done for the case in which the module structure is from Theorem \ref{thm:uqsl2ondpg}.
For the case in which the module structure is from Theorem \ref{thm:uqsl2ondpg_2}, the argument is similar using Lemma \ref{lem:uqsl2onw_2} and Lemma \ref{lem:uqsl2_v_to_w_2} instead of Lemma \ref{lem:uqsl2onw} and Lemma \ref{lem:uqsl2_v_to_w}.
\hfill $\Box$\\

\noindent
In the beginning of this section we discussed two $\C$-algebra homomorphisms $U_q(\mathfrak{sl}_2) \rightarrow T$.
We now discuss these homomorphisms.

\begin{theorem}
For either of our two $\C$-algebra homomorphisms, let $U$ denote the image.
Then the algebra $T$ is generated by $U$ together with the elements $\Upsilon^{\pm 1}, \Psi^{\pm 1}$ where $\Upsilon$ is from Definition \ref{def:Upsilon} and $\Psi$ is from Definition \ref{def:Psi}.
\end{theorem}

\noindent {\it Proof:} \rm
By Theorem \ref{thm:uqsl2ondpg}, Theorem \ref{thm:uqsl2ondpg_2}, and since $A,A^*$ generate $T$.
\hfill $\Box$\\

\noindent
We finish this section with a comment.  Recall the Casimir element $\Delta$ of $U_q(\mathfrak{sl}_2)$ from Definition \ref{def:casimir}.

\begin{lemma}
For the $U_q(\mathfrak{sl}_2)$-module structure on $V$ from either Theorem \ref{thm:uqsl2ondpg} or Theorem \ref{thm:uqsl2ondpg_2},
the Casimir element $\Delta$ acts on $V$ as the element $\Lambda$ from Definition \ref{def:Lambda}.
\end{lemma}

\noindent {\it Proof:}
By (\ref{eq:diameter_decomp}), for $0 \leq d \leq D$ it suffices to show that $\Delta, \Lambda$ agree on $\rho_d V$.
%Let $d$ be given.
Recall from Lemma \ref{lem:Lambda_scalar} that the element $\Lambda$ acts on $\rho_d V$ as (\ref{eq:Lambda_scalar}).
By Lemma \ref{lem:ut_homog_comp_d} and Lemma \ref{lem:casimir-scalar}, the element $\Delta$ acts on $\rho_d V$ as (\ref{eq:Lambda_scalar}).
Therefore they agree on $\rho_d V$.
The result follows.
\hfill $\Box$\\

\section{The matrices $L, F, R, K$}\label{sec:lfrk}

We continue to discuss the dual polar graph $\Gamma$ from Section \ref{sec:dpg}.
In this section we define some nice matrices that generate $T$ and find relations among them.

\begin{definition} \rm
\label{def:rfl}
We define the matrices $L, F, R$ in $T$ by%\in {\rm Mat}_X(\C)$ by
\begin{align}
L &= \sum_{i=1}^{D} E^*_{i-1} A E^*_{i},\label{eq:thel}\\
F &= \sum_{i=0}^{D} E^*_{i} A E^*_{i},\label{eq:thef}\\
R &= \sum_{i=0}^{D-1} E^*_{i+1} A E^*_{i}.\label{eq:ther}
\end{align}
We call $L$ (resp. $F$) (resp. $R$) the {\em lowering matrix} (resp. {\em flattening matrix}) (resp. {\em raising matrix}).
\end{definition}

\noindent
Observe that $L^t = R$ and $F^t = F$.  Moreover $A = L + F + R$.

\begin{lemma}
\label{lem:rfl_temp2}
The following {\rm (i)--(iii)} hold.
\begin{enumerate}
\item[\rm (i)] $L E^*_i = E^*_{i-1} L E^* _i = E^*_{i-1} L = E^*_{i-1} A E^*_i \qquad (1 \leq i \leq D), \qquad LE_0^* = 0$.
\item[\rm (ii)] $F E^*_i = E^*_i F E^* _i = E^*_i F = E^*_i A E^*_i \qquad (0 \leq i \leq D)$.
\item[\rm (iii)] $R E^*_i = E^*_{i+1} R E^* _i = E^*_{i+1} R = E^*_{i+1} A E^*_i \qquad (0 \leq i \leq D-1), \qquad RE_D^* = 0$.
\end{enumerate}
\end{lemma}

\noindent {\it Proof:} \rm
Routine verification using Definition \ref{def:rfl} and $E^*_i E^*_j = \delta_{ij}E^*_i$ $(0 \leq i, j \leq D)$.
\hfill $\Box$\\

\begin{proposition}
\label{prop:rfl_reln}
The following {\rm (i), (ii)} hold.
\begin{enumerate}
\item[\rm (i)] $LF - q^2 FL = (q^{2e}-1) L$.
\item[\rm (ii)] $FR - q^2 RF = (q^{2e}-1) R$.
\end{enumerate}
\end{proposition}

\noindent {\it Proof:} \rm
(i) Let $y, z$ denote vertices in $X$.
Let $i = \partial(x,z)$.
For the equation in (i) we show that the $(y,z)$-entries of both sides are equal.
By construction
\begin{align*}
L_{yz}
&=
\begin{cases}
1 &\mbox{if } \partial(y,z) = 1 \mbox{ and } \partial(x,y) = i-1;\\
0 &\mbox{otherwise.}
\end{cases}
\end{align*}
By Definition \ref{def:rfl} and Lemma \ref{lem:yz-entry},
\begin{align*}
(LF)_{yz}
&=
\begin{cases}
0 &\mbox{if } \partial(x,y) \neq i-1;\\
|\Gamma_i(x) \cap \Gamma(y) \cap \Gamma(z)| &\mbox{if } \partial(x,y) = i-1,
\end{cases}\\
(FL)_{yz}
&=
\begin{cases}
0 &\mbox{if } \partial(x,y) \neq i-1;\\
|\Gamma_{i-1}(x) \cap \Gamma(y) \cap \Gamma(z)| &\mbox{if } \partial(x,y) = i-1.
\end{cases}\\
\end{align*}
We split our argument into cases.\\
Case $y = z$, or $\partial(y,z) > 2$, or $\partial(x,y) \neq i-1$:  The $(y,z)$-entries of $LF, FL$ and $L$ are zero.
Therefore the $(y,z)$-entry of each side of the equation in (i) is zero.\\
Case $\partial(y,z) = 2$ and $\partial(x,y) = i-1$:
By Lemma \ref{lem:dpg_yz2} the $(y,z)$-entry of $LF$ is $q^2$ and the $(y,z)$-entry of $FL$ is $1$.
The $(y,z)$-entry of $L$ is zero.
Therefore the $(y,z)$-entry of each side of the equation in (i) is zero.\\
Case $\partial(y,z) = 1$ and $\partial(x,y) = i-1$:
By Lemma \ref{lem:np_yz1} and Corollary \ref{cor:dpg_ai} the $(y,z)$-entry of $LF$ is $q^{2e}-1$ and the $(y,z)$-entry of $FL$ is $0$.
The $(y,z)$-entry of $L$ is 1.
Therefore the $(y,z)$-entry of each side of the equation in (i) is $q^{2e}-1$.\\
(ii) Take the transpose of each term in (i).
\hfill $\Box$\\

\begin{proposition}
\label{prop:rfl_reln2}
The following {\rm (i), (ii)} hold.
\begin{enumerate}
\item[\rm (i)] $\displaystyle \frac{q^4}{q^2+1} RL^2 - LRL + \frac{q^{-2}}{q^2+1} L^2 R = - q^{2e+2D-2} L$.
\item[\rm (ii)] $\displaystyle \frac{q^4}{q^2+1} R^2 L - RLR + \frac{q^{-2}}{q^2+1} L R^2 = - q^{2e+2D-2} R$.
\end{enumerate}
\end{proposition}

\noindent {\it Proof:} \rm
(i) Let
\begin{eqnarray*}
S = \frac{q^4}{q^2+1} RL^2 - LRL + \frac{q^{-2}}{q^2+1} L^2 R + q^{2e+2D-2} L.
\end{eqnarray*}
We show $S = 0$.
Since $I = \sum_{i=0}^D E^*_i$,
\begin{eqnarray*}
S = \left( \sum_{i=0}^D E^*_i \right) S \left( \sum_{j=0}^D E^*_j \right) = \sum_{i=0}^D \sum_{j=0}^D E^*_i S E^*_j.
\end{eqnarray*}
By Lemma \ref{lem:rfl_temp2} and since $E^*_lE^*_m = \delta_{lm}E^*_l$ for $0 \leq l,m \leq D$, we have $E^*_i S E^*_j = 0$ if $i \neq j-1$.
Therefore it suffices to show $E^*_{j-1} S E^*_j = 0$ for $1 \leq j \leq D$.
Let $j$ be given.
By Lemma \ref{lem:dpgdualeval},
\begin{eqnarray}
\theta^*_l - \theta^*_{l-1} = q^2 (\theta^*_{l+1} - \theta^*_{l}) \qquad (1 \leq l \leq D-1). \label{eq:ratiothetastar}
\end{eqnarray}
Expanding (\ref{eq:drg-td-reln1}) we get
\begin{align*}
0 &= A^3 A^* - (\beta + 1) A^2 A^* A + (\beta + 1) A A^* A^2 - A^* A^3\\
&\quad + \gamma (A^* A^2 - A^2 A^*) + \varrho (A^* A - A A^*).
\end{align*}
In the above equation multiply each term on the left by $E^*_{j-1}$ and on the right by $E^*_j$.
Simplify the result using $A = L + F + R$, and $A^*E_l^* = \theta_l^* E_l^*$ $(0 \leq l \leq D$) along with Lemma \ref{lem:dpgconsts}(i), Lemma \ref{lem:rfl_temp2} and (\ref{eq:ratiothetastar}).  This yields
\begin{align}
\label{eq:ese_temp}
\begin{split}
0
&= (q+q^{-1})^2E^*_{j-1} \left( -\frac{q^4}{q^2+1} RL^2 + LRL - \frac{q^{-2}}{q^2+1} L^2 R \right) E_j^*\\
&\quad + E^*_{j-1} \left( \vphantom{\frac{1}{1}} LF^2 - \beta FLF + F^2L - \gamma (FL + LF) - \varrho L \right) E_j^*.
\end{split}
\end{align}
Using Lemma \ref{lem:dpgconsts} and Proposition \ref{prop:rfl_reln}(i) one checks 
\begin{eqnarray}
LF^2 - \beta FLF + F^2L - \gamma (FL + LF) - \varrho L = - q^{2D+2e-2} (q+q^{-1})^2 L.\label{eq:ese_temp2}
\end{eqnarray}
Simplifying (\ref{eq:ese_temp}) using (\ref{eq:ese_temp2}) and $q+q^{-1} \neq 0$, we have $E^*_{j-1} S E^*_j = 0$.
We have now shown that $S = 0$.
The result follows.\\
(ii) Take the transpose of each term in (i).
\hfill $\Box$\\

\begin{definition} \rm
\label{def:thek}
We define the matrix $K \in T$ by%\in {\rm Mat}_X(\C)$ by
\begin{eqnarray}
K = \sum_{i=0}^D q^{-2i}E_i^*.\label{eq:thek}
\end{eqnarray}
\end{definition}

\noindent
Observe that $K$ is diagonal.  Moreover $\displaystyle K = \xi^{-1}(A^*-\zeta I)$ where $\zeta,\xi$ are from Lemma \ref{lem:dpgdualeval}.\\

\noindent
The next two lemmas follow from Definition \ref{def:thek}.

\begin{lemma}
\label{lem:thek1}
The matrix $K$ is invertible and $\displaystyle K^{-1} = \sum_{i=0}^D q^{2i} E_i^*$.
\end{lemma}

\begin{lemma}
\label{lem:thek2}
For $0 \leq i \leq D$ we have $KE_i^* = E_i^*K = q^{-2i}E_i^*$.
\end{lemma}

\begin{lemma}
\label{lem:rflk-gen-t}
The algebra $T$ is generated by $L,F,R,K$.
\end{lemma}

\noindent {\it Proof:}
Recall that $T$ is generated by $A, A^*$.
The result follows from this, the comments after Definition \ref{def:rfl}, and the comments after Definition \ref{def:thek}.
\hfill $\Box$

\begin{proposition}
\label{prop:krfl_reln}
The following {\rm (i)--(iii)} hold.
\begin{enumerate}
\item[\rm (i)] $KL = q^2 LK$.
\item[\rm (ii)] $KF = FK$.
\item[\rm (iii)] $KR = q^{-2}RK$.
\end{enumerate}
\end{proposition}

\noindent {\it Proof:} \rm
(i) By (\ref{eq:thek}) we have
$KL = \sum_{i=0}^D q^{-2i} E^*_iL$
and
$q^2 LK = \sum_{i=0}^d q^{2-2i} LE^*_i$.
Comparing these equations using Lemma \ref{lem:rfl_temp2}(i) we find $KL = q^2 LK$.\\
(ii) Similar to the proof of (i).\\
(iii) Take the transpose of each term in (i).
\hfill $\Box$\\

\noindent
We mention a consequence of the relations from Proposition \ref{prop:rfl_reln}, Proposition \ref{prop:rfl_reln2} and Proposition \ref{prop:krfl_reln}.

\begin{lemma}
\label{lem:RLFKcommute}
The matrices $LR, RL, F, K$ mutually commute.
\end{lemma}

\noindent {\it Proof:} \rm
By Proposition \ref{prop:krfl_reln} the matrix $K$ commutes with each of $LR, RL, F$.
By Proposition \ref{prop:rfl_reln} the matrix $F$ commutes with each of $LR, RL$.
It remains to show that $LR, RL$ commute.
In Proposition \ref{prop:rfl_reln2} multiply each side of equation (i) on the right by $R$ and multiply each side of equation (ii) on the left by $L$.
Taking the difference between the resulting equations and simplifying we get $LR^2L = RL^2R$.
Therefore $LR, RL$ commute.
The result follows.
\hfill $\Box$\\

\noindent
The matrices $L,F,R$ can be expressed in terms of $A,K,K^{-1}$ as follows.

\begin{lemma}
The following {\rm (i)--(iii)} hold.
\begin{enumerate}
\item[\rm (i)] $\displaystyle L = \frac{q^{-1}K^{-1}AK + qKAK^{-1} - (q+q^{-1})A}{(q-q^{-1})^2(q+q^{-1})}$.
\item[\rm (ii)] $\displaystyle F = \frac{(q^2+q^{-2})A - K^{-1}AK - KAK^{-1}}{(q-q^{-1})^2}$.
\item[\rm (iii)] $\displaystyle R = \frac{qK^{-1}AK + q^{-1}KAK^{-1} - (q+q^{-1})A}{(q-q^{-1})^2(q+q^{-1})}$.
\end{enumerate}
\end{lemma}

\noindent {\it Proof:}
In each equation eliminate $A$ using $A = L + F + R$ and simplify the result using Proposition \ref{prop:krfl_reln}.
\hfill $\Box$\\

\section{The central elements $\Omega,G,G^*$ in terms of $L,F,R,K$}
We continue to discuss the dual polar graph $\Gamma$ from Section \ref{sec:dpg}.
Recall  the central elements $\Omega, G, G^*$ of $T$ from Lemma \ref{lem:omega-z-zstar}.
Recall the elements $L,F,R,K$ of $T$ from Section \ref{sec:lfrk}.
In this section we display each of $\Omega, G, G^*$ in terms of $L,F,R,K$.

\begin{proposition}
\label{prop:ogg-to-lfrk}
Let $\xi, \zeta$ denote the scalars from Lemma \ref{lem:dpgdualeval},
and let $\gamma,\varrho$ denote the scalars from Lemma \ref{lem:dpgconsts}.
Then $\Omega,G,G^*$ can be expressed in terms of $L,F,R,K$ as follows.
\begin{enumerate}
\item[\rm (i)] $\displaystyle \Omega = -\frac{\xi(q^2-1)^2}{q^2} KF 
-\frac{\xi(q^2-1)(q^{2e}-1)}{q^2} K - 2\gamma\zeta I$.
\item[\rm (ii)] $\displaystyle G = \xi (1-q^4)KRL + \xi (1-q^{-4})KLR + \xi (q^{-2}-1)(q^{2e}-1)KF - \varrho \xi K - \varrho \zeta I$.
\item[\rm (iii)] $\displaystyle G^* = \frac{\zeta\xi(q^2-1)^2}{q^2} KF 
+\frac{\zeta\xi(q^2-1)(q^{2e}-1)}{q^2} K + \gamma\zeta^2 I$. 
\end{enumerate}
\end{proposition}

\noindent {\it Proof:} \rm
(i) Combine Lemma \ref{lem:Omega}, Definition \ref{def:rfl} and Definition \ref{def:thek}.\\
(ii) Combine Lemma \ref{lem:Z}, Definition \ref{def:rfl} and Definition \ref{def:thek}.\\
(iii) Combine Lemma \ref{lem:zstar} and (i).
\hfill $\Box$\\

\section{The central elements $C_0,C_1,C_2$}\label{sec:c012}

We continue to discuss the dual polar graph $\Gamma$ from Section \ref{sec:dpg}.
In this section we define three matrices $C_0, C_1, C_2$ in $T$ which involve $L,F,R,K$ from Section \ref{sec:lfrk}.
We show that $C_0, C_1, C_2$ are in $Z(T)$.
Then we display the actions of $C_0, C_1, C_2$ on each irreducible $T$-module.
Using this data we show that $C_0, C_1, C_2$ generate $Z(T)$.

\begin{definition} \rm
\label{def:misc_central_elements}
We define the matrices $C_0,C_1,C_2$ in $T$ as follows.
\begin{enumerate}
\item[\rm (i)] $\displaystyle C_0 = KF + \frac{q^{2e}-1}{q^2-1}K$.
\item[\rm (ii)] $\displaystyle C_1 = \frac{q^{2}}{q^2+1} KRL
- \frac{q^{-2}}{q^2+1} KLR
+ \frac{q^{2e+2D-2}}{q^2-1} K$.
\item[\rm (iii)] $\displaystyle C_2 = \frac{1}{q^2+1} K^2RL
-\frac{q^{-2}}{q^2+1} K^2LR
+ \frac{q^{2e+2D-2}}{q^4-1} K^2$.
\end{enumerate}
\end{definition}

\noindent
We have an observation.

\begin{lemma}
\label{lem:misc_central_elements_symmetric}
The matrices $C_0, C_1, C_2$ are symmetric.
\end{lemma}

\noindent {\it Proof:}
By Lemma \ref{lem:RLFKcommute} and Definition \ref{def:misc_central_elements} together with the fact that $F, K$ are symmetric and $R^t = L$.
\hfill $\Box$\\

\begin{lemma}
\label{lem:central_elements_c}
Each of $C_0, C_1, C_2$ is in $Z(T)$.
\end{lemma}

\noindent {\it Proof:} \rm
We first show that $C_0 \in Z(T)$.
By Lemma \ref{lem:rflk-gen-t} it suffices to show that $C_0$ commutes with each of $L,F,R,K$.
By Lemma \ref{lem:RLFKcommute} the matrix $C_0$ commutes with $F, K$.
Using Proposition \ref{prop:rfl_reln}(ii) and Proposition \ref{prop:krfl_reln}(iii)
one checks that $C_0$ commutes with $R$.
By this, Lemma \ref{lem:misc_central_elements_symmetric} and since $R^t = L$, we have
$C_0L - LC_0 = (RC_0 - C_0R)^t = 0.$
Therefore $C_0$ commutes with $L$.
We have now shown that $C_0$ commutes with each of $L,F,R,K$.
Therefore $C_0 \in Z(T)$.
By a similar argument using Proposition \ref{prop:rfl_reln2}, Proposition \ref{prop:krfl_reln} and Lemma \ref{lem:RLFKcommute} we find that $C_1, C_2$ are in $Z(T)$.
\hfill $\Box$\\

\begin{lemma}
\label{lem:central_elements_scalars}
Let $W$ denote an irreducible $T$-module.
Let $r, t, d$ denote the endpoint, dual endpoint and diameter of $W$, respectively.
Then the following {\rm (i)--(iii)} hold.
\begin{enumerate}
\item[\rm (i)] $C_0$ acts on $W$ as $\chi_0(r,t,d)I$ where
\begin{eqnarray}
\chi_0(r,t,d) = \frac{q^{2e+2D-2d-2r-2t} - q^{2t-2r}}{q^2-1}.\label{eq:chi0scalar}
\end{eqnarray}
\item[\rm (ii)] $C_1$ acts on $W$ as $\chi_1(r,t,d)I$ where
\begin{eqnarray}
\chi_1(r,t,d) = \frac{q^{2e+2D-1-d-2r}(q^{d+1} + q^{-d-1})}{q^4-1}.\label{eq:chi1scalar}
\end{eqnarray}
\item[\rm (iii)] $C_2$ acts on $W$ as $\chi_2(r,t,d)I$ where
\begin{eqnarray}
\chi_2(r,t,d) = \frac{q^{2e+2D-2-2d-4r}}{q^4-1}.\label{eq:chi2scalar}
\end{eqnarray}
\end{enumerate}
\end{lemma}

\noindent {\it Proof:} \rm
Fix an integer $i$ $(0 \leq i \leq d)$.
By Definition \ref{def:misc_central_elements}(i), (\ref{eq:thef}), (\ref{eq:thek}) and Definition \ref{def:stand-basis},
the element $C_0$ acts on $E_{r+i}^*W$ as
\begin{eqnarray}
q^{-2r-2i} a_i(W) + \frac{q^{-2r-2i}(q^{2e}-1)}{q^2-1}\label{eq:c0scalar}
\end{eqnarray}
times $I$.
Using Lemma \ref{lem:dpgintw} one checks that (\ref{eq:c0scalar}) equals $\chi_0(r,t,d)$.
Therefore (i) holds.
By Definition \ref{def:misc_central_elements}(ii), (\ref{eq:thel}), (\ref{eq:ther}), (\ref{eq:thek}) and Definition \ref{def:stand-basis},
the element $C_1$ acts on $E_{r+i}^*W$ as
\begin{eqnarray}
\frac{q^{2-2r-2i}}{q^2+1} c_{i}(W)b_{i-1}(W) -\frac{q^{-2-2r-2i}}{q^2+1} b_{i}(W)c_{i+1}(W) + \frac{q^{2e+2D-2-2r-2i}}{q^2-1}\label{eq:c1scalar}
\end{eqnarray}
times $I$.
Using Lemma \ref{lem:dpgintw} one checks that (\ref{eq:c1scalar}) equals $\chi_1(r,t,d)$.
Therefore (ii) holds.
By Definition \ref{def:misc_central_elements}(iii), (\ref{eq:thel}), (\ref{eq:ther}), (\ref{eq:thek}) and Definition \ref{def:stand-basis},
the element $C_2$ acts on $E_{r+i}^*W$ as
\begin{eqnarray}
\frac{q^{-4r-4i}}{q^2+1} c_{i}(W)b_{i-1}(W) -\frac{q^{-2-4r-4i}}{q^2+1} b_{i}(W)c_{i+1}(W) + \frac{q^{2e+2D-2-4r-4i}}{q^4-1}\label{eq:c2scalar}
\end{eqnarray}
times $I$.
Using Lemma \ref{lem:dpgintw} one checks that (\ref{eq:c2scalar}) equals $\chi_2(r,t,d)$.
Therefore (iii) holds.
\hfill $\Box$\\

\begin{theorem}\label{thm:center-c012}
The algebra $Z(T)$ is generated by $C_0,C_1,C_2$.
\end{theorem}

\noindent {\it Proof:} \rm
Let $Z'$ denote the subalgebra of $T$ generated by $C_0,C_1,C_2$.
We show $Z' = Z(T)$.
By Lemma \ref{lem:central_elements_c}, the algebra $Z'$ is contained in $Z(T)$.
We now show the reverse inclusion.
To do this, by Lemma \ref{lem:centerofT} it suffices to show $E_\lambda \in Z'$ for all $\lambda \in {\rm Feas}$.
Let $(r,t,d)$ and $(r',t',d')$ denote distinct elements of ${\rm Feas}$.
We claim that there exists an integer $i \in \{0,1,2\}$ such that $\chi_i(r,t,d) \neq \chi_i(r',t',d')$.
Suppose this is not the case.
By (\ref{eq:chi0scalar})--(\ref{eq:chi2scalar}),
\begin{align}
q^{2e+2D-2d-2r-2t}-q^{2t-2r} &= q^{2e+2D-2d'-2r'-2t'}-q^{2t'-2r'}, \label{eq:c0prime}\\
q^{-d-2r}(q^{d+1}+q^{-d-1}) &= q^{-d'-2r'}(q^{d'+1}+q^{-d'-1}), \label{eq:c1prime}\\
q^{-2d-4r} &= q^{-2d'-4r'}. \label{eq:c2prime}
\end{align}
We show $d=d'$.
By (\ref{eq:c2prime}) we have $q^{-d-2r} = q^{-d'-2r'}$.
By this and (\ref{eq:c1prime}) we have $q^{d+1}+q^{-d-1} = q^{d'+1}+q^{-d'-1}$.
Simplifying this we get $(q^{d+d'+2}-1)(q^{d+1}-q^{d'+1}) = 0$.
Since $d+d'+2 \neq 0$ and $q$ is not a root of unity, we have $q^{d+1} = q^{d'+1}$ so $d = d'$.
Next we show $r=r'$.
This is immediate from (\ref{eq:c2prime}) and the fact that $d = d'$.
Now we show $t=t'$.
Evaluate  (\ref{eq:c0prime}) using $r=r'$ and $d=d'$, and simplify the result 
we get $(q^{2e+2D-2d}+q^{2t+2t'})(q^{2t'}-q^{2t})=0$.
But $q^{2e+2D-2d}+q^{2t+2t'} \neq 0$ since $q$ is not a root of unity.
Therfore $q^{2t'} = q^{2t}$ and thus $t=t'$.
We have shown $(r,t,d) = (r',t',d')$ for a contradiction.
Therefore the claim holds.
By the claim and Lemma \ref{lem:central_elements_scalars} we have $E_\lambda \in Z'$ for all $\lambda \in {\rm Feas}$.
The result follows.
\hfill $\Box$\\

\noindent
We finish this section with a comment.

\begin{lemma}
\label{lem:rl_c1c2}
The following {\rm (i), (ii)} hold.
\begin{enumerate}
\item[\rm (i)] $\displaystyle RL = \frac{q^2+1}{q^2-1}K^{-1}C_1 - \frac{q^2+1}{q^2-1}K^{-2}C_2 - \frac{q^{2e+2D}}{(q^2-1)^2}I$.
\item[\rm (ii)] $\displaystyle LR = \frac{q^2(q^2+1)}{q^2-1}K^{-1}C_1 - \frac{q^4(q^2+1)}{q^2-1}K^{-2}C_2 - \frac{q^{2e+2D}}{(q^2-1)^2}I$.
\end{enumerate}
\end{lemma}

\noindent {\it Proof:} \rm
Solve for $RL$ and $LR$ using Definition \ref{def:misc_central_elements}(ii), (iii).
\hfill $\Box$\\

\section{How $C_0,C_1,C_2$ relate to $\Omega,G,G^*$ and $\Upsilon,\Psi,\Lambda$}

We continue to discuss the dual polar graph $\Gamma$ from Section \ref{sec:dpg}.
So far we obtained a number of central elements of $T$.
We have $\Omega, G, G^*$ from Lemma \ref{lem:omega-z-zstar}, $\Upsilon, \Psi, \Lambda$ from Section \ref{sec:up-psi-lamb}, and $C_0, C_1, C_2$ from Section \ref{sec:c012}.
In Proposition \ref{prop:omegazzstar-upsilonpsilambda} we expressed $\Omega, G, G^*$ in terms of $\Upsilon, \Psi, \Lambda$.
In this section we express $\Omega, G, G^*$ in terms of $C_0, C_1, C_2$
and express $C_0, C_1, C_2$ in terms of $\Upsilon, \Psi, \Lambda$.

\begin{proposition}
\label{prop:OGGandCs}
Let $\zeta, \xi$ denote the scalars from Lemma \ref{lem:dpgdualeval} and let $\gamma,\varrho$ denote the scalars from Lemma \ref{lem:dpgconsts}.
Then the following {\rm (i)}--{\rm (iii)} hold.
\begin{enumerate}
\item[\rm (i)] $\Omega = - \xi q^{-2} (q^2-1)^2 C_0 -2\gamma \zeta I$.
\item[\rm (ii)] $G = \xi (q^{-2}-1) (q^{2e}-1)C_0 + \xi (q^{-2}-1) (q^2+1)^2 C_1 - \varrho \zeta I$.
\item[\rm (iii)] $G^* = \zeta \xi q^{-2} (q^2-1)^2 C_0 + \gamma \zeta^2 I$.
\end{enumerate}
\end{proposition}

\noindent {\it Proof:} \rm
(i) Evaluate $\Omega$ using Proposition \ref{prop:ogg-to-lfrk}(i).  Evaluate $C_0$ using Definition \ref{def:misc_central_elements}.\\
(ii) Evaluate $G$ using Proposition \ref{prop:ogg-to-lfrk}(ii).  Evaluate $C_0, C_1$ using Definition \ref{def:misc_central_elements}.\\
(iii) Combine Lemma \ref{lem:zstar} and (i).
\hfill $\Box$\\

\begin{proposition}
\label{prop:CsandPhiPsi}
The matrices $C_0, C_1, C_2$ can be expressed in terms of $\Upsilon, \Psi,\Lambda$ as follows.
\begin{enumerate}
\item[\rm (i)] $C_0 = (q^2-1)^{-1} (q^{2e} \Upsilon^{-2} - \Psi^{-2})$.
\item[\rm (ii)] $C_1 = q^{2e+D-3} (q^2-1) (q^2+1)^{-1} \Upsilon^{-1} \Psi^{-1} \Lambda$.
\item[\rm (iii)] $C_2 = q^{2e-2}(q^4-1)^{-1} \Upsilon^{-2} \Psi^{-2}$.
\end{enumerate}
\end{proposition}

\noindent {\it Proof:} \rm
(i) By construction $V$ is a direct sum of irreducible $T$-modules.
Let $W$ denote an irreducible $T$-module in the sum.
It suffices to show that for the equation in (i) the two sides agree on $W$.
Let $r, t, d$ denote the endpoint, dual endpoint, and diameter of $W$, respectively.
By Lemma \ref{lem:central_elements_scalars}, the element $C_0$ acts on $W$ as $\chi_0(r,t,d)I$.
By Lemma \ref{lem:displacement1-decomp} we have $W \subseteq \sigma_{r+t+d-D}V$.
By this and Lemma \ref{lem:Upsilon-action} the matrix $\Upsilon$ acts on $W$ as $q^{r+t+d-D}1$.
By Lemma \ref{lem:displacement2-decomp} we have $W \subseteq \psi_{r-t}V$.
By this and Lemma \ref{lem:Psi-action} the matrix $\Psi$ acts on $W$ as $q^{r-t}1$.
By these comments and (\ref{eq:chi0scalar}), the two sides agree on $W$.
The result follows\\
(ii), (iii) Similar to the proof of (i).
\hfill $\Box$\\

\section{$L,F,R,K$ and $U_q(\mathfrak{sl}_2)$}

We continue to discuss the dual polar graph $\Gamma$ from Section \ref{sec:dpg}.
In this section we display some relationships between the $U_q(\mathfrak{sl}_2)$-module structures from Theorem \ref{thm:uqsl2ondpg} and Theorem \ref{thm:uqsl2ondpg_2} and the matrices $L,F,R,K$ from Section \ref{sec:lfrk}.

\medskip
\noindent
Recall the central elements $\Upsilon, \Psi$ of $T$ from Section \ref{sec:up-psi-lamb}.

\begin{lemma}
\label{lem:xyzondpg}
The actions of $x,y,z$ on the $U_q(\mathfrak{sl}_2)$-module $V$ from Theorem \ref{thm:uqsl2ondpg} are given by
\begin{align*}
x & = q^{-D} \Upsilon^{-1} \Psi^{-1} K^{-1} + q^{-2e-D}(q^2-1) \Upsilon \Psi^{-1} R,\\
y & = q^{-D} \Upsilon^{-1} \Psi^{-1} K^{-1} - q^{-D}(q^2-1) \Upsilon^{-1} \Psi L,\\
z & = q^D \Upsilon \Psi K.
\end{align*}
\end{lemma}

\noindent {\it Proof:} \rm
Let $h$ denote the scalar from (\ref{eq:dpg-params}) and let $\xi$ denote the scalar from Lemma \ref{lem:dpgdualeval}.
By Lemma \ref{def:misc_central_elements} we have $F - hI = K^{-1}C_0$.
By this and Proposition \ref{prop:CsandPhiPsi} we have
\begin{eqnarray}
F - hI =  (q^2-1)^{-1} (q^{2e} \Upsilon^{-2} - \Psi^{-2}) K^{-1}.\label{eq:f_up_psi}
\end{eqnarray}
In each of (\ref{eq:x1_bbstar})--(\ref{eq:z_bbstar}) evaluate $B,B^*$ using $B = L + F + R - hI$ and $B^* = \xi K$.  Simplify the result using Proposition \ref{prop:krfl_reln} and (\ref{eq:f_up_psi}).
\hfill $\Box$\\

\begin{lemma}
\label{lem:xyzprimeondpg}
The actions of $x,y,z$ on the $U_q(\mathfrak{sl}_2)$-module $V$ from Theorem \ref{thm:uqsl2ondpg_2} are given by
\begin{align*}
x & = q^{-D} \Upsilon^{-1} \Psi^{-1} K^{-1} - q^{-D}(q^2-1) \Upsilon^{-1} \Psi R,\\
y & = q^{-D} \Upsilon^{-1} \Psi^{-1} K^{-1} + q^{-2e-D}(q^2-1) \Upsilon \Psi^{-1} L,\\
z & = q^D \Upsilon \Psi K.
\end{align*}
\end{lemma}

\noindent {\it Proof:} \rm
Similar to the proof of Lemma \ref{lem:xyzondpg} but use (\ref{eq:xprime1_bbstar})--(\ref{eq:zprime_bbstar}) instead of (\ref{eq:x1_bbstar})--(\ref{eq:z_bbstar}).
\hfil $\Box$\\

\noindent
We now give reformulations of Lemma \ref{lem:xyzondpg} and Lemma \ref{lem:xyzprimeondpg} in terms of the generators $k, e, f$ for $U_q(\mathfrak{sl}_2)$.

\begin{lemma}
\label{lem:kefondpg}
The actions of $k,e,f$ on the $U_q(\mathfrak{sl}_2)$-module $V$ from Theorem \ref{thm:uqsl2ondpg} are given by
\begin{align*}
k &= q^D \Upsilon \Psi K,\\
e &= \Psi^2 K L,\\
f &= q^{1-2e-D} \Upsilon \Psi^{-1} R.
\end{align*}
\end{lemma}

\noindent {\it Proof:} \rm
Combine Lemma \ref{lem:uqsl2equitable} and Lemma \ref{lem:xyzondpg}.
\hfill $\Box$\\

\begin{lemma}
\label{lem:kefprimeondpg}
The actions of $k,e,f$ on the $U_q(\mathfrak{sl}_2)$-module $V$ from Theorem \ref{thm:uqsl2ondpg_2} are given by
\begin{align*}
k &= q^D \Upsilon \Psi K,\\
e &= -q^{-2e} \Upsilon^2 K L,\\
f &= -q^{1-D} \Upsilon^{-1} \Psi R.
\end{align*}
\end{lemma}

\noindent {\it Proof:} \rm
Combine Lemma \ref{lem:uqsl2equitable} and Lemma \ref{lem:xyzprimeondpg}.
\hfill $\Box$\\

\noindent
We finish this section with a comment decscribing how the two actions of $U_q(\mathfrak{sl}_2)$ on $V$ from Theorem \ref{thm:uqsl2ondpg}
and Theorem \ref{thm:uqsl2ondpg_2}
are related.

\begin{theorem}
Consider the table below.
In the first column the three displayed elements each induces an element in ${\rm End}(V)$ using the $U_q(\mathfrak{sl}_2)$-module structure from Theorem \ref{thm:uqsl2ondpg_2}.
In the second column the three displayed elements each induces an element in ${\rm End}(V)$ using the $U_q(\mathfrak{sl}_2)$-module structure from Theorem \ref{thm:uqsl2ondpg}.
For each row the two elements induce the same element of ${\rm End}(V)$.\\
\begin{center}
\begin{tabular}{c|c}
$U_q(\mathfrak{sl}_2)$-module structure & $U_q(\mathfrak{sl}_2)$-module structure\\
from Theorem \ref{thm:uqsl2ondpg_2} & from Theorem \ref{thm:uqsl2ondpg}\\
\hline
$k$ & $k$ \\
$e$ & $-q^{-2e} \Upsilon^{2} \Psi^{-2} e$ \\
$f$ & $-q^{2e} \Upsilon^{-2} \Psi^{2} f$\\
\end{tabular}
\end{center}
\end{theorem}

\noindent {\it Proof:} \rm
Compare the $U_q(\mathfrak{sl}_2)$-actions from Lemma \ref{lem:kefondpg} and Lemma \ref{lem:kefprimeondpg}.
\hfill $\Box$\\

\section{Acknowledgements}
This paper was the author's Ph.D. thesis at the University of Wisconsin--Madison.  The author would like to thank his advisor Paul Terwilliger for offering many valuable ideas and suggestions.

\noindent
Chalermpong Worawannotai\\
Department of Mathematics\\
University of Wisconsin\\
480 Lincoln Drive\\
Madison, WI 53706-1388 USA\\
email: \texttt{worawann@math.wisc.edu}


\begin{thebibliography}{99}


\bibitem{artin} E. Artin. {\it Geometric Algebra.} Interscience, New York, 1957.

\bibitem{bannai} E. Bannai and T. Ito. {\it Algebraic Combinatorics I:
      Association Schemes.} Benjamin/Cummings, London, 1984.

\bibitem{biggs} N. Biggs.
{\it Algebraic Graph Theory. Second edition.} Cambridge University Press, Cambridge, 1993.

\bibitem{bcn} A. E. Brouwer, A. M. Cohen, and A. Neumaier. 
      {\it Distance-Regular Graphs.} Springer-Verlag, Berlin, 1989.

\bibitem{quads} A. E. Brouwer and H. A. Wilbrink.
\newblock The structure of near polygons with quads.
\newblock{\em Geommetriae Dedicata.} {\bf14} (1983), 145--176.


\bibitem{cameron} P. J. Cameron.
{\it Projective and Polar Spaces.} QWM Math Notes, London, 1992.


\bibitem{cerzo} D. Cerzo.
\newblock Structure of thin irreducible modules of a $Q$-polynomial distance-regular graph.
\newblock {\em Linear Algebra Appl.} {\bf 433} (2010), no. 8--10, 1573--1613.
{\tt arXiv:1003.5368}.


\bibitem{CR} C. Curtis and I. Reiner.
{\it {R}epresentation {T}heory of {F}inite {G}roups and
{A}ssociative {A}lgebras}.
Interscience, New York, 1962.


%\bibitem{gasper} G. Gasper and M. Rahman.
%{\it Basic Hypergeometric Series,} volume 35 of {\it Encyclopedia of Mathematics and its Applications.} Cambridge University Press, Cambridge, 1990. 


\bibitem{go2} J. T. Go and P. Terwilliger.
 Tight distance-regular graphs and the subconstituent algebra.
\newblock{\em
 European J. Combin.} {\bf 23} (2002),
   793--816. 

\bibitem{godsil}
C. D. Godsil.
{\it Algebraic Combinatorics.} Chapman and Hall, Inc., New York, 1993.

\bibitem{horn}
R.A. Horn and C.R. Johnson.
{\it Matrix Analysis.} Cambridge University Press, Cambridge, MA, 1985.


\bibitem{somealg}
T. Ito, K. Tanabe and P. Terwilliger.
Some algebra related to $P$- and $Q$-polynomial
association schemes.
{\em Codes and association schemes (Piscataway NJ, 1999)},
167--192, DIMACS Ser. Discrete Math. Theoret. Comput. Sci. 56,
{\em Amer. Math. Soc., Providence RI,} 2001.


\bibitem{drgqtet}
T. Ito and P. Terwilliger.
\newblock Distance-regular graphs and the $q$-tetrahedron algebra.
{\em European J. Combin.}
{\bf 30} (2009), no. 3, 682--697.
{\tt arXiv:math.CO/0608694.}


\bibitem{q-inv}
T. Ito and P. Terwilliger.
\newblock $q$-Inverting pairs of linear transformations and the $q$-tetrahedron algebra.
\newblock {\em Linear Algebra Appl.} {\bf 426} (2007), no. 2--3, 516--532. 
{\tt arXiv:math.RT/0606237.}


\bibitem{equit1}
T. Ito, P. Terwilliger, and C. W. Weng.
\newblock The quantum algebra 
$U_q(\mathfrak{sl}_2)$ and its equitable presentation.
\newblock {\em J. Algebra} {\bf 298} (2006), 284--301.
{\tt arXiv:math.QA/0507477}. 


\bibitem{jantzen}
J. C. Jantzen.
{\it Lectures on quantum groups,} Graduate Studies in Mathematics 6, Amer. Math. Soc., Providence RI, 1996.


\bibitem{kassel}
C. Kassel.
{\it Quantum Groups,} Springer-Verlag, New York, 1995.


\bibitem{kim}
Joohyung Kim.
\newblock Some matrices associated with the split decomposition for a $Q$-polynomial distance-regular graph.
\newblock {\em European J. Combin.} {\bf 30} (2009), no. 1, 96--113.


%\bibitem{koekoek}
%R. Koekoek and R. Swarttouw.
%{\it The Askey-scheme of hypergeometric orthogonal polynomials and its q-analog,} volume 98-17 of {\it Reports of the faculty of Technical Mathematics and Informatics.}  Delft, The Netherlands, 1998.


\bibitem{aap1}
A. A. Pascasio. 
 On the multiplicities of the primitive idempotents of a $Q$-polynomial distance-regular graph.
\newblock{\em
European J.
   Combin.} {\bf 23} (2002),  1073--1078. 



\bibitem{shult}
E. E. Shult and A. Yanushka.
\newblock Near $n$-gons and line systems.
\newblock {\em Geometriae Dedicata.} {\bf 9} (1980), 1--72.


\bibitem{terw-lp24}
P. Terwilliger.
\newblock Leonard pairs from 24 points of view.
\newblock {\em Rocky mountain Journal of Mathematics.}, 32(2):827--888, 2002.


\bibitem{terwSub1} P. Terwilliger. 
The subconstituent algebra of
an association scheme I.
\newblock{\em
J. Algebraic Combin.}
{\bf 1} (1992), 363--388.  

\bibitem{terwSub2} P. Terwilliger. 
The subconstituent algebra of
an association scheme II.
\newblock{\em
J. Algebraic Combin.}
{\bf 2} (1993), 73--103.  


\bibitem{terwSub3} P. Terwilliger.  The subconstituent algebra of
an association scheme III. 
\newblock{ \em
J. Algebraic Combin.}
{\bf 2} (1993), 177--210.  


\bibitem{ds}
P. Terwilliger.
\newblock The displacement and split decompositions
for a $Q$-polynomial distance-regular graph.
\newblock{\em Graphs Combin.} {\bf 21} (2005), 263--276.
{\tt arXiv:math.CO/0306142}.


\bibitem{LS99}
P. Terwilliger.
\newblock Two linear transformations each tridiagonal with respect to an
  eigenbasis of the other.
  \newblock {\em Linear Algebra Appl.}  {\bf 330} (2001), 149--203.
{\tt arXiv:math.RA/0406555}.


\bibitem{madrid09}
P. Terwilliger.
\newblock Two linear transformations each tridiagonal with respect to an eigenbasis of the other; an algebraic approach to the Askey scheme of orthogonal polynomials. Lecture notes for the summer school on orthogonal polynomials and special
functions. Universidad Carlos III de Madrid, Leganes, Spain. July 8--July 18, 2004.
{\tt arXiv:math.QA/0408390}


\bibitem{tdd-lbub}
P. Terwilliger.
\newblock Two linear transformations each tridiagonal with respect to an eigenbasis of the other; the TD-D canonical form and the LB-UB canonical form.
\newblock {\em J. Algebra.} {\bf 291} (2005), no. 1, 1--45.
{\tt arXiv:math.RA/0304077}.


\bibitem{uaw-uqsl2}
P. Terwilliger.
\newblock The universal Askey-Wilson algebra and the equitable presentation of $U_q(\mathfrak{sl}_2)$.
\newblock {\em SIGMA} {\bf 7} (2011) 099, 26 pages.
{\tt arXiv:1107.3544}.


\bibitem{lpaw}
P. Terwilliger and R. Vidunas.
\newblock Leonard pairs and the Askey-Wilson relations.
\newblock{\em J. Alegbra Appl.}
{\bf 3} (2004) 411--426.
{\tt arXiv:math.QA/0305356}.


\bibitem{weakgeod} C. Weng.  Weak-geodetically closed subgraphs in distance-regular graphs.
\newblock{ \em Graphs. Combin.}
{\bf 14} (1998), no.3, 275--304.


\end{thebibliography}
\end{document}